\documentclass[preprint,11pt]{elsarticle}

\usepackage[margin=1in]{geometry}
\usepackage{amsmath,amssymb,bm}
\usepackage{amsthm}
\theoremstyle{plain}
\newtheorem{theorem}{Theorem}
\newtheorem{lemma}{Lemma}
\usepackage{mathtools}
\usepackage{booktabs}
\usepackage{longtable}
\usepackage{graphicx}
\usepackage{tikz}
\usetikzlibrary{arrows.meta}
\usepackage{placeins}
\usepackage{enumitem}
\usepackage{hyperref}
\usepackage{xcolor}
\usepackage{tabularx}
\usepackage{array}
\usepackage{microtype}
\usepackage{etoolbox}

\newcommand{\tablefont}{\fontsize{9}{11}\selectfont}
\AtBeginEnvironment{tabular}{\tablefont}
\AtBeginEnvironment{tabular*}{\tablefont}
\AtBeginEnvironment{tabularx}{\tablefont}
\AtBeginEnvironment{longtable}{\tablefont}

\newcommand{\R}{\mathbb{R}}
\newcommand{\bx}{x}
\newcommand{\bv}{v}

\begin{document}

\begin{frontmatter}
    
\title{Split-and-Learn: Decoupled Learning for Reactive Transport \\
in Heterogeneous Marine Sediments}

\author[tamucc-math]{Maria Vasilyeva}
\ead{maria.vasilyeva@tamucc.edu}
% \cortext[cor1]{Corresponding author.}
\author[llnl]{Siu Wun Cheung}
\ead{cheung26@llnl.gov}
\author[tamu]{Jihoon Kim}
\ead{jihoon.kim@tamu.edu}
\author[tamucc-pens]{Richard B. Coffin}
\ead{richard.coffin@tamucc.edu}
\author[tamucc-pens]{Ingo Pecher}
\ead{ingo.pecher@tamucc.edu}

\affiliation[tamucc-math]{organization={Department of Mathematics and Statistics, Texas A\&M University-Corpus Christi},
  city={Corpus Christi}, state={TX}, country={USA}}
\affiliation[llnl]{organization={Center for Applied Scientific Computing, Lawrence Livermore National Laboratory},
  city={Livermore}, state={CA}, country={USA}}
\affiliation[tamu]{organization={Department of Petroleum Engineering, Texas A\&M University},
  city={College Station}, state={TX}, country={USA}}
\affiliation[tamucc-pens]{organization={Department of Physical and Environmental Sciences, Texas A\&M University-Corpus Christi},
  city={Corpus Christi}, state={TX}, country={USA}}

\begin{abstract}
Reactive transport in heterogeneous marine sediments couples flow and multicomponent transport with nonlinear chemistry and evolving porosity and permeability, and its simulation over thousands of years is computationally expensive, especially when uncertain heterogeneity requires many simulations.
We develop Split-and-Learn, a framework that keeps the sequential operator-splitting scheme of the full-order model and replaces each stage by a data-driven surrogate with its own structure. 
Transport is learned on a local data-driven multiscale space, either as a compressed reduced solution map (global inverse atoms) or as coarse-block pieces of the forward operator assembled into a Galerkin or least-squares Petrov–Galerkin system (local atoms). 
Chemistry is learned through bounded reaction extents with stoichiometric reconstruction, which preserves reaction invariants and nonnegativity. 
Flow uses learned interface traces and a reduced mixed multiscale solve that enforces coarse-block mass balance. 
{We prove the stability of the split scheme with the evolving medium and bound the error of the learned solver by the accumulated one-step errors of its stages.}
Numerical experiments in two- and three-dimensional heterogeneous sediment domains show the accuracy and computational cost of the individual and combined learned stages against the full-order split solver. 
\end{abstract}


\end{frontmatter}

\section{Introduction}

Reactive transport in marine sediments is essential for predicting CO$_2$ residence time \cite{coffin2024offshore} and, more broadly, for evaluating whether injected or naturally occurring CO$_2$ remains stored, migrates through the sediment column, or participates in long-term biogeochemical cycling within geologically complex offshore formations \cite{akam2020dissolved,coffin2023importance,hommel2018poroperm}.
{Chemical reactions between rock and fluids during CO$_2$ storage induce mineral dissolution and precipitation, thereby altering porosity, permeability, and CO$_2$ injectivity \cite{XU2004917}.}
{In marine sediments, microbially mediated reactions also control the fate of the carbon, since the degradation of organic matter by sulfate reduction produces dissolved inorganic carbon and thereby drives carbonate precipitation \cite{akam2020dissolved,chatterjee2011pore}.}
The fate of stored CO$_2$ is strongly influenced by uncertain and spatially heterogeneous subsurface properties, which necessitates repeated coupled flow, transport, and reaction simulations across a range of geological realizations and parameter scenarios.
The computational cost of such simulations arises from repeated spatial flow and transport solves together with nonlinear chemical calculations over long time horizons.

Coupled reactive transport is commonly solved by sequential operator splitting, which separates flow, transport, and
local chemistry and treats each subproblem with a suited solver \cite{yeh1989critical,steefel2015codes,hundsdorfer2003numerical,carrayrou2004operatorsplitting,lagneau2010porosity}. 
{For example, employing a sequential approach, we can first solve the flow problem implicitly and then solve the solute transport problem along with chemical reaction implicitly as well \cite{Xu2008}. This sequential framework can also be extended to incorporate coupling with geomechanics \cite{Taron2009aTHMC,KIM201559}.}
Machine learning has been used in two ways to accelerate such simulations. Emulators, neural networks, and neural operators predict model responses or solution fields directly
\cite{laloy2019emulation,han2024surrogate,li2022artificial,badawi2025neural,lu2021deeponet,li2020fno}.
Used as time steppers, however, their errors accumulate over many steps and through the coupling, training remedies give no bounds \cite{brandstetter2022message}, and a single map for the coupled system \cite{tang2022coupled,wen2022ufno} loses the structure of each process. 
Process decomposition before learning is used in deep splitting and physical-operator learning \cite{beck2021deepsplitting,gopakumar2026physicaloperators}.
Learned chemistry inside a split solver \cite{leal2020odml,prasianakis2020nn,laloy2022surrogate,delucia2021dectree} keeps the numerical scheme, but once chemistry is cheap, the flow and transport solves dominate the cost.

The flow and transport solves can be reduced by multiscale methods \cite{hou1997multiscale,jenny2003msfv,efendiev2009multiscale,
efendiev2013gmsfem,chung2016mixed} and projection-based model reduction \cite{benner2015survey,hesthaven2016certified,
carlberg2011lspg,negri2015mdeim}. Both inherit the structure of the full model, but each new heterogeneous field needs
fine-scale work online for local basis computations, assembly of the fine-grid operator, and projection. Non-intrusive methods learn reduced dynamics from data \cite{schmid2010dmd,peherstorfer2016opinf,amsallem2011interpolation}, usually for
low-dimensional parameters rather than heterogeneous fields, and local neural compression of coarse operators has been analyzed for steady elliptic problems \cite{kropfl2022compression}. 
{Neural networks were also trained to predict the local coarse matrices and multiscale basis functions of a multiscale method from the local permeability \cite{wang2019gmsfemdl}, and to advance reduced multiscale models of porous media flow in time \cite{cheung2020deepglobal,wang2020reduceddl}.}
Domain-decomposition reduced-order models \cite{hoang2021ddlspg,diaz2024ddnmrom,chung2024trainsmall} and the data-driven finite element method \cite{choi2026rigor,cheung2026partition} assemble global reduced models from bases trained on local subdomains. 
A data-driven multiscale space with local snapshot bases was combined with operator learning for steady nonlinear flow in random heterogeneous media \cite{vasilyeva2026operator}, but a learned solution operator used as a time stepper can become unstable as its errors accumulate over many steps.

In this work we develop \textit{Split-and-Learn}, a framework for multicomponent reactive transport in heterogeneous marine sediment
columns. We keep the sequential scheme of \cite{vasilyeva2024decoupled}, including the porosity and permeability updates,
and replace each stage by a data-driven surrogate that keeps the structure of its process.
Each surrogate is trained and evaluated on the intermediate states of the split solver. 
The contributions are as follows.
\begin{itemize}[leftmargin=*,itemsep=1pt]
\item \textit{Two transport models on a local data-driven multiscale space}. Global inverse atoms learn the reduced solution map of the transport step and need only a matrix--vector product online. Local atoms learn coarse-block pieces of the forward operator and assemble a reduced system.
\item \textit{Structure-preserving chemistry and flow surrogates.} Chemistry is learned through bounded reaction extents with stoichiometric reconstruction, which preserves reaction invariants and nonnegativity. Flow uses learned interface traces and a reduced mixed multiscale {solver} with coarse-block mass balance.
\item {\textit{Stability and error analysis.} We prove the stability of the split scheme with the evolving medium, with a growth rate independent of the time step, and bound the error of the multiscale and learned solvers by the accumulated one-step errors of their stages.}
\item \textit{Numerical experiments.} In 2D and 3D, we compare individual and combined learned stages with direct and iterative full-order references and with the same multiscale spaces without learning.
The fully learned 3D solver is up to $57$ times faster than the direct reference, and the transport error of local atoms is three to six times smaller than that of global inverse atoms.
\end{itemize}

The paper is organized as follows. 
Section~\ref{sec:model} describes the model, and Section~\ref{sec:numerics} the operator splitting and finite-volume discretization. 
Section~\ref{sec:learning} develops the learned transport, chemistry, flow.
{Section~\ref{sec:learned-error-propagation} analyzes the stability of the split scheme and the propagation of learned-stage errors.} 
Section~\ref{sec:results} presents the numerical results. The paper ends with a conclusion.

\section{Reactive Transport in Random Heterogeneous Media}
\label{sec:model}

We consider multicomponent reactive transport in a heterogeneous porous medium.  
The motivating application is the evolution of geochemical profiles in a marine-sediment column \cite{akam2020dissolved, coffin2023importance, vasilyeva2024decoupled}.  
The model contains pressure-driven flow, species-wise advection--diffusion, and local nonlinear geochemical reactions, which are three mechanisms with different mathematical structure \cite{steefel2018approaches,regnier2002modeling}.

\subsection{Multicomponent transport}

Let $\Omega\subset\R^d$ ($d=2,3$) be the spatial domain and $t\in[0,T]$ the time.
The concentration $c_i(x,t)$ of component $i=1,\ldots,N_{\mathrm{comp}}$ satisfies the mass balance
\begin{equation*}
  \frac{\partial}{\partial t}\big(\eta_i c_i\big)
  +\nabla\cdot\big(\eta_i\bv_i c_i\big)
  -\nabla\cdot\big(\eta_iD_i(\phi)\nabla c_i\big)
  =\eta_i R_i(c_1,\ldots,c_{N_{\mathrm{comp}}}),
  \quad x \in \Omega, \quad t\in(0,T],
\end{equation*}
where the phase weight $\eta_i$ accounts for the volume occupied by the
component and $\bv_i$ is its phase velocity
\begin{equation}
  \eta_i(\phi)=
  \begin{cases}
    1-\phi, & i\in\mathcal I_s,\\
    \phi,   & i\in\mathcal I_f,
  \end{cases}
  \quad
  \bv_i=
  \begin{cases}
    \bv_s, & i\in\mathcal I_s,\\
    \bv_f, & i\in\mathcal I_f,
  \end{cases}
  \label{eq:phase-weights}
\end{equation}
where $\mathcal I_s$ and
$\mathcal I_f$ are the index sets for solid and dissolved components, $\phi\in(0,1)$ is the porosity, $\bv_s$ is the downward solid-burial velocity and  $\bv_f=q/\phi$ is the pore-fluid velocity{, and $D_i(\phi)$ is the effective diffusion coefficient, given below}.

{The pore fluid of density $\rho$ satisfies the mass balance
\begin{equation*}
 \partial_t(\rho\phi)+\nabla\cdot(\rho q)=0,\quad
 \frac{1}{\rho}\partial_t(\rho\phi)=\phi c_f\,\partial_tp+\partial_t\phi,
\end{equation*}
where $c_f=\rho^{-1}d\rho/dp$ is the fluid compressibility.
We assume an incompressible pore fluid, $c_f=0$, so that $\rho$ is constant.}
The Darcy flux $q$ and the pressure $p$ satisfy Darcy's law and the mass balance of an incompressible fluid
\begin{equation}
  q(\bx,t)=-\frac{k(\bx,t)}{\mu}\nabla p(\bx,t),
  \quad
  \nabla\cdot q(\bx,t)=-\partial_t\phi(\bx,t),
  \label{eq:darcy}
\end{equation}
where $\mu$ is the dynamic viscosity, $k$ is the intrinsic permeability, and the right-hand side accounts for the evolving pore space. 
Pressure is prescribed at the top and bottom of the column, no-flow conditions are imposed on the lateral boundaries, and a constant downward burial contribution is added to the Darcy discharge.

The effective diffusion coefficient is
\begin{equation*}
  D_i(\phi)=\frac{D_{i,0}}{\theta^2(\phi)}+D_b,
  \quad
  \theta^2(\phi)=1-2\ln\phi,
\end{equation*}
where $D_b$ represents mixing or biodiffusion, and $D_{i,0}=0$ for solid components.

The initial condition is $c_i(x,0)=c_{i,0}(x)$, and the boundary conditions are
\begin{equation}
  c_i=c_{i,D} \quad\text{on }\Gamma_{i,D},
  \quad
  \left(
    \eta_i v_i c_i
    -\eta_i D_i(\phi)\nabla c_i
  \right)\cdot n
  = j_i
  \quad \text{on } \Gamma_{i,N},
  \quad
  \partial\Omega=\Gamma_{i,D}\cup\Gamma_{i,N},
  \label{eq:concentration-boundaries}
\end{equation}
where $c_{i,D}$ is a prescribed concentration, $j_i$ a prescribed normal flux, and $n$ the outward unit normal vector.
On inflow boundaries, where $\eta_i v_i\cdot n<0$, the advective flux uses a prescribed inflow concentration $c_{i,\mathrm{in}}$, and all boundary data may differ among components.

\subsection{Geochemical reactions}
\label{sec:reactions}

In this work, we consider five components in marine sediment ($N_{\mathrm{comp}}=5$)
\begin{equation}
  c=(c_1,c_2,c_3,c_4,c_5)
  =\big(c_{\mathrm{CaCO_3}},c_{\mathrm{POC}},c_{\mathrm{SO_4}},
         c_{\mathrm{Ca}},c_{\mathrm{DIC}}\big),
  \label{eq:species}
\end{equation}
where calcium carbonate ($\mathrm{CaCO_3}$) and particulate organic carbon (POC) are solid ($\mathcal I_s=\{1,2\}$), and sulfate ($\mathrm{SO_4}$), calcium ($\mathrm{Ca}$), and dissolved inorganic carbon (DIC) are dissolved ($\mathcal I_f=\{3,4,5\}$). 

Two reactions act on them: organoclastic sulfate reduction  and carbonate precipitation and dissolution  with rates $R_{\mathrm{sr}}$ and $R_{\mathrm{cal}}$ \cite{chatterjee2011pore,lichtschlag2021impact,blouet2021seep,rahmati2021new,zander2020new}
\begin{equation*}
  R_{\mathrm{sr}}=k_{\mathrm{sr}}(z) c_2
    \min\!\left(\frac{c_3}{K_{\mathrm{MM}}},1\right),
  \quad
  R_{\mathrm{cal}}=\big(k_p-k_dc_1\big)
   \left(\frac{c_4c_5}{K_{\mathrm{sp}}}-1\right),
\end{equation*}
where $K_{\mathrm{MM}}$ is the sulfate saturation constant, $K_{\mathrm{sp}}$ is the solubility product, and $k_p,k_d$ are kinetic constants. 
{Sulfate reduction is microbially mediated, and its dependence on sulfate is commonly described by Michaelis--Menten kinetics, $c_3/(K_m+c_3)$ with the half-saturation constant $K_m$ \cite{boudreau1984sulfate}. As in \cite{blouet2021seep}, we use the piecewise-linear simplification $\min(c_3/K_{\mathrm{MM}},1)$ with the limiting concentration $K_{\mathrm{MM}}$, the linearized Monod kinetics of multicomponent models of early diagenesis in which the rate is proportional to sulfate below $K_{\mathrm{MM}}$ and does not depend on it above.}
The sulfate-reduction rate decreases with sediment age as $k_{\mathrm{sr}}(z)=k_{\mathrm{sr}}^0/(\tau_{\mathrm{age},0}+a\max\{0,-z-z_0\})$\cite{chatterjee2011pore,blouet2021seep},
and $R_{\mathrm{cal}}>0$ precipitates carbonate, whereas $R_{\mathrm{cal}}<0$ dissolves it
\cite{blouet2021seep,lichtschlag2021impact}. The component sources follow from the stoichiometry
\begin{equation}
  R=(R_1,\ldots,R_5)^T=S_{\mathrm{rxn}}\begin{pmatrix}R_{\mathrm{cal}}\\R_{\mathrm{sr}}\end{pmatrix},
  \quad
  S_{\mathrm{rxn}}=\begin{pmatrix}
    1&0&0&-1&-1\\
    0&-1&-\tfrac12&0&1
  \end{pmatrix}^T,
  \label{eq:stoichiometric-matrix}
\end{equation}
whose columns give the species changes per unit extent of each reaction. In a reaction-only step,
$\dot c_i=R_i$, so the reactions conserve $c_1+c_4$ and $2c_3-c_2$.

\begin{figure}[htbp]
\centering
\resizebox{0.7\textwidth}{!}{
\begin{tikzpicture}[
  pool/.style={draw=black!65,rounded corners=2pt,fill=blue!4,minimum width=1.7cm,minimum height=0.8cm,align=center,font=\small},
  solid/.style={draw=black!65,rounded corners=2pt,fill=gray!12,minimum width=1.7cm,minimum height=0.8cm,align=center,font=\small},
  rxn/.style={draw=black!65,rounded corners=6pt,fill=white,inner sep=3pt,align=center,font=\footnotesize},
  arrow/.style={-{Latex[length=2.2mm]},thick},
  both/.style={{Latex[length=2.2mm]}-{Latex[length=2.2mm]},thick},
  lab/.style={font=\footnotesize,align=center}
]
\draw[draw=black!35,rounded corners=3pt] (-1.5,2.4) rectangle (7.9,3.4);
\draw[draw=black!35,rounded corners=3pt] (-1.5,-0.2) rectangle (7.9,0.8);
\node[lab,anchor=south west] at (-1.3,3.4) {Pore fluid (dissolved)};
\node[lab,anchor=north west] at (-1.3,-0.1) {Sediment matrix (solid)};
\node[pool] (so4) at (0,2.9) {SO$_4$};
\node[pool] (dic) at (3.2,2.9) {DIC};
\node[pool] (ca) at (6.4,2.9) {Ca};
\node[solid] (poc) at (0,0.3) {POC};
\node[solid] (cal) at (6.4,0.3) {CaCO$_3$};
\node[rxn] (rsr) at (1.6,1.575) {sulfate reduction\\(microbial), $R_{\mathrm{sr}}$};
\node[rxn] (rcal) at (4.8,1.575) {precipitation,\\dissolution, $R_{\mathrm{cal}}$};
\draw[arrow] (so4) -- (rsr);
\draw[arrow] (poc) -- (rsr);
\draw[arrow] (rsr) -- (dic);
\draw[both] (dic) -- (rcal);
\draw[both] (ca) -- (rcal);
\draw[both] (rcal) -- (cal);
\draw[arrow] (-3.2,2.9) -- (-1.5,2.9) node[lab,midway,above] {inflow};
\draw[arrow] (7.9,2.9) -- (9.3,2.9) node[lab,midway,above] {outflow};
\node[lab] at (-2.3,2.5) {advection,\\diffusion};
\draw[arrow] (-2.9,0.3) -- (-1.5,0.3) node[lab,midway,above] {burial};
\draw[arrow] (7.9,0.3) -- (9.3,0.3) node[lab,midway,above] {burial};
\node[lab] at (-2.1,-0.3) {mixing};
\draw[arrow,dashed] (cal.east) ++(0,0.0) -- ++(1.3,1.75) node[lab,right] {porosity $\phi$,\\permeability $k$};
\end{tikzpicture}
}
\caption{{Reduced reaction network in a sediment cell. Sulfate reduction, which is microbially mediated, converts POC and sulfate to DIC, and carbonate precipitation and dissolution exchange Ca and DIC with CaCO$_3$, which changes the porosity and permeability. The dissolved components are transported with the pore fluid by advection and diffusion, and the solid components by burial and mixing.}}
\label{fig:reaction-network}
\end{figure}
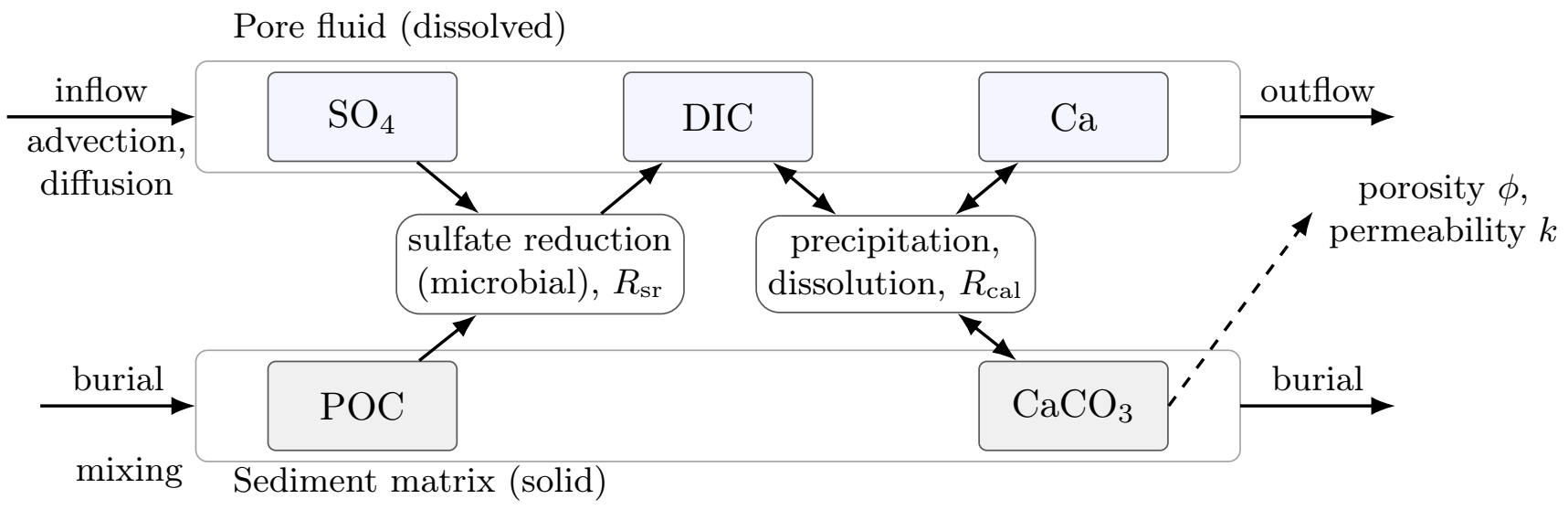

{Figure~\ref{fig:reaction-network} summarizes this reduced reaction network.
DIC is the sum of CO$_2$, HCO$_3^-$, and CO$_3^{2-}$, and the carbonate speciation and the pH are not modeled, in contrast to \cite{blouet2021seep}.
The saturation state in $R_{\mathrm{cal}}$ therefore uses the Ca and DIC concentrations with an effective solubility product $K_{\mathrm{sp}}$.
The degradation of POC to DIC is one overall reaction \cite{arndt2013quantifying}, so dissolved organic carbon (DOC) is not a separate component.
The focus of this work is the numerical method, and the extension to carbonate speciation and DOC is left for future work.}

\subsection{Random medium and porosity update}

The reference permeability and the initial porosity are derived from the same random field \cite{vasilyeva2021preconditioning} \begin{equation}
  k_{\mathrm{ref}}(\bx,\omega_{\mathrm{r}})
  = k_0\exp\!\big(\sigma_k Y_k(\bx,\omega_{\mathrm{r}})\big),
  \quad
  \phi_0(\bx,\omega_{\mathrm{r}})
  = 
    \bar\phi+\sigma_\phi Y_k(\bx,\omega_{\mathrm{r}}).
  \label{eq:random-fields}
\end{equation}
Here $\omega_{\mathrm{r}}$ labels a realization and $Y_k$ is a standardized Gaussian field obtained from a truncated Karhunen--Lo\`eve expansion (KLE) with the $N_{\mathrm{KL}}$ leading eigenpairs $(\lambda_j,e_j)$ of the discretized covariance operator
\begin{equation*}
  \widetilde{Y}_k(\bx,\omega_{\mathrm{r}})
  =\sum_{j=1}^{N_{\mathrm{KL}}}\sqrt{\lambda_j}\,\zeta_j(\omega_{\mathrm{r}})e_j(\bx),
  \quad
  Y_k(\bx,\omega_{\mathrm{r}})
  =\frac{\widetilde{Y}_k(\bx,\omega_{\mathrm{r}})
          -\overline{\widetilde{Y}_k}(\omega_{\mathrm{r}})}
         {\operatorname{std}_{\bx}[\widetilde{Y}_k(\bx,\omega_{\mathrm{r}})]},
  \quad
  \zeta_j\overset{\mathrm{iid}}{\sim}\mathcal{N}(0,1),
\end{equation*}
where the overline and $\operatorname{std}_{\bx}$ denote the spatial sample mean and standard deviation. Since $Y_k$ has zero spatial mean and unit spatial standard deviation, $k_0$ is the geometric mean permeability, $\sigma_k$ is the standard deviation of $\ln k_{\mathrm{ref}}$, and $\bar\phi$ and $\sigma_\phi$ are the mean and standard deviation of the initial porosity. The covariance is anisotropic exponential
\begin{equation*}
  C(\bx,\bx')=
  \exp\!\left[-\left(
  \sum_{q=1}^{d}\frac{(x_q-x'_q)^2}{\ell_q^2}
  \right)^{1/2}\right],
\end{equation*}
with correlation lengths $\ell_q$. 

Porosity evolves with mineral precipitation and dissolution  \cite{steefel2015codes, minto2019porosity}
\begin{equation}
\frac{\partial \phi}{\partial t}(\bx,t,\omega_{\mathrm{r}})
=
-\gamma_\phi
R_1\big(c(\bx,t,\omega_{\mathrm{r}})\big),
  \label{eq:porosity-update}
\end{equation}
where $\gamma_\phi>0$ converts a mineral-concentration change to a pore-volume
change.  

Permeability follows from the updated porosity through a normalized power-law relation \cite{hommel2018poroperm}
\begin{equation}
  k(\bx,t,\omega_{\mathrm{r}})
  = k_{\mathrm{ref}}(\bx,\omega_{\mathrm{r}})
    \left(\frac{\phi(\bx,t,\omega_{\mathrm{r}})}{\phi_0(\bx,\omega_{\mathrm{r}})}\right)^{m_k}
    \left(\frac{1-\phi_0(\bx,\omega_{\mathrm{r}})}{1-\phi(\bx,t,\omega_{\mathrm{r}})}\right)^{n_k},
  \label{eq:permeability-update}
\end{equation}
where $m_k$ and $n_k$ are constitutive exponents.

\section{Split: Physics-Based Operator Splitting and Discretization}
\label{sec:numerics}

We first split the coupled evolution into flow, species-wise transport, local reaction, and porosity--permeability update stages. The decomposition defines the sequence of subproblems advanced at each time step and provides the structure used later for split learning. To solve the resulting subproblems and generate the simulation datasets, we discretize the governing equations in space with a finite-volume method (FVM).

\subsection{Operator splitting and time approximation}

Operator splitting is widely used in reactive-transport simulation to separate the flow and transport solves from the local, chemically coupled reaction solve, allowing each subproblem to be treated with a method suited to its structure \cite{vasilyeva2024decoupled,steefel2015codes,carrayrou2004splitting}.  
At the continuum level, each time interval is split into four sequential stages.
First, the Darcy problem is solved for the current porosity and permeability.
Second, the reaction-free transport equation is solved for each component with the flow and porosity fixed.
Third, the local reaction system $\partial_t(\eta_ic_i)=\eta_iR_i(c)$ is solved.
Finally, porosity and permeability are updated for the next flow solve. 
We note that sequential composition introduces operator-splitting error because the subproblems are advanced one after another rather than simultaneously. In general, this error depends on the time step and on the interaction between the operators.  In first-order non-iterative sequential scheme used here, porosity and permeability are updated after reaction and affect the flow and transport solve at the next time step.

\begin{figure}[htbp]
\centering
\resizebox{0.8\linewidth}{!}{
\begin{tikzpicture}[
  stagebox/.style={draw=black!65,rounded corners=2pt,fill=blue!4,
    text width=3.65cm,minimum height=1.15cm,align=center,font=\small},
  state/.style={draw=black!65,rounded corners=2pt,fill=gray!4,
    text width=3.65cm,minimum height=0.75cm,align=center,font=\small},
  update/.style={draw=black!65,rounded corners=2pt,fill=gray!4,
    text width=6.6cm,minimum height=1.05cm,align=center,font=\small},
  arrow/.style={-{Latex[length=2.1mm]},thick},
  lab/.style={font=\footnotesize,fill=white,inner sep=2pt}
]
\node[state] (current) at (0,2.0)
  {Current state\\$(\phi^n,k^n,c^n)$};
\node[stagebox] (flow) at (0,0)
  {Flow\\Solve for $p^n$ and $q^n$};
\node[stagebox] (transport) at (4.8,0)
  {Transport\\Advance $c^n\longmapsto c^*$};
\node[stagebox] (reaction) at (9.6,0)
  {Local reaction\\Advance $c^*\longmapsto c^{n+1}$};
\draw[arrow] (current.south) -- node[lab,right] {$\phi^n,k^n$} (flow.north);
\draw[arrow] (current.east) -- (4.8,2.0) --
  node[lab,right] {$c^n,\phi^n,k^n$} (transport.north);
\draw[arrow] (flow.east) -- node[lab,above] {$q^n$} (transport.west);
\draw[arrow] (transport.east) -- node[lab,above] {$c^*$} (reaction.west);
\node[update] (next) at (4.8,-1.7)
  {Porosity and permeability update from mineral change $(\phi^{n+1},k^{n+1},c^{n+1})$};
\draw[arrow] (reaction.south) -- node[lab,right] {$c^*,c^{n+1}$} (9.6,-1.7) -- (next.east);
\draw[arrow] (next.west) -- (-2.6,-1.7) -- (-2.6,2.0) -- (current.west);
\node[lab] at (-2.6,0.8) {Next step};
\end{tikzpicture}
}
\caption{One step of the sequential splitting algorithm. The current medium determines the Darcy flux, which drives transport. The post-transport state is advanced by local reaction, and the mineral change updates porosity and permeability for the next time step.}
\label{fig:splitting}
\end{figure}
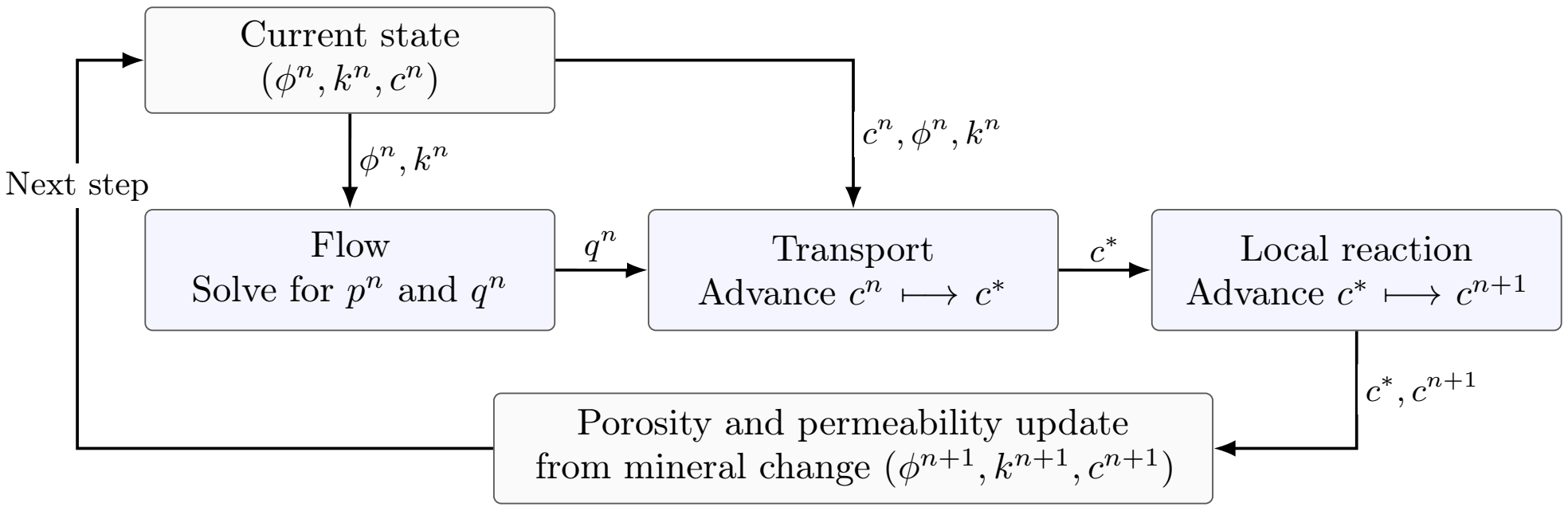

For the time approximation, we use first-order backward Euler for both the transport and reaction stages.  Transport is solved implicitly at the new time level, followed by an implicit reaction solve.  Let $t_n=n\Delta t$, with current porosity, permeability, and concentrations $\phi^n$, $k^n$, and $c^n$.  
A single time step proceeds as follows
\begin{itemize}
  \item \textit{Flow}: solve Darcy's problem using the current $\phi^n$ and $k^n$,
        with the prescribed pressure and no-flow boundary conditions
        \begin{equation*}
          q^n=-\frac{k^n}{\mu}\nabla p^n,
            \quad
            \nabla\cdot q^n
            =
            -\frac{\phi^n-\phi^{n-1}}{\Delta t}.
        \end{equation*}
        Set $\bv_f^n=q^n/\phi^n$ and use the
        phase velocities defined in \eqref{eq:phase-weights}.
  \item \textit{Transport}:  apply backward
        Euler to the reaction-free continuum transport equation
        \begin{equation}
          \frac{
            \eta_i(\phi^n)c_i^*
            -\eta_i(\phi^n)c_i^n}{\Delta t}
          +\nabla\cdot\left(
            \eta_i(\phi^n)\bv_i^n c_i^*
          \right)
          -\nabla\cdot\left(
            \eta_i(\phi^n)D_i(\phi^n)\nabla c_i^*
          \right)=0,\quad i=1,\ldots,5.
          \label{eq:transport-be}
        \end{equation}
        Impose the boundary conditions in
        \eqref{eq:concentration-boundaries} at the new time level to obtain
        $c_i^*$.
  \item \textit{Reaction}: solve the backward-Euler reaction system locally
        \begin{equation}
          c^{n+1}-c^*
          =\Delta t\,R(c^{n+1}).
          \label{eq:reaction-be}
        \end{equation}
  \item \textit{Porosity and permeability update} from the mineral change
        \begin{align*}
          \phi^{n+1}(\bx)
            &=\phi^n(\bx)-\gamma_\phi
              \big(c_1^{n+1}(\bx)-c_1^*(\bx)\big),\\
          k^{n+1}(\bx)
            &=k_{\mathrm{ref}}(\bx)
              \left(\frac{\phi^{n+1}(\bx)}{\phi_0(\bx)}\right)^{m_k}
              \left(\frac{1-\phi_0(\bx)}
                         {1-\phi^{n+1}(\bx)}\right)^{n_k}.
        \end{align*}
        The updated porosity and permeability are used in the next step.
\end{itemize}
Figure~\ref{fig:splitting} summarizes one time step. The concentrations advance through a single flow--transport--reaction pass, followed by the porosity and permeability update for the next step.
For a fixed medium, energy estimates as in \cite{vasilyeva2024decoupled} give stability and first-order accuracy in
time of the transport--reaction composition under monotonicity assumptions on the transport and reaction operators.

\subsection{Finite-volume approximation}

We use a {uniform} Cartesian finite-volume mesh, where each cell $K$ has volume $|K|$ and each face $e$ has area $|e|$ and outward unit normal $n_{K,e}$. 
The unknown $C_{i,K}^n$ represents the cell average of concentration $c_i$ at time $t_n$.
Integrating the transport equation over cell $K$ gives
\begin{equation}
  \frac{|K|}{\Delta t}
    \left[\eta_i(\phi_K^n)C_{i,K}^*
          -\eta_i(\phi_K^n)C_{i,K}^n\right]
  +\sum_{e\subset\partial K}
    \left(F^{\mathrm{adv}}_{i,K,e}+F^{\mathrm{diff}}_{i,K,e}\right)
  =0,
  \label{eq:fv-balance}
\end{equation}
where reaction is advanced separately by the cell-local step \eqref{eq:reaction-be}. 
On an interior face shared by cells $K$ and $L$, the diffusive flux is approximated by a two-point flux 
\begin{equation}
F^{\mathrm{diff}}_{i,K,e}
= T^{D,n}_{i,KL}
(C_{i,K}^*-C_{i,L}^*),
\quad
T^{D,n}_{i,KL}=\frac{|e|}{d_{KL}}
(\eta_iD_i)_{KL}^{n},
\label{eq:fv-diffusion}
\end{equation}
where $d_{KL}$ is the distance between cell centers, $(\eta_iD_i)_{KL}^{n}$ is the harmonic mean of the cell values of $\eta_iD_i(\phi^n)$ in $K$ and $L$, and the same face flux with opposite signs in $K$ and $L$ makes the scheme locally conservative.
The oriented advective flux for the current flow and porosity is
\[
F^{\mathrm{adv}}_{i,K,e}
=Q_{i,K,e}^{n}C^{\mathrm{up},*}_{i,e},
\]
where $Q_{i,K,e}^n$ is the outward phase flux $\eta_i(\phi^n)\bv_i^n$ integrated over face $e$, which is the Darcy face flux $q_{K,e}^n$ for dissolved species, since $\eta_i\bv_i=\phi\bv_f=q$, and the solid-burial face flux of $(1-\phi)\bv_s$ for solid species.
On an interior face shared by $K$ and $L$,
$C^{\mathrm{up},*}_{i,e}=C_{i,K}^*$ when $Q_{i,K,e}^n\geq 0$ and
$C^{\mathrm{up},*}_{i,e}=C_{i,L}^*$ when $Q_{i,K,e}^n<0$.  
Dirichlet, diffusive-flux, and advective-inflow boundary conditions enter the face fluxes and cell right-hand sides.

The Darcy equation is discretized by the same conservative face-balance principle, with the two-point Darcy flux on an interior face shared by cells $K$ and $L$
\begin{equation}
  \sum_{e\subset\partial K}q_{K,e}^{n}
  =-\frac{|K|}{\Delta t} (\phi_K^n-\phi_K^{n-1}),
  \quad
  q_{K,e}^{n}=T_{KL}^{n}(p_K^{n}-p_L^{n}),
  \quad
  T_{KL}^{n}=\frac{|e|}{\mu d_{KL}}k_{KL}^{n},
  \label{eq:fv-darcy-flux}
\end{equation}
where $q_{K,e}^{n}$ is the outward Darcy flux on face $e$ of cell $K$ and $k_{KL}^{n}$ 
is the harmonic mean of the cell permeabilities.

Let $N_h$ denote the number of cells and $\mathcal N(K)$ the cells sharing an interior face with $K$. Collecting \eqref{eq:fv-darcy-flux} over all cells gives the sparse pressure system, $p^n\in\R^{N_h}$
\begin{equation}
  A_p(k^n)p^n=b_p^n,
  \label{eq:fv-flow-matrix}
\end{equation}
where the matrix $A_p(k^n)\in\R^{N_h\times N_h}$ is the two-point stencil of the transmissibilities
$T_{KL}^n$ on interior faces and $T_{K,e}^n=|e|k_K^n/(\mu d_{K,e})$ on boundary faces with prescribed pressure $p_e^n$, where $d_{K,e}$ is the cell-center-to-face distance. The vector $b_p^n$ contains the prescribed-pressure contributions $T_{K,e}^np_e^n$ and the porosity-change source $-|K|(\phi_K^n-\phi_K^{n-1})/\Delta t$, and no-flow faces contribute zero flux.

For each species, $C_i=(C_{i,K})_K\in\R^{N_h}$ and the diagonal storage matrix is $M_i(\phi)=\operatorname{diag}(|K|\eta_i(\phi_K))$. 
The sparse transport matrix $A_i(k,\phi,q)\in\R^{N_h\times N_h}$ collects the upwind advective and diffusive face fluxes of \eqref{eq:fv-balance}--\eqref{eq:fv-diffusion}. The boundary-face fluxes contribute to the diagonal of $A_i$ (advective outflow and prescribed diffusive concentrations) and to the boundary vector $b_i^{n+1}\in\R^{N_h}$ (prescribed advective inflow and diffusive concentrations). 
The fully discrete transport step is the sparse linear system
\begin{equation}
  \left[
    \frac{M_i(\phi^n)}{\Delta t}
    +A_i(k^n,\phi^n,q^n)
  \right]C_i^*
  =\frac{M_i(\phi^n)}{\Delta t}C_i^n+b_i^{n+1},
  \quad i=1,\ldots,5,
  \label{eq:operator-form}
\end{equation}
The reaction step \eqref{eq:reaction-be} is then solved in each cell and couples the components locally
\begin{equation}
  C_K^{n+1}-C_K^*
  =\Delta t\,R(C_K^{n+1}).
  \label{eq:reaction-disc}
\end{equation}
Each time step thus solves the pressure system \eqref{eq:fv-flow-matrix}, the transport systems \eqref{eq:operator-form}, and the local reaction systems \eqref{eq:reaction-disc}, and Section~\ref{sec:learning} replaces each of them by a data-driven surrogate.

\section{Learn: Data-Driven Approximation of the Split Operators}
\label{sec:learning}

We construct a surrogate for each stage of the split reactive-transport update.
Transport is approximated in a reduced data-driven spatial space, for which we consider two approaches, a learned solution operator (global inverse atoms) and a learned forward transport operator assembled from local pieces and solved in the reduced space (local atoms).
Local chemistry is learned through bounded reaction extents and stoichiometric reconstruction. 
Flow is reconstructed from learned interface traces through a mixed multiscale solve.

\subsection{Learning the reduced transport}
\label{sec:data-driven-space}

The transport surrogate is built in three parts.
Subsection~\ref{sec:ms-space} constructs a local data-driven multiscale space and encoder from concentration snapshots, which reduce the state dimension and are shared by both transport models, and encodes the fine-grid transport step on this space in two ways, as a projected Galerkin or LSPG system that is solved online, or as a reduced solution map whose training data are fine-grid solutions projected onto the space.
Subsection~\ref{sec:global-atoms} learns the reduced solution map of the transport step by compressing reference maps into matrix atoms (global inverse atoms).
Subsection~\ref{sec:local-atoms} learns coarse-block pieces of the forward transport operator and solves a reduced projected system online (local atoms).

\subsubsection{Data-driven multiscale space and reduced systems}
\label{sec:ms-space}

We first build the data-driven multiscale space locally \cite{vasilyeva2026operator}.
We introduce a coarse grid $\mathcal T_H$ that partitions the domain $\Omega$ into $N_H$ coarse cells $K_c^H$, each a union of fine cells, with coarse nodes $x_j$.
The coarse neighborhood $\omega_j$ is the union of the coarse cells that share the node $x_j$, so the neighborhoods overlap and cover $\Omega$.
The species-$i$ snapshots restricted to each $\omega_j$ form the matrix $X_{i,j}${, whose columns are sampled over the training realizations and time steps}. 
A POD, computed by $X_{i,j}=U_{i,j}\Sigma_{i,j}V_{i,j}^{T}$, supplies local spatial modes ordered by decreasing singular value.
Let $U_{i,j}^{(m)}$ contain the first $m$ left singular vectors.
We select the smallest count $m$ for which the relative squared snapshot projection error satisfies
\begin{equation*}
 \frac{\|X_{i,j}-U_{i,j}^{(m)}(U_{i,j}^{(m)})^T X_{i,j}\|_F^2}
      {\|X_{i,j}\|_F^2}
 =\frac{\sum_{\ell>m}\sigma_{i,j,\ell}^2}
        {\sum_{\ell}\sigma_{i,j,\ell}^2}
 =1-\frac{\sum_{\ell=1}^{m}\sigma_{i,j,\ell}^{2}}{\sum_{\ell}\sigma_{i,j,\ell}^{2}}
 \leq1-\rho_{\mathrm{tr}},
\end{equation*}
where $0<\rho_{\mathrm{tr}}\leq1$ is the retained energy fraction.
A prescribed mode cap may prevent attainment of this threshold.
When it is attained, the relative local snapshot error is at most $\sqrt{1-\rho_{\mathrm{tr}}}$. 
The bound concerns the training snapshots before spatial weighting.
On each $\omega_j$, we add the constant function to the selected modes, orthonormalize them, and discard linearly dependent vectors. 
Multiplication by {piecewise linear} partition-of-unity weights $\chi_j$, with $\sum_j\chi_j=1$, followed by zero extension and normalization gives the multiscale functions $\varphi_{i,j,\ell}$. 
The multiscale functions span the space $V_i^{\mathrm{ms}}$, in which the species-$i$ concentration is approximated by
\begin{equation*}
V_i^{\mathrm{ms}}
 =\operatorname{span}\{\varphi_{i,j,\ell}\}_{j,\ell},
 \quad c_i\approx\sum_{j}\sum_{\ell}a_{i,j,\ell}\,\varphi_{i,j,\ell}.
\end{equation*}
Collecting the functions $\varphi_{i,j,\ell}$ as columns gives the sparse reconstruction matrix $P_i\in\mathbb R^{N_h\times r_i}$, where $a_i\in\mathbb R^{r_i}$ collects the coefficients $a_{i,j,\ell}$, $N_h$ is the number of fine cells, and $r_i$ is the reduced dimension. 
The vector $C_i$ is reconstructed from $a_i$, and the least-squares encoder $E_i$ determines $a_i$ from $C_i$
\begin{equation*}
C_i\approx P_ia_i,
 \quad 
 a_i=E_iC_i,
 \quad E_iC_i=\underset{a\in\mathbb R^{r_i}}{\operatorname{argmin}}
 \|P_i a-C_i\|_2^2,
\end{equation*}
where the minimizer solves the normal equations $(P_i^TP_i)\,a_i=P_i^TC_i$, whose sparse matrix $P_i^TP_i$ is factorized once offline with a small diagonal shift.
We note that the basis and encoder are fixed offline and shared across realizations.
Furthermore, we use this space to reduce the storage of the transport operators, which act on the $r_i$ coefficients $a_i$ instead of the $N_h$ fine-grid values, and to represent the solution, which is advanced in the coefficients $a_i$ and reconstructed as $P_ia_i$.

Both transport models use reduced operators that are generated offline from the fine-grid transport step \eqref{eq:operator-form}.
For the initial medium $\kappa=(k,\phi,q)$ of each training realization, we assemble the finite-volume matrix $A_i(\kappa)=A_i(k,\phi,q)$, the storage matrix $M_i(\phi)$, and the boundary vector $b_i(\kappa)$ of \eqref{eq:operator-form}, with the flow and boundary conditions of the coupled solver, and write the step as
\begin{equation}
  \mathcal L_i(\kappa)\,C_i^*=\tfrac{1}{\Delta t}M_i(\phi)C_i^n+b_i(\kappa),
  \quad
  \mathcal L_i(\kappa)=\tfrac{1}{\Delta t}M_i(\phi)+A_i(\kappa).
  \label{eq:training-system}
\end{equation}
The system \eqref{eq:training-system} is the fine-grid system \eqref{eq:operator-form} of size $N_h$, evaluated at $\kappa$. Next, together with $P_i$ and $E_i$, these matrices define the reduced systems of dimension $r_i$ learned in the next two subsections, in two different ways.
The first classical way projects the equation \eqref{eq:training-system} onto $V_i^{\mathrm{ms}}$ and solves the projected system with the Galerkin (G) or the least-squares Petrov--Galerkin (LSPG) projections
\begin{equation}
\begin{aligned}
  \text{G:}\quad &
  \bigl(P_i^T\mathcal L_i(\kappa)P_i\bigr)\,a_i^*
  =P_i^T\bigl(\tfrac{1}{\Delta t}M_i(\phi)C_i^n+b_i(\kappa)\bigr),\\
  \text{LSPG:}\quad &
  \bigl(\mathcal L_i(\kappa)P_i\bigr)^T\bigl(\mathcal L_i(\kappa)P_i\bigr)\,a_i^*
  =\bigl(\mathcal L_i(\kappa)P_i\bigr)^T\bigl(\tfrac{1}{\Delta t}M_i(\phi)C_i^n+b_i(\kappa)\bigr),
\end{aligned}
  \label{eq:reduced-local}
\end{equation}
with matrices in $\R^{r_i\times r_i}$. 
The local atoms learn these systems, whose matrices and right-hand sides are assembled from coarse-block pieces, and the projected system is solved online on the reduced space. 
The second way projects the fine-grid solution, so the step \eqref{eq:training-system} is solved with the initial state $C_i^n=P_ia_i^n$ and the result is encoded by $E_i$, the least-squares projection onto $V_i^{\mathrm{ms}}$, which gives
\begin{equation}
  a_i^*=G_i(\kappa)\,a_i^n+g_i(\kappa),
  \quad
  G_i(\kappa)=E_i\mathcal L_i(\kappa)^{-1}\tfrac{1}{\Delta t}M_i(\phi)P_i,
  \quad
  g_i(\kappa)=E_i\mathcal L_i(\kappa)^{-1}b_i(\kappa),
  \label{eq:reduced-global}
\end{equation}
with $G_i\in\R^{r_i\times r_i}$ and $g_i\in\R^{r_i}$. 
The fine-grid solves in \eqref{eq:reduced-global} are used only offline to generate the training data, and the global inverse atoms learn the map $(G_i,g_i)$, so the online step is a reduced matrix--vector product without a fine-grid solve.

\subsubsection{Global inverse atoms}
\label{sec:global-atoms}

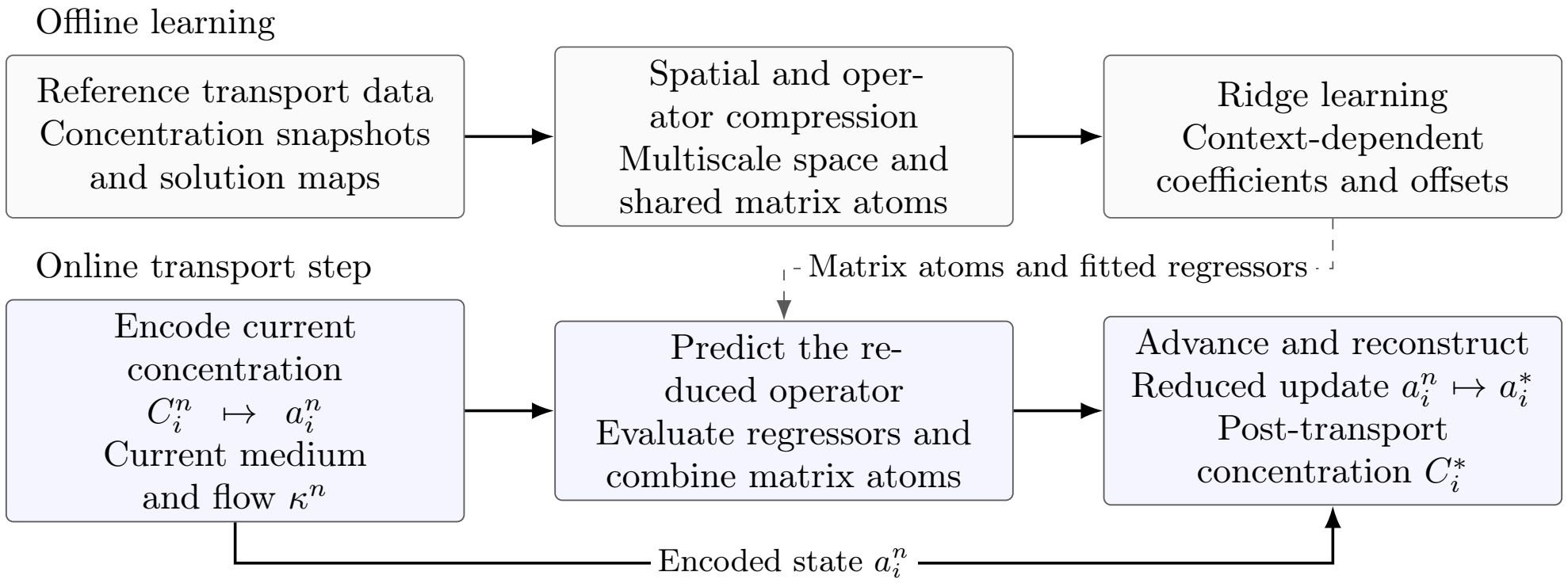
\begin{figure}[htbp]
\centering
\resizebox{0.8\textwidth}{!}{
\begin{tikzpicture}[
  offline/.style={draw=black!60,rounded corners=2pt,fill=gray!4,
    text width=4.25cm,minimum height=1.6cm,align=center,font=\small},
  stage/.style={draw=black!65,rounded corners=2pt,fill=blue!4,
    text width=4.25cm,minimum height=1.65cm,align=center,font=\small},
  arrow/.style={-{Latex[length=2.2mm]},thick},
  trained/.style={-{Latex[length=2mm]},dashed,black!65},
  lab/.style={font=\footnotesize,fill=white,inner sep=2pt}
]
\node[anchor=west,font=\small] at (-2.1,4.3) {Offline learning};
\node[offline] (data) at (0,3.2)
  {Reference transport data\\Concentration snapshots and solution maps};
\node[offline] (compress) at (5.4,3.2)
  {Spatial and operator compression\\Multiscale space and shared matrix atoms};
\node[offline] (train) at (10.8,3.2)
  {Ridge learning\\Context-dependent coefficients and offsets};
\draw[arrow] (data.east) -- (compress.west);
\draw[arrow] (compress.east) -- (train.west);
\node[anchor=west,font=\small] at (-2.1,1.9) {Online transport step};
\node[stage] (input) at (0,0.5)
  {Encode current concentration\\$C_i^n\mapsto a_i^n$\\Current medium and flow $\kappa^n$};
\node[stage] (predict) at (5.4,0.5)
  {Predict the reduced operator\\Evaluate regressors and combine matrix atoms};
\node[stage] (output) at (10.8,0.5)
  {Advance and reconstruct\\Reduced update $a_i^n\mapsto a_i^*$\\Post-transport concentration $C_i^*$};
\draw[arrow] (input.east) -- (predict.west);
\draw[arrow] (predict.east) -- (output.west);
\draw[trained] (train.south) -- (10.8,1.9) -- (5.4,1.9) -- (predict.north);
\node[lab] at (8.1,1.9) {Matrix atoms and fitted regressors};
\draw[arrow] (input.south) -- (0,-1.) -- (10.8,-1.) -- (output.south);
\node[lab] at (5.4,-1.) {Encoded state $a_i^n$};
\end{tikzpicture}
}
\caption{Transport learning and prediction. 
A multiscale representation reduces the state dimension, while operator compression and ridge regression provide the context-dependent reduced map. 
Reconstruction supplies concentrations for local chemistry.}
\label{fig:transport-training-overview}
\end{figure}

The global inverse atoms learn the reduced map \eqref{eq:reduced-global} as a function of the medium and flow.
They are constructed on the multiscale space of Subsection~\ref{sec:ms-space} in two offline stages and one online stage in every time step, summarized in Figure~\ref{fig:transport-training-overview}.

\noindent\textit{Stage 1: Reference maps.}
For each training medium $\kappa$, the reference map $(G_i(\kappa),g_i(\kappa))$ of \eqref{eq:reduced-global} is computed from \eqref{eq:training-system} by one sparse LU factorization of $\mathcal L_i(\kappa)$ and $r_i+1$ forward and backward substitutions followed by the encoder $E_i$.
Here $G_i$ propagates the encoded concentration and $g_i$ accounts for the boundary forcing.
Each training medium gives one pair $(G_i,g_i)$ and its context $x_{\mathrm{ctx}}(\kappa)$, defined below, and these pairs form the data for Stage 2.
The maps could also be generated from the Galerkin or LSPG systems \eqref{eq:reduced-local}, for example $G_i^{\mathrm{G}}=(P_i^T\mathcal L_iP_i)^{-1}P_i^T\tfrac{1}{\Delta t}M_iP_i$, with one reduced solve per medium instead of a fine-grid factorization and $r_i+1$ fine-grid substitutions, and Stages 2 and 3 do not depend on this choice.
We use \eqref{eq:reduced-global}, which gives the best approximation of the fine-grid step in $V_i^{\mathrm{ms}}$.

\noindent\textit{Stage 2: Operator compression and regression.}
We compress the increments $J_i=G_i-I$ of the training maps into $K_i$ global reduced matrix atoms $\Theta_{i,k}\in\R^{r_i\times r_i}$, obtained by singular-value decomposition of the vectorized training increments, so that each map is represented by $K_i$ coefficients.
Subtracting $I$ removes the unchanged-state contribution, so that only the transport increment is compressed. 
We predict the atom coefficients and boundary offset by ridge regression
\begin{equation*}
  \begin{aligned}
    \widehat J_i(\kappa)&=\sum_{k=1}^{K_i}
       \widehat\alpha_{i,k}(\kappa)\Theta_{i,k}, \quad 
    \widehat\alpha_i(\kappa)=W_i x_{\mathrm{ctx}}(\kappa)+b_{\alpha,i},\quad
    \widehat g_i(\kappa)=V_i x_{\mathrm{ctx}}(\kappa)+b_{g,i},
  \end{aligned}
\end{equation*}
where $K_i$ is the atom count and the context $x_{\mathrm{ctx}}(\kappa)$ contains local projections of log permeability, porosity, and flow, $\kappa=(k,\phi,q)$.
The fitted matrices $W_i,V_i$ and biases $b_{\alpha,i},b_{g,i}$ predict atom weights and offsets.

\noindent\textit{Stage 3: Online transport advancement.}
At each time step, the context $x_{\mathrm{ctx}}(\kappa^n)$ is computed from the current medium and flow $\kappa^n=(k^n,\phi^n,q^n)$ and each species is advanced by
\begin{equation*}
  \begin{aligned}
    a_i^n&=E_iC_i^n,\quad 
    a_i^*=(I+\widehat J_i(\kappa^n))a_i^n
                    +\widehat g_i(\kappa^n),
  \end{aligned}
\end{equation*}
where the reconstruction $C_i^*=\max\{P_i a_i^*,0\}$ enforces cellwise nonnegativity. 
The basis, atoms, and regressors remain fixed, while the predicted map changes with the current context. 
Accuracy depends on spatial reduction, atom truncation, and regression.
Online transport requires no fine-grid assembly and no linear solve, only a matrix--vector product per species.
The shared atom expansion relates to empirical operator approximation \cite{negri2015mdeim} and reduced-operator interpolation \cite{amsallem2011interpolation}.
It also relates to operator inference \cite{peherstorfer2016opinf}, but uses reference solution maps as targets.

\subsubsection{Local atoms}
\label{sec:local-atoms}

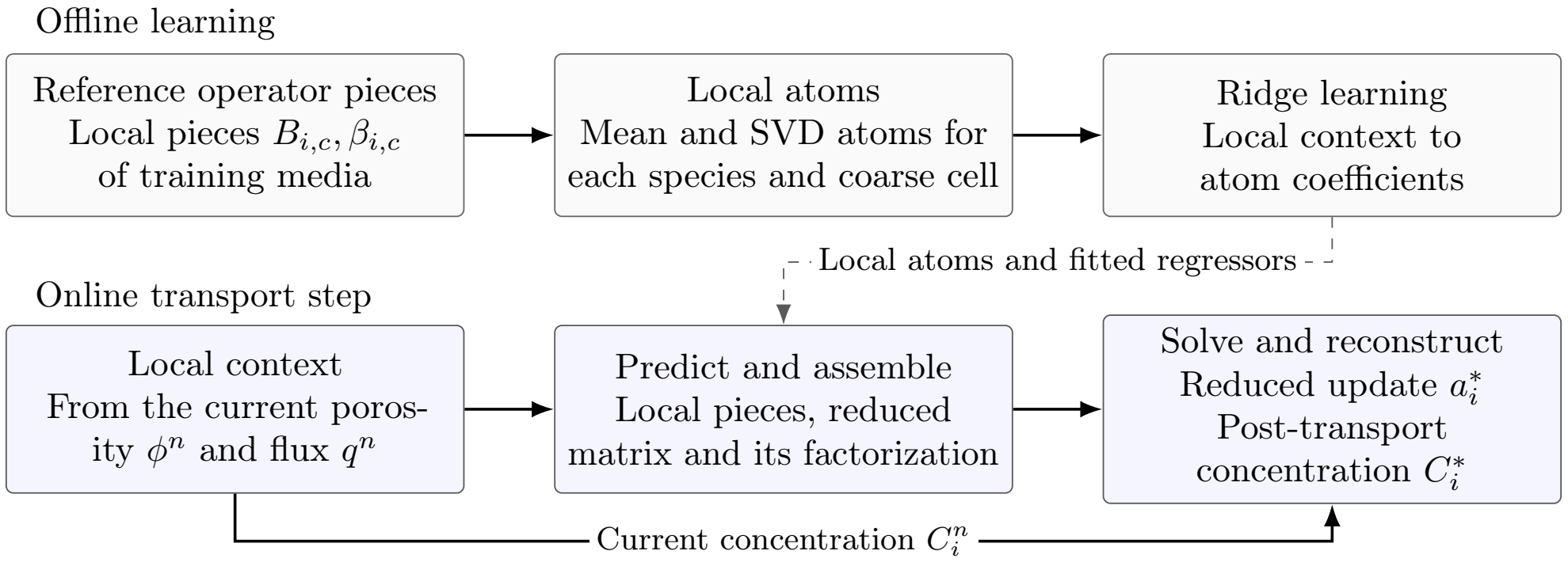
\begin{figure}[htbp]
\centering
\resizebox{0.8\textwidth}{!}{
\begin{tikzpicture}[
  offline/.style={draw=black!60,rounded corners=2pt,fill=gray!4,
    text width=4.25cm,minimum height=1.6cm,align=center,font=\small},
  stage/.style={draw=black!65,rounded corners=2pt,fill=blue!4,
    text width=4.25cm,minimum height=1.65cm,align=center,font=\small},
  arrow/.style={-{Latex[length=2.2mm]},thick},
  trained/.style={-{Latex[length=2mm]},dashed,black!65},
  lab/.style={font=\footnotesize,fill=white,inner sep=2pt}
]
\node[anchor=west,font=\small] at (-2.1,4.3) {Offline learning};
\node[offline] (data) at (0,3.2)
  {Reference operator pieces\\Local pieces $B_{i,c},\beta_{i,c}$ of training media};
\node[offline] (compress) at (5.4,3.2)
  {Local atoms\\Mean and SVD atoms for each species and coarse cell};
\node[offline] (train) at (10.8,3.2)
  {Ridge learning\\Local context to atom coefficients};
\draw[arrow] (data.east) -- (compress.west);
\draw[arrow] (compress.east) -- (train.west);
\node[anchor=west,font=\small] at (-2.1,1.6) {Online transport step};
\node[stage] (input) at (0,0.5)
  {Local context\\From the current porosity $\phi^n$ and flux $q^n$};
\node[stage] (predict) at (5.4,0.5)
  {Predict and assemble\\Local pieces, reduced matrix and its factorization};
\node[stage] (output) at (10.8,0.5)
  {Solve and reconstruct\\Reduced update $a_i^*$\\Post-transport concentration $C_i^*$};
\draw[arrow] (input.east) -- (predict.west);
\draw[arrow] (predict.east) -- (output.west);
\draw[trained] (train.south) -- (10.8,1.95) -- (5.4,1.95) -- (predict.north);
\node[lab] at (8.1,1.95) {Local atoms and fitted regressors};
\draw[arrow] (input.south) -- (0,-0.8) -- (10.8,-0.8) -- (output.south);
\node[lab] at (5.4,-0.8) {Current concentration $C_i^n$};
\end{tikzpicture}
}
\caption{Local atoms for transport. 
Local pieces of the transport operator are compressed into atoms and regressed on the local medium and flow offline. In each time step, the predicted pieces give a reduced system, whose solution gives the concentrations for local chemistry.}
\label{fig:local-atoms-overview}
\end{figure}

The global inverse atoms of Subsection~\ref{sec:global-atoms} approximate the reduced inverse transport map $G_i(\kappa)$.
The map $G_i(\kappa)$ is dense and depends on the medium and flow in the whole domain, so its accurate approximation requires several hundred global reduced matrix atoms.
In the local atoms, we compress the discrete transport operator before inversion.
The finite-volume matrix couples only neighboring cells, and its entries depend on the local coefficients.
Its coarse-block pieces are therefore local functions of the medium and admit compact local expansions.
The local atoms are constructed on the same multiscale space in two offline stages and one online stage in every time step, summarized in Figure~\ref{fig:local-atoms-overview}.

\noindent\textit{Stage 1: Local operator pieces.} 
For the coarse cells $K_c^H$ of Subsection~\ref{sec:ms-space}, let $R_c$ restrict a fine-grid vector to the cells of $K_c^H$. 
Let $\mathcal J_{i,c}$ contain the indices of the multiscale functions $\varphi_{i,j,\ell}$ that do not vanish in $K_c^H$ or in the fine-cell layers adjacent to it, and let $P_{i,c}=R_cP_i|_{\mathcal J_{i,c}}$ denote the corresponding columns of $P_i$ on the rows of $K_c^H$.
Since $A_i$ has the two-point stencil, $R_cA_iP_i$ vanishes outside the columns $\mathcal J_{i,c}$.
We define the local operator and boundary pieces
\begin{equation}
 B_{i,c}(\kappa)=R_cA_i(\kappa)P_i|_{\mathcal J_{i,c}},
 \quad
 \beta_{i,c}(\kappa)=R_cb_i(\kappa).
 \label{eq:local-pieces}
\end{equation}
The storage term $M_i(\phi)/\Delta t$ is diagonal and is evaluated exactly from the current porosity.
For a Galerkin projection, only the reduced projected pieces are needed
\begin{equation*}
 X_{i,c}=P_{i,c}^TB_{i,c}\in\R^{|\mathcal J_{i,c}|\times|\mathcal J_{i,c}|},
 \quad
 y_{i,c}=P_{i,c}^T\beta_{i,c}\in\R^{|\mathcal J_{i,c}|},
\end{equation*}
where the pieces are taken directly from the system \eqref{eq:training-system} of each training realization ($A_i(\kappa)$ and $b_i$) and  $M_{i,c}=R_cM_i(\phi^n)R_c^T/\Delta t$.
Given the pieces, the Galerkin reduced system for \eqref{eq:operator-form} is
\begin{equation}
 \Bigl[\sum_{c=1}^{{N_H}}S_{i,c}\bigl(P_{i,c}^TM_{i,c}P_{i,c}+X_{i,c}\bigr)S_{i,c}^T\Bigr]a_i^*
 =\sum_{c=1}^{{N_H}}S_{i,c}\bigl(P_{i,c}^TM_{i,c}R_cC_i^n+y_{i,c}\bigr),
 \label{eq:local-galerkin}
\end{equation}
where $S_{i,c}\in\{0,1\}^{r_i\times|\mathcal J_{i,c}|}$ embeds the local index set $\mathcal J_{i,c}$ into the reduced space.
Alternatively, the LSPG system in \eqref{eq:reduced-local} is assembled from the local residual factors $\Gamma_{i,c}=M_{i,c}P_{i,c}+B_{i,c}$ as
\begin{equation}
 \Bigl[\sum_{c=1}^{{N_H}}S_{i,c}\Gamma_{i,c}^T\Gamma_{i,c}S_{i,c}^T\Bigr]a_i^*
 =\sum_{c=1}^{{N_H}}S_{i,c}\Gamma_{i,c}^T\bigl(M_{i,c}R_cC_i^n+\beta_{i,c}\bigr).
 \label{eq:local-lspg}
\end{equation}
With exact pieces, \eqref{eq:local-galerkin} and \eqref{eq:local-lspg} are the multiscale Galerkin and LSPG approximations \eqref{eq:reduced-local} of the finite-volume transport step, so the learned model approximates a projection method rather than the encoded inverse map.
The Galerkin model learns the pieces $X_{i,c}$ and $y_{i,c}$, whose size depends only on the number of local multiscale functions, whereas the LSPG model learns the block rows $B_{i,c}$ and $\beta_{i,c}$, whose size grows with the number of fine cells in the coarse block.

\noindent\textit{Stage 2: Local atoms and regression.}
For each species $i$ and coarse cell $K_c^H$, we collect the entries of the learned pieces into one vector $v_{i,c}$, which contains $X_{i,c}$ and $y_{i,c}$ for the Galerkin model and $B_{i,c}$ and $\beta_{i,c}$ for the LSPG model. 
Pieces from $N_{\mathrm{tr}}$ training media give the mean $\bar v_{i,c}$ and, by singular-value decomposition of the centered snapshots, ${K_{\mathrm a}}$ orthonormal local atoms $\psi_{i,c,k}$.
Each piece is approximated as
\begin{equation}
 \widehat v_{i,c}(\kappa)=\bar v_{i,c}
 +\sum_{k=1}^{{K_{\mathrm a}}}\widehat w_{i,c,k}(\kappa)\psi_{i,c,k},
 \quad
 \widehat w_{i,c}(\kappa)=W_{i,c}x_{\mathrm{ctx},c}(\kappa)+\bar w_{i,c}, 
 \label{eq:local-atoms}
\end{equation}
where the local context $x_{\mathrm{ctx},c}(\kappa)$ contains compressed porosity and flow in $K_c^H$ and its adjacent fine-cell layers.
The ridge matrices $W_{i,c}$ are fitted on the training media, and the ridge penalty is selected on validation media separately for each species and coarse cell.

\noindent\textit{Stage 3: Online transport advancement.}
At each time step, the local contexts are evaluated from the current porosity and flux, the pieces are predicted by \eqref{eq:local-atoms}, and the reduced matrix of \eqref{eq:local-galerkin} or \eqref{eq:local-lspg} is assembled and factorized.
The new concentration is $C_i^*=\max\{P_ia_i^*,0\}$, as for the global inverse atoms.
We note that, compared with the global inverse atoms, the local atoms need only a few atoms per coarse block and local regressions, since their pieces depend only on the local medium, at the price of a block-sparse reduced solve at every step instead of a matrix--vector product.

\subsection{Learning local chemistry}
\label{sec:ml-chemistry}

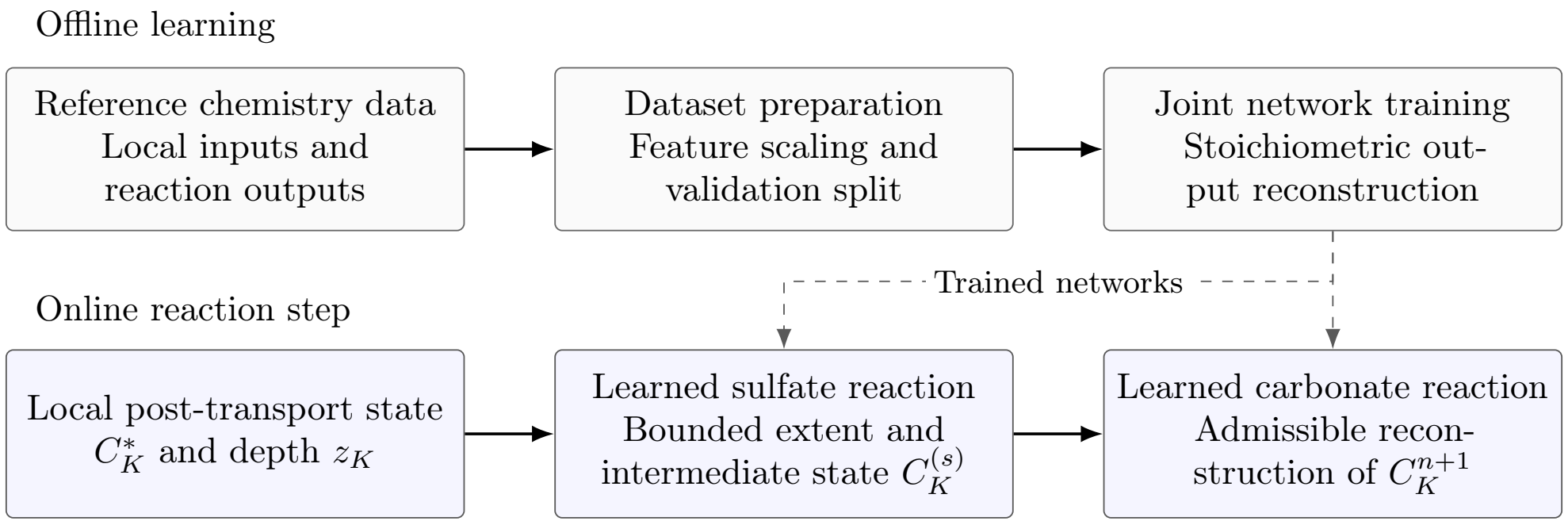
\begin{figure}[htbp]
\centering
\resizebox{0.8\textwidth}{!}{
\begin{tikzpicture}[
  offline/.style={draw=black!60,rounded corners=2pt,fill=gray!4,
    text width=4.25cm,minimum height=1.6cm,align=center,font=\small},
  stage/.style={draw=black!65,rounded corners=2pt,fill=blue!4,
    text width=4.25cm,minimum height=1.65cm,align=center,font=\small},
  arrow/.style={-{Latex[length=2.2mm]},thick},
  trained/.style={-{Latex[length=2mm]},dashed,black!65},
  lab/.style={font=\footnotesize,fill=white,inner sep=2pt}
]
\node[anchor=west,font=\small] at (-2.1,4.) {Offline learning};
\node[offline] (data) at (0,2.8)
  {Reference chemistry data\\Local inputs and reaction outputs};
\node[offline] (prepare) at (5.4,2.8)
  {Dataset preparation\\Feature scaling and validation split};
\node[offline] (train) at (10.8,2.8)
  {Joint network training\\Stoichiometric output reconstruction};
\draw[arrow] (data.east) -- (prepare.west);
\draw[arrow] (prepare.east) -- (train.west);
\node[anchor=west,font=\small] at (-2.1,1.2) {Online reaction step};
\node[stage] (input) at (0,0)
  {Local post-transport state\\$C_K^*$ and depth $z_K$};
\node[stage] (sulfate) at (5.4,0)
  {Learned sulfate reaction\\Bounded extent and intermediate state $C_K^{(s)}$};
\node[stage] (carbonate) at (10.8,0)
  {Learned carbonate reaction\\Admissible reconstruction of $C_K^{n+1}$};
\draw[arrow] (input.east) -- (sulfate.west);
\draw[arrow] (sulfate.east) -- (carbonate.west);
\draw[trained] (train.south) -- (10.8,1.5) -- (5.4,1.5) -- (sulfate.north);
\draw[trained] (train.south) -- (carbonate.north);
\node[lab] at (8.1,1.5) {Trained networks};
\end{tikzpicture}
}
\caption{Local chemistry learning and prediction. Reference pairs train two
sequential reaction maps with stoichiometric and admissibility constraints.
The predicted concentrations enter the porosity update.}
\label{fig:chemistry-learning}
\end{figure}

Operator splitting makes chemistry a cell-local update, evaluated after transport and before the porosity update. Learning this map can reduce reactive-transport costs \cite{delucia2021dectree,leal2020accelerating}.
We encode stoichiometry and admissibility in the output reconstruction, following conservation-aware learning principles \cite{sturm2022conservation}. 
Figure~\ref{fig:chemistry-learning} summarizes reference-data construction, joint training against the reference chemistry output, and local prediction.

For cell $K$, let $C_K^*$ be the post-transport concentrations, in the ordering of \eqref{eq:species}. The
reference map $C_K^{n+1}=\mathcal R_{\Delta t}(C_K^*)$ solves the local problem \eqref{eq:reaction-be}, with depth $z_K$
entering the sulfate-reduction rate. By \eqref{eq:stoichiometric-matrix}, its five-component update is determined by two
reaction extents, $C_K^{n+1}-C_K^*=S_c\xi_c+S_s\xi_s$, where $S_c$ and $S_s$ are the columns of $S_{\mathrm{rxn}}$,
$\xi_s=\Delta t\,R_{\mathrm{sr}}$ is the consumed POC, and $\xi_c=\Delta t\,R_{\mathrm{cal}}$ is the net CaCO$_3$ change.
The surrogate applies sulfate reduction followed by carbonate reaction,
$C_K^*\xrightarrow{\mathcal R_{\theta_s}} C_K^{(s)}\xrightarrow{\mathcal R_{\theta_c}} C_K^{n+1}$, each submap
reconstructing all species from one extent, and both are trained jointly against the reference output. We omit cell
subscripts below.

Reactant availability bounds the sulfate extent by $\xi_{s,\max}=\min(C_2^*,2C_3^*)$. A neural network $N_s$ with parameters $\theta_s$ takes scaled POC, sulfate, and depth as features $\psi_s$ 
and predicts a reaction fraction through a bounded transformation
$\mathcal H_s$ with values in $[0,1]$
\begin{equation*}
  \xi_s=\xi_{s,\max}\,\mathcal H_s\!\left(N_s(\psi_s)\right),
  \quad C^{(s)}=C^*+S_s\xi_s.
\end{equation*}
For carbonate reaction, the network $N_c$ with parameters $\theta_c$ takes the scaled
conserved combinations $I_{\mathrm{Ca}}=C_1^{(s)}+C_4^{(s)}$ and $I_{\mathrm{DIC}}=C_5^{(s)}-C_4^{(s)}$ and the
intermediate Ca and CaCO$_3$ concentrations as features $\psi_c$. The transformation $\mathcal H_c$ gives a positive Ca
candidate, which $\mathcal P_c$ restricts to values consistent with nonnegative CaCO$_3$, Ca, and DIC
\begin{equation*}
    C_4^{n+1}=\mathcal P_c\bigl(\mathcal H_c(N_c(\psi_c))\bigr),\quad
    \xi_c=C_4^{(s)}-C_4^{n+1},\quad C^{n+1}=C^{(s)}+S_c\xi_c.
\end{equation*}
The two bounded extents preserve the reaction invariants $C_1+C_4$ and $2C_3-C_2$ and keep the reconstructed
concentrations nonnegative, provided the post-transport input is admissible.

\subsection{Learning the flow approximation}
\label{sec:ml-flow}

\begin{figure}[htbp]
\centering
\resizebox{0.8\textwidth}{!}{
\begin{tikzpicture}[
  offline/.style={draw=black!60,rounded corners=2pt,fill=gray!4,
    text width=4.25cm,minimum height=1.6cm,align=center,font=\small},
  stage/.style={draw=black!65,rounded corners=2pt,fill=blue!4,
    text width=4.25cm,minimum height=1.65cm,align=center,font=\small},
  arrow/.style={-{Latex[length=2.2mm]},thick},
  trained/.style={-{Latex[length=2mm]},dashed,black!65},
  lab/.style={font=\footnotesize,fill=white,inner sep=2pt}
]
\node[anchor=west,font=\small] at (-2.1,6.05) {Offline learning};
\node[offline] (data) at (0,4.8)
  {Reference local flow data\\Permeability patches and interface traces};
\node[offline] (compress) at (5.4,4.8)
  {Trace compression\\Shared patterns for each interface type};
\node[offline] (train) at (10.8,4.8)
  {MLP/CNN training\\Permeability patches to trace coefficients};
\draw[arrow] (data.east) -- (compress.west);
\draw[arrow] (compress.east) -- (train.west);
\node[anchor=west,font=\small] at (-2.1,3.2) {Once per realization};
\node[stage] (input) at (0,2.2)
  {Initial medium\\Local permeability features};
\node[stage] (predict) at (5.4,2.2)
  {Predict interface traces\\Reconstruct independent interface modes};
\node[stage] (basis) at (10.8,2.2)
  {Build the flux space\\Local numerical extensions\\Fixed basis $\Psi$};
\draw[arrow] (input.east) -- (predict.west);
\draw[arrow] (predict.east) -- (basis.west);
\draw[trained] (train.south) -- (10.8,3.5) -- (5.4,3.5) -- (predict.north);
\node[lab] at (8.1,3.5) {Trained trace models};
\node[anchor=west,font=\small] at (-2.1,0.5) {At each time step};
\node[stage] (current) at (0,-0.5)
  {Evolving medium\\Permeability and porosity-rate source};
\node[stage] (solve) at (5.4,-0.5)
  {Mixed coarse flow solve\\Update coefficients and impose block balance};
\node[stage] (output) at (10.8,-0.5)
  {Reconstruct Darcy fluxes\\$q^n=\Psi a_q^n$\\Supply flow to transport};
\draw[arrow] (current.east) -- (solve.west);
\draw[arrow] (solve.east) -- (output.west);
\draw[trained] (basis.south) -- (10.8,0.8) -- (5.4,0.8) -- (solve.north);
\draw[trained] (basis.south) -- (output.north);
\node[lab] at (8.1,0.8) {Initialized flux space};
\end{tikzpicture}
}
\caption{Flow learning and prediction. Compressed reference traces train
MLP or CNN predictors. Their outputs define a flux space once per realization.
A mixed coarse solve uses the evolving medium at every time step and
reconstructs the Darcy fluxes for transport.}
\label{fig:flow-training-overview}
\end{figure}
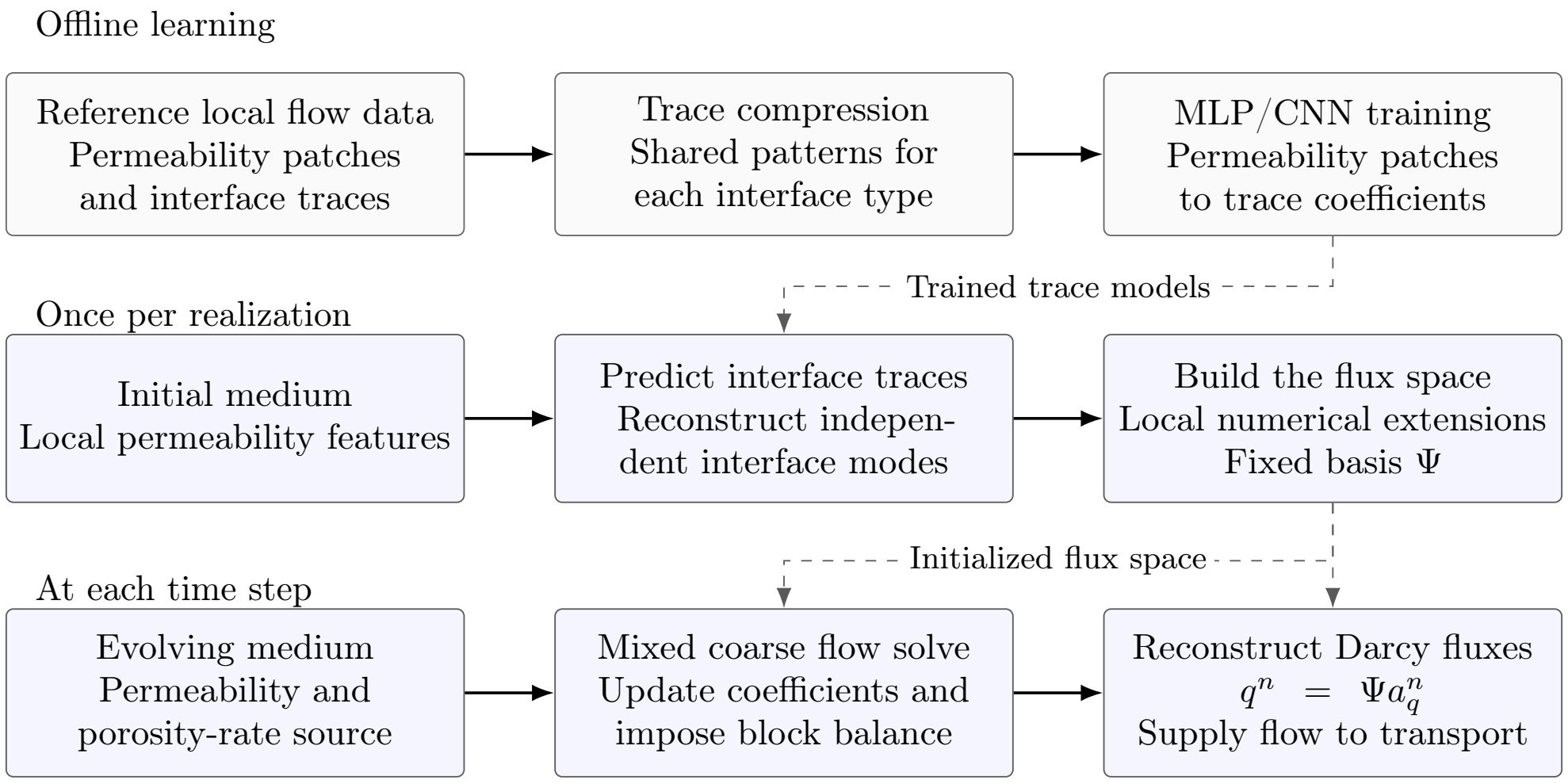

We approximate the Darcy flux in a mixed multiscale space \cite{chen2003mixed,efendiev2009multiscale,chung2016mixed}. 
The flow uses the same $N_H$ coarse cells as the transport multiscale space of Subsection~\ref{sec:ms-space}, which we call coarse blocks $B$ here, but its basis functions are associated with the coarse interfaces instead of the coarse nodes.
Coarse interfaces $e$ separate neighboring blocks or coincide with the pressure boundaries, and a trace on $e$ is the normal flux integrated over its fine-grid faces. 
Local Neumann solves extend each interface trace into the two adjacent blocks, and the extended modes form the columns of the flux basis $\Psi\in\mathbb R^{N_F\times m_q}$, where $N_F$ is the number of active fine faces and $m_q$ the number of flux modes. 
The Darcy flux is represented by $q^n=\Psi a_q^n$, where $q^n$ contains oriented integrated fluxes and $a_q^n$ their coarse coefficients, and a coarse mixed solve determines $a_q^n$ with exact mass balance in each block. 
Standard mixed multiscale methods compute the traces by local fine-scale problems for each new permeability field, and we learn them from the permeability instead.
The flow approximation is constructed in an offline stage, an initialization stage for each realization, and an online stage in every time step, summarized in Figure~\ref{fig:flow-training-overview}.

\noindent\textit{Stage 1: Offline trace construction and learning.} 
Local finite-volume problems on oversampled patches, with linear pressure data in each coordinate direction, generate reference traces.
After normalization, the directional traces are concatenated into a vector $y_e$. 
The input $x_e$ is standardized log permeability on the surrounding patch. 
For each type, we subtract the training mean $\overline y_{\iota}$ and apply POD to obtain $K_f$ orthonormal trace basis vectors $a_{\iota,\ell}$.
Each basis vector represents a recurring directional-flux pattern over the interface's fine faces. 
The network $N_{f,\iota(e)}$ predicts their coefficients $\alpha_e=N_{f,\iota(e)}(x_e)$ from the permeability patch,
and the predicted trace is $\widehat y_e=\overline y_{\iota(e)}+\sum_{\ell=1}^{K_f}\alpha_{e,\ell}a_{\iota(e),\ell}$.
We use multilayer perceptron (MLP) predictors in 2D and convolutional neural network (CNN) predictors in 3D.

\noindent\textit{Stage 2: Initialization for each realization.}
The trained networks predict traces from the initial permeability.
Extracting their independent modes gives interface modes, which local Neumann solves extend into adjacent blocks using the initial permeability.
The extended modes form the flux basis $\Psi$, which is fixed for the realization. 
We also preconstruct $\mathsf B\in\mathbb R^{N_B\times m_q}$, which integrates each mode's divergence over the $N_B$ coarse blocks. 

\noindent\textit{Stage 3: Online flow advancement.}
At each time step, the current transmissibilities $T_f^n$ and porosity rate define the reduced matrix $M_q^n=\Psi^T\operatorname{diag}\bigl((T_f^n)^{-1}\bigr)\Psi$ and the block source $S^n=(S_B^n)_B$ with $S_B^n=-\sum_{K\subset B}|K|(\phi_K^n-\phi_K^{n-1})/\Delta t$, together with the boundary-pressure vector $g_D^n$. 
Then flux coefficients and block pressures satisfy
\begin{equation}
 \begin{pmatrix}M_q^n&-\mathsf B^T\\ \mathsf B&0\end{pmatrix}
 \begin{pmatrix}a_q^n\\P^n\end{pmatrix}
 =\begin{pmatrix}-\Psi^Tg_D^n\\S^n\end{pmatrix}.
 \label{eq:flow-reduced-system}
\end{equation}
where the reduced matrix $M_q^n$ is  the projection of the current fine-face inverse transmissibilities onto the fixed flux space. Each mode is supported on the two coarse blocks adjacent to its interface, so $M_q^n$ is sparse with a coarse-block stencil, and its entries $(M_q^n)_{ab}=\sum_f\Psi_{fa}\Psi_{fb}/T_f^n$ are evaluated by one precomputed sparse map applied to $(T_f^n)^{-1}$.
Each time step therefore requires only this evaluation, which needs the fine-face transmissibilities computed from $k^n$ but no fine-scale solve or fine flow matrix, one reduced saddle-point solve \eqref{eq:flow-reduced-system}, and the flux reconstruction $q^n=\Psi a_q^n$. The learned traces enter only through $\Psi$, which is fixed at initialization, while $M_q^n$ follows the evolving permeability. Local fine-scale solves are needed only for initialization.

\section{Stability and Error Propagation}
\label{sec:learned-error-propagation}

We compare the learned trajectory with the reference split solver of Section~\ref{sec:numerics}, which uses the same stage order, boundary data, initial medium, and porosity--permeability update.
Following the theory of additive (splitting) schemes \cite{samarskii2001additive,samarskii2002operator,vabishchevich2013additive}, we write each stage in canonical form, estimate it from one time layer to the next, and chain the estimates over the split step and the time steps, so that the stability of the reference scheme and the learned-stage errors give the error bound.
Section~\ref{sec:stability} proves the stability of the split scheme, and Section~\ref{sec:error-propagation} extends it to the learned transport and bounds the propagation of the learned-stage errors.

The state at time $t_n$ is $U^n=(C^n,\phi^n,\phi^{n-1})$, with the porosity increment $\Delta\phi^n=\phi^n-\phi^{n-1}$. 
For a general state $U$ of the same form, $C_U$, $\phi_U$, and $\Delta\phi_U$ denote its concentrations, porosity, and porosity increment. 
One reference split step $U^{n+1}=\mathcal S_{\Delta t}U^n$ solves 
the flow problem \eqref{eq:fv-darcy-flux}--\eqref{eq:fv-flow-matrix} for $q^n$, 
the transport step \eqref{eq:operator-form} for $C^*$, 
the reaction step \eqref{eq:reaction-disc} for $C^{n+1}$, and it updates the porosity {to $\phi^{n+1}$ and the permeability to $k^{n+1}$ by \eqref{eq:permeability-update}}. The reference flow map $\mathcal F(U)$ returns the Darcy fluxes, the transport map $\mathcal T_{\Delta t}(C,U,q)$ returns the solution $C^*$ of the transport equations with initial concentrations $C$, the coefficients of $U$, and the fluxes $q$, and the reaction map $\mathcal R_{\Delta t}(C^*)$ returns the solution $C^{n+1}$ of the reaction equations in each cell, so that $q^n=\mathcal F(U^n)$, the post-transport concentrations are $C^{*,n}=\mathcal T_{\Delta t}(C^n,U^n,q^n)$, and $C^{n+1}=\mathcal R_{\Delta t}(C^{*,n})$.
The learned step $\widehat{\mathcal S}_{\Delta t}$ replaces $\mathcal F$, $\mathcal T_{\Delta t}$, and $\mathcal R_{\Delta t}$ by their learned approximations $\widehat{\mathcal F}$, $\widehat{\mathcal T}_{\Delta t}$, and $\widehat{\mathcal R}_{\Delta t}$, 
so that the learned trajectory is $\widehat U^{n+1}=\widehat{\mathcal S}_{\Delta t}\widehat U^n$ with $\widehat U^0=U^0$.

\subsection{Stability of the split scheme}
\label{sec:stability}

With $C=(C_1,\ldots,C_5)$ and $M(\phi)=\operatorname{diag}(M_1(\phi),\ldots,M_5(\phi))$ from \eqref{eq:operator-form}, 
we measure the concentrations in the phase norm 
$\|C\|_\phi^2=(M(\phi)C,C)=\sum_i\|C_i\|_\phi^2$, 
with $\|C_i\|_\phi^2=(M_i(\phi)C_i,C_i)=\sum_K|K|\,\eta_i(\phi_K)\,C_{i,K}^2$, the discrete counterpart of $\sum_i\int_\Omega\eta_i(\phi)c_i^2\,dx$, in which each species is weighted by the volume of its phase.
Its cell part is $|y|_{K,\phi}^2=\sum_i\eta_i(\phi_K)y_i^2$, so that $\|C\|_\phi^2=\sum_K|K|\,|C_K|_{K,\phi}^2$.
Cell fields, such as the porosity, are measured in the discrete $L^2$ norm $\|u\|^2=\sum_K|K|\,u_K^2$, and the Darcy fluxes \eqref{eq:fv-darcy-flux} by $\|q\|^2=\sum_eq_e^2$.
Since $\eta_i\leq1$,
\begin{equation}
\|C\|_\phi\leq\|C\|\leq \underline{\eta}(\phi)^{-1/2}\|C\|_\phi,\quad \underline{\eta}(\phi) =\min_{i,K}\eta_i(\phi_K)>0.
 \label{eq:norm-equivalence}
\end{equation}
We use two assumptions on the reference trajectory, for $0\leq n<N_t$ with $N_t={t_{\max}}/\Delta t$.
\begin{itemize}[leftmargin=*,itemsep=1pt]
\item[(A1)] (slowly changing porosity) $|\phi_K^{n+1}-\phi_K^n|\leq2\mu_{\phi}\Delta t\,\eta_i(\phi_K^n)$ for all $i$ and $K$, with $\mu_{\phi}\Delta t\leq1/2$ and $\phi^{-1}=\phi^0$; by the porosity update $\phi^{n+1}-\phi^n=-\gamma_\phi\Delta t\,R_1(C^{n+1})$, this holds with $\mu_{\phi}=\gamma_\phi\max|R_1|/(2\eta_{\min})$
{, where $\eta_{\min}=\min_{n} \underline{\eta}(\phi^n)>0$.}
\item[(A2)] (reaction) $|R(x)-R(y)|_{{K,\phi^n}}\leq L_r\,|x-y|_{{K,\phi^n}}$ for all $x,y$ in a set $\mathcal C\subset\R^5$ and all $K$, with $L_r\Delta t\leq1/2$.
\end{itemize}
Here $2\mu_{\phi}$ is the largest relative rate of change of the pore and solid volumes, while $1/L_r$ is the fastest reaction time scale.
Both assumptions thus bound the time step by a physical time scale.

The proof of stability chains one-step estimates of the stages, which are stated and proved in Appendix~\ref{app:stage-estimates}.
In the reference medium, one transport step increases the phase norm of the difference of two solutions at most by the factor $1+2\mu_{\phi}\Delta t$ (Lemma~\ref{lem:transport-stage}), and one reaction step at most by the factor $1+2L_r\Delta t$ (Lemma~\ref{lem:reaction-stage}).
The CaCO$_3$ precipitated or dissolved in one step, $m(C)=\Delta t\,R_1(\mathcal R_{\Delta t}(C))$, which gives the porosity increment $\Delta\phi^{n+1}=-\gamma_\phi\,m(C^{*,n})$, changes between two inputs at most by $2\eta_{\min}^{-1/2}L_r\Delta t$ times their difference (Lemma~\ref{lem:reaction-stage}), and the change of the phase volumes from step $n$ to step $n+1$ increases the phase norm at most by the factor $1+\mu_{\phi}\Delta t$ (Lemma~\ref{lem:phase-change}).
The stages are coupled through the medium, which enters the flow through the permeability and the porosity-change source, and the transport through the storage, the diffusion coefficients, and the Darcy fluxes.
For the transport of the same concentrations $C_U$ in the medium of a state $U$ and in the reference medium, Lemma~\ref{lem:medium} gives
\begin{equation}
 \|\mathcal T_{\Delta t}(C_U,U,\mathcal F(U))-\mathcal T_{\Delta t}(C_U,U^n,q^n)\|_{\phi^n}\leq\Lambda_\phi\Delta t\,\|\phi_U-\phi^n\|+\Lambda_{\Delta\phi}\|\Delta\phi_U-\Delta\phi^n\|,
 \label{eq:medium-bound}
\end{equation}
where the constants $\Lambda_\phi$ and $\Lambda_{\Delta\phi}$ depend on the bounds of the transmissibilities, of the reference pressure gradient, and of the transported concentrations and their discrete gradients, and not on $\Delta t$.

\begin{theorem}[Stability of the split scheme]
\label{thm:split-stability}
Let (A1) and (A2) hold, and let $\mathcal U$ be a set of states $U$ for which \eqref{eq:medium-bound} holds at every step with the same constants $\Lambda_\phi$ and $\Lambda_{\Delta\phi}$ and the cell values of $\mathcal S_{\Delta t}U$ and of the reference trajectory lie in $\mathcal C$.
Then the distance of two states,
\begin{equation*}
 \mathcal E_n(U,V)=\|C_U-C_V\|_{\phi^n}+\|\phi_U-\phi_V\|+\Lambda_{\Delta\phi}\|\Delta\phi_U-\Delta\phi_V\|,
\end{equation*}
satisfies, for all $U\in\mathcal U$,
\begin{equation}
 \mathcal E_{n+1}(\mathcal S_{\Delta t}U,U^{n+1})\leq e^{\beta\Delta t}\,\mathcal E_n(U,U^n),
 \quad
\beta=\max(2\mu_{\phi},\Lambda_\phi)+\mu_{\phi}+2L_r\bigl(1+\gamma_\phi\eta_{\min}^{-1/2}(1+\Lambda_{\Delta\phi})\bigr).
 \label{eq:split-stability}
\end{equation}
In particular, if a perturbed reference trajectory $\widetilde U^{n+1}=\mathcal S_{\Delta t}\widetilde U^n$ stays in $\mathcal U$, then $\mathcal E_n(\widetilde U^n,U^n)\leq e^{\beta {t_{\max}}}\mathcal E_0(\widetilde U^0,U^0)$ for $n\leq N_t$.
\end{theorem}

\begin{proof}
Let $U\in\mathcal U$, $V=\mathcal S_{\Delta t}U$, and $C_U^*=\mathcal T_{\Delta t}(C_U,U,\mathcal F(U))$, and let $z_C=\|C_U-C^n\|_{{\phi^n}}$, $z_\phi=\|\phi_U-\phi^n\|$, $z_{\Delta\phi}=\|\Delta\phi_U-\Delta\phi^n\|$, and $z_*=\|C_U^*-C^{*,n}\|_{{\phi^n}}$.

\noindent\textit{Transport.}
We pass through the transport of $C_U$ in the reference medium,
\begin{equation*}
 C_U^*-C^{*,n}
 =\bigl[C_U^*-\mathcal T_{\Delta t}(C_U,U^n,q^n)\bigr]
 +\bigl[\mathcal T_{\Delta t}(C_U,U^n,q^n)-\mathcal T_{\Delta t}(C^n,U^n,q^n)\bigr].
\end{equation*}
The first difference has the same concentrations in different media and is bounded by \eqref{eq:medium-bound}.
The second difference $w$ has different concentrations in the same medium, so each $w_i$ solves the transport step \eqref{eq:operator-form} with the right-hand side $M_i(\phi^n)(C_{U,i}-C_i^n)/\Delta t$, and Lemma~\ref{lem:transport-stage} with $v=C_{U,i}-C_i^n$ and $g=0$, summed over the species, gives $\|w\|_{{\phi^n}}\leq(1+2\mu_{\phi}\Delta t)z_C$.
Hence
\begin{equation*}
 z_*\leq(1+2\mu_{\phi}\Delta t)z_C+\Lambda_\phi\Delta t\,z_\phi+\Lambda_{\Delta\phi}z_{\Delta\phi}.
\end{equation*}
\textit{Reaction and porosity.}
Lemmas~\ref{lem:phase-change} and~\ref{lem:reaction-stage} give $\|C_V-C^{n+1}\|_{{\phi^{n+1}}}\leq(1+\mu_{\phi}\Delta t)(1+2L_r\Delta t)\,z_*$.
The increments are $\Delta\phi_V=-\gamma_\phi m(C_U^*)$ and $\Delta\phi^{n+1}=-\gamma_\phi m(C^{*,n})$, so Lemma~\ref{lem:reaction-stage} gives $\|\Delta\phi_V-\Delta\phi^{n+1}\|\leq2\gamma_\phi\eta_{\min}^{-1/2}L_r\Delta t\,z_*$, and $\phi_V=\phi_U+\Delta\phi_V$ and $\phi^{n+1}=\phi^n+\Delta\phi^{n+1}$ give $\|\phi_V-\phi^{n+1}\|\leq z_\phi+\|\Delta\phi_V-\Delta\phi^{n+1}\|$.

\noindent\textit{One split step.}
Adding the three bounds with the weights $1$, $1$, and $\Lambda_{\Delta\phi}$ gives
\begin{equation*}
 \mathcal E_{n+1}(V,U^{n+1})\leq z_\phi+\alpha\,z_*,\quad
 \alpha=(1+\mu_{\phi}\Delta t)(1+2L_r\Delta t)+2\gamma_\phi\eta_{\min}^{-1/2}L_r(1+\Lambda_{\Delta\phi})\Delta t\leq e^{(\beta-\max(2\mu_{\phi},\Lambda_\phi))\Delta t}.
\end{equation*}
After inserting the bound for $z_*$, the coefficients of $z_C$, $z_\phi$, and $\Lambda_{\Delta\phi}z_{\Delta\phi}$ are $\alpha(1+2\mu_{\phi}\Delta t)$, $1+\alpha\Lambda_\phi\Delta t$, and $\alpha$, which are all at most $\alpha\bigl(1+\max(2\mu_{\phi},\Lambda_\phi)\Delta t\bigr)\leq e^{\beta\Delta t}$.
This gives \eqref{eq:split-stability}, and applying it $n$ times with $U=\widetilde U^j$ gives the bound for the perturbed trajectory, since $n\Delta t\leq {t_{\max}}$.
\end{proof}

Theorem~\ref{thm:split-stability} bounds the transition operator of the whole split step \cite{vabishchevich2012transition,vabishchevich2014systems}.
The exponent $\beta$ is a sum of rates, in units of 1/time, of the porosity change ($\mu_{\phi}$), of the response of the transport to the medium ($\Lambda_\phi$), and of the reactions ($L_r$), and it does not depend on $\Delta t$, although the flow source contains $\Delta\phi/\Delta t$.
The reason is that the reaction changes the increment only by $\Delta t$ times the post-transport difference (Lemma~\ref{lem:reaction-stage}), which the weight $\Lambda_{\Delta\phi}$ in $\mathcal E_n$ absorbs.
In a fixed medium, $\gamma_\phi=0$, the porosity terms vanish and the factor of one step is $1+2L_r\Delta t$.

\subsection{Error of the multiscale and learned schemes}
\label{sec:error-propagation}

We compare three trajectories that start from the same state $U^0$, 
the reference trajectory $U^{n+1}=\mathcal S_{\Delta t}U^n$, 
the multiscale trajectory $U_{\mathrm{ms}}^{n+1}=\mathcal S^{\mathrm{ms}}_{\Delta t}U_{\mathrm{ms}}^n$, whose step replaces only the transport by the multiscale step $\mathcal T^{\mathrm{ms}}_{\Delta t}$, which solves the Galerkin or LSPG system \eqref{eq:reduced-local} assembled from the fine-grid operator, and 
the learned trajectory $\widehat U^{n+1}=\widehat{\mathcal S}_{\Delta t}\widehat U^n$.
Let $U'$ be a multiscale or learned state at step $n$ and $\mathcal S'_{\Delta t}$ its step.
The distance $\mathcal E_{n+1}$ of Theorem~\ref{thm:split-stability} 
satisfies the triangle inequality, so inserting the reference step from the approximate state gives
\begin{equation}
 \mathcal E_{n+1}(\mathcal S'_{\Delta t}U',U^{n+1})\leq
 \underbrace{\mathcal E_{n+1}(\mathcal S'_{\Delta t}U',\mathcal S_{\Delta t}U')}_{\text{local error}}
 +\underbrace{\mathcal E_{n+1}(\mathcal S_{\Delta t}U',U^{n+1})}_{\text{propagated error}}.
 \label{eq:error-splitting}
\end{equation}
Here Theorem~\ref{thm:split-stability} with $U=U'$ bounds the propagated error by $e^{\beta\Delta t}\mathcal E_n(U',U^n)$. 
It uses the stability of the reference step and not of the learned step, so the multiscale and learned stages enter only through the local error, the error of one step from the same state $U'$.

The local error consists of four one-step errors, which are bounded in Appendix~\ref{app:local-errors}.
The Galerkin multiscale step is the projection of the fine-grid system \eqref{eq:training-system} onto the multiscale space $V_i^{\mathrm{ms}}$ of Section~\ref{sec:ms-space}, so it is stable in the phase norm with the same factor $1+2\mu_\phi \Delta t$ as the fine-grid step, and its error $\delta_{\mathrm{ms}}^n$ is bounded by the projection error of the fine-grid solution in the mass and energy norms (Lemmas~\ref{lem:ms-projection} and~\ref{lem:ms-step}),
{and for the LSPG system in the mass norm only.}
As for multiscale spaces built from local snapshots \cite{efendiev2013gmsfem,calo2016randomized}, the projection error in each coarse neighborhood $\omega_j$ is the sum of the POD truncation error, given by the first discarded singular value $\sigma_{i,j,m+1}$ of the snapshot matrix $X_{i,j}$, and of the snapshot error, which measures how well the snapshots represent the test media.
The learned step of the Galerkin local atoms solves \eqref{eq:local-galerkin} with the predicted transport matrix and boundary vector, with the errors $\varepsilon_A^n$ and $\varepsilon_b^n$ relative to the storage.
It is stable with the factor $1+2( \mu_\phi + \varepsilon_A^n)\Delta t$, and its distance $\delta_\ell^n$ from the multiscale step is proportional to $\Delta t$ times these errors (Lemma~\ref{lem:learned-step}).
{The learned step of the global inverse atoms differs from the encoded fine-grid step by the error of the predicted map $(\widehat J_i,\widehat g_i)$ applied to the current coefficients, without a reduced solve, and the error of the encoded step is the least-squares projection error of the fine-grid solution and of the current state.}
For both atom models, the error of the predicted operator is the sum of the atom truncation error, which has the same two terms as the projection error of the multiscale space, and of the regression error of the atom coefficients (Lemma~\ref{lem:atoms}).
The learned flow is the mixed multiscale solve \eqref{eq:flow-reduced-system} with the flux basis $\Psi$ of Section~\ref{sec:ml-flow}, and its flux error $\varepsilon_f^n$ is the best approximation error of the fine-grid fluxes in $\operatorname{range}\Psi$ plus the oscillation of the porosity-change source inside the coarse blocks \eqref{eq:flow-error}.
For the learned chemistry map of Section~\ref{sec:ml-chemistry}, the one-step error $\delta_r^n$ is assumed to be given.

\begin{theorem}[Error of the multiscale and learned schemes]
\label{thm:learned-errors}
At the learned state $\widehat U^n$, let $\widehat q^n=\widehat{\mathcal F}(\widehat U^n)$ and $\widetilde q^n=\mathcal F(\widehat U^n)$ be the learned and the fine-grid fluxes, and let
\begin{equation*}
\begin{aligned}
 \widehat C^{*,n}&=\widehat{\mathcal T}_{\Delta t}(\widehat C^n,\widehat U^n,\widehat q^n),&\quad
 \widehat C_{\mathrm{ms}}^{*,n}&=\mathcal T^{\mathrm{ms}}_{\Delta t}(\widehat C^n,\widehat U^n,\widehat q^n),\\
 \widetilde C_q^{*,n}&=\mathcal T_{\Delta t}(\widehat C^n,\widehat U^n,\widehat q^n),&\quad
 \widetilde C^{*,n}&=\mathcal T_{\Delta t}(\widehat C^n,\widehat U^n,\widetilde q^n)
\end{aligned}
\end{equation*}
be the learned, the multiscale, and the fine-grid transport steps with the learned fluxes, and the fine-grid transport step with the fine-grid fluxes.
Assume that, for $0\leq n<N_t$,
\begin{itemize}[leftmargin=*,itemsep=1pt]
\item[(i)] the assumptions of Theorem~\ref{thm:split-stability} hold with a set $\mathcal U$ that contains $\widehat U^n$,
\item[(ii)] the one-step errors of the multiscale transport, the learned transport, the learned flow, and the learned chemistry satisfy
\begin{equation*}
 \|\widehat C_{\mathrm{ms}}^{*,n}-\widetilde C_q^{*,n}\|\leq\delta_{\mathrm{ms}}^n,\quad
 \|\widehat C^{*,n}-\widehat C_{\mathrm{ms}}^{*,n}\|\leq\delta_\ell^n,\quad
 \|\widehat q^n-\widetilde q^n\|\leq\varepsilon_f^n,\quad
 \|\widehat{\mathcal R}_{\Delta t}(\widehat C^{*,n})-\mathcal R_{\Delta t}(\widehat C^{*,n})\|\leq\delta_r^n,
\end{equation*}
\item[(iii)] the transport depends on the fluxes as $\|\widetilde C_q^{*,n}-\widetilde C^{*,n}\|_{{\phi^n}}\leq\Delta t\,L_q\|\widehat q^n-\widetilde q^n\|$, with $L_q$ of Lemma~\ref{lem:medium},
\item[(iv)] the cell values of $\mathcal R_{\Delta t}(\widehat C^{*,n})$ and $\mathcal R_{\Delta t}(\widetilde C^{*,n})$ lie in $\mathcal C$.
\end{itemize}
Then, with $\beta$ of Theorem~\ref{thm:split-stability} and $\varkappa=2+\gamma_\phi\eta_{\min}^{-1/2}(1+\Lambda_{\Delta\phi})$,
\begin{equation}
 \max_{n\leq N_t}\mathcal E_n(\widehat U^n,U^n)\leq\varkappa\,e^{\beta t_{\max}}\sum_{j=0}^{N_t-1}\delta^j,\quad
 \delta^j=\delta_{\mathrm{ms}}^j+\delta_\ell^j+\delta_r^j+\Delta t\,L_q\varepsilon_f^j,
 \label{eq:learned-global-error}
\end{equation}
and $\|\widehat C^n-C^n\|\leq\eta_{\min}^{-1/2}\mathcal E_n(\widehat U^n,U^n)$ in the $L^2$ norm.
For the multiscale solver, with the reference flow and chemistry, the same bound holds for $\mathcal E_n(U_{\mathrm{ms}}^n,U^n)$ with $\delta^j=\delta_{\mathrm{ms}}^j$.
\end{theorem}

\begin{proof}
Let $E^n=\mathcal E_n(\widehat U^n,U^n)$, and let $V=\widehat{\mathcal S}_{\Delta t}\widehat U^n=\widehat U^{n+1}$ and $W=\mathcal S_{\Delta t}\widehat U^n$ be the learned and the reference steps from $\widehat U^n$.
By \eqref{eq:error-splitting} and Theorem~\ref{thm:split-stability} with (i), $E^{n+1}\leq\mathcal E_{n+1}(V,W)+e^{\beta\Delta t}E^n$, and it remains to bound the local error $\mathcal E_{n+1}(V,W)$.

\noindent\textit{Transport.}
The post-transport concentrations of $V$ and $W$ are $\widehat C^{*,n}$ and $\widetilde C^{*,n}$.
Passing through $\widehat C_{\mathrm{ms}}^{*,n}$ and $\widetilde C_q^{*,n}$, (ii), (iii), and $\|\cdot\|_{{\phi^n}}\leq\|\cdot\|$ give $\|\widehat C^{*,n}-\widetilde C^{*,n}\|_{{\phi^n}}\leq\delta_\ell^n+\delta_{\mathrm{ms}}^n+\Delta t\,L_q\varepsilon_f^n=d^n$.

\noindent\textit{Reaction and porosity.}
The concentrations are $C_V=\widehat{\mathcal R}_{\Delta t}(\widehat C^{*,n})$ and $C_W=\mathcal R_{\Delta t}(\widetilde C^{*,n})$.
Adding and subtracting $\mathcal R_{\Delta t}(\widehat C^{*,n})$, Lemma~\ref{lem:phase-change}, (ii), and Lemma~\ref{lem:reaction-stage} with (iv) give
\begin{equation*}
 \|C_V-C_W\|_{{\phi^{n+1}}}\leq(1+\mu_\phi \Delta t)\bigl(\delta_r^n+(1+2L_r\Delta t)\,d^n\bigr).
\end{equation*}
Both steps start from $\widehat\phi^n$ and update the porosity by the reaction step,
 $\Delta\phi_V=-\gamma_\phi[\widehat{\mathcal R}_{\Delta t}(\widehat C^{*,n})-\widehat C^{*,n}]_1$ and 
$\Delta\phi_W=-\gamma_\phi [\mathcal R_{\Delta t}(\widetilde C^{*,n})-\widetilde C^{*,n}]_1$, which is $-\gamma_\phi\,m(\widetilde C^{*,n})$ of Lemma~\ref{lem:reaction-stage} by \eqref{eq:reaction-disc}, so that, by (ii) and Lemma~\ref{lem:reaction-stage},
\begin{equation*}
 \|\phi_V-\phi_W\|=\|\Delta\phi_V-\Delta\phi_W\|\leq\gamma_\phi\bigl(\delta_r^n+2\eta_{\min}^{-1/2}L_r\Delta t\,d^n\bigr).
\end{equation*}
\textit{Summation.}
Adding the three parts of $\mathcal E_{n+1}(V,W)$ with the weights $1$, $1$, and $\Lambda_{\Delta\phi}$, the coefficients of $\delta_r^n$ and $d^n$ are
$1+\mu_\phi\Delta t+\gamma_\phi(1+\Lambda_{\Delta\phi})$ and $(1+\mu_\phi\Delta t)(1+2L_r\Delta t)+2\gamma_\phi\eta_{\min}^{-1/2}L_r\Delta t\,(1+\Lambda_{\Delta\phi})$, which are both at most $e^{\mu_\phi\Delta t}\varkappa$ by $2L_r\Delta t\leq1$ and $\eta_{\min}^{-1/2}\geq1$.
Hence $E^{n+1}\leq e^{\beta\Delta t}E^n+e^{\mu_\phi\Delta t}\varkappa\,\delta^n$ with $E^0=0$, and, since $\mu_\phi\leq\beta$ and $n\Delta t\leq t_{\max}$,
\begin{equation*}
 E^n\leq\varkappa\sum_{j=0}^{n-1}e^{\beta(n-1-j)\Delta t+\mu_\phi\Delta t}\,\delta^j\leq\varkappa\,e^{\beta t_{\max}}\sum_{j=0}^{N_t-1}\delta^j.
\end{equation*}
The $L^2$ bound follows from \eqref{eq:norm-equivalence}, and the multiscale solver is the case $\delta_\ell^n=\delta_r^n=\varepsilon_f^n=0$ with $U_{\mathrm{ms}}^n$ in place of $\widehat U^n$.
\end{proof}

The factor $\varkappa$ counts how one local error enters the concentrations and, through the change of CaCO$_3$, the porosity and its increment, and $e^{\beta t_{\max}}$ is the largest growth of these errors over the simulated time.
Since $\|\widehat C^n\|_{{\phi^n}}\leq\|C^n\|_{{\phi^n}}+\mathcal E_n(\widehat U^n,U^n)$, the learned scheme inherits the stability of the reference scheme up to the accumulated one-step errors.
The errors $\delta_{\mathrm{ms}}^n$ and $\delta_r^n$ are errors of complete one-step maps, so a uniform bound $\delta$ on them contributes $\mathcal O(t_{\max}\delta/\Delta t)$.
The learning error of the Galerkin local atoms and the flux error enter with the factor $\Delta t$, so they contribute $\mathcal O(t_{\max})$ times the largest operator and flux errors for every $\Delta t$.
{For the LSPG local atoms, which learn the block rows $B_{i,c}$ and $\beta_{i,c}$ instead of the reduced transport matrix, $\delta_\ell^n$ is also proportional to $\Delta t$, with one more term that contains the residual of the learned step (Appendix~\ref{app:local-errors}), and for the global inverse atoms it is bounded by the error of the predicted map \eqref{eq:delta-ell-global}, without a  factor $\Delta t$.}

\section{Numerical Results}
\label{sec:results}

We consider 2D and 3D test problems and report results for a representative test realization, whose full-order solution is shown first, and for an evaluation set of 20 realizations in each dimension.
The 2D domain is $10\times100$ m with $10\times100$ cells, and the 3D domain is $10\times10\times100$ m with $10\times10\times100$ cells.
Permeability and initial porosity share the Gaussian field of \eqref{eq:random-fields}, with $k_0=10^{-13}\,\mathrm{m}^2$, $N_{\mathrm{KL}}=100$ retained KLE modes, $\sigma_k=1.5$, $\bar\phi=0.45$, $\sigma_\phi=0.05$, and correlation lengths of 3 m in depth and 8 m horizontally, so the fields are layered.
The porosity and permeability updates use $\gamma_\phi=10^{-3}$, $m_k=3$, and $n_k=2$, and we take $N_t=100$ time steps of $\Delta t=50$ years, ${t_{\max}}=5{,}000$ years.
The evaluation realizations are generated in the same way as the test realization, and porosity and permeability evolve in all coupled runs.
Before local equilibration, the initial concentrations are $c_{2,0}=0.2$, $c_{5,0}=3$, and $c_{1,0}=c_{3,0}=0$, while $c_{4,0}$ decreases linearly from 10 at the surface to 0 at the base \cite{vasilyeva2024decoupled}.
The prescribed concentration boundary values are $c_{2,D}=0.2$ and $c_{3,D}=26.2$ at the top, and $c_{3,D}=c_{4,D}=0$ and $c_{5,D}=20$ at the base.
At the top, burial carries CaCO$_3$, POC, Ca, and DIC into the domain, and their inflow concentrations are set to the mean equilibrated initial concentrations in the top cells.
For the Darcy solve, we use the scaled pressure $\widetilde p=(k_0/\mu)(p-p_{\mathrm{top}})$, where $p_{\mathrm{top}}$ is the reference pressure at the surface.
The boundary conditions are $\widetilde p=0$ at the top and $\widetilde p=u_0H$ at the base, with column depth $H$ and $u_0=1\,\mathrm{cm\,yr^{-1}}$.
The zero boundary value denotes a pressure reference, while the physical pressure difference is $\mu u_0H/k_0$.
Lateral boundaries have no flow.
The prescribed deep burial speed is $0.1\,\mathrm{cm\,yr^{-1}}$.

\begin{figure}[h!]
\centering
\begin{minipage}[b]{0.49\textwidth}\centering
\includegraphics[width=\linewidth]{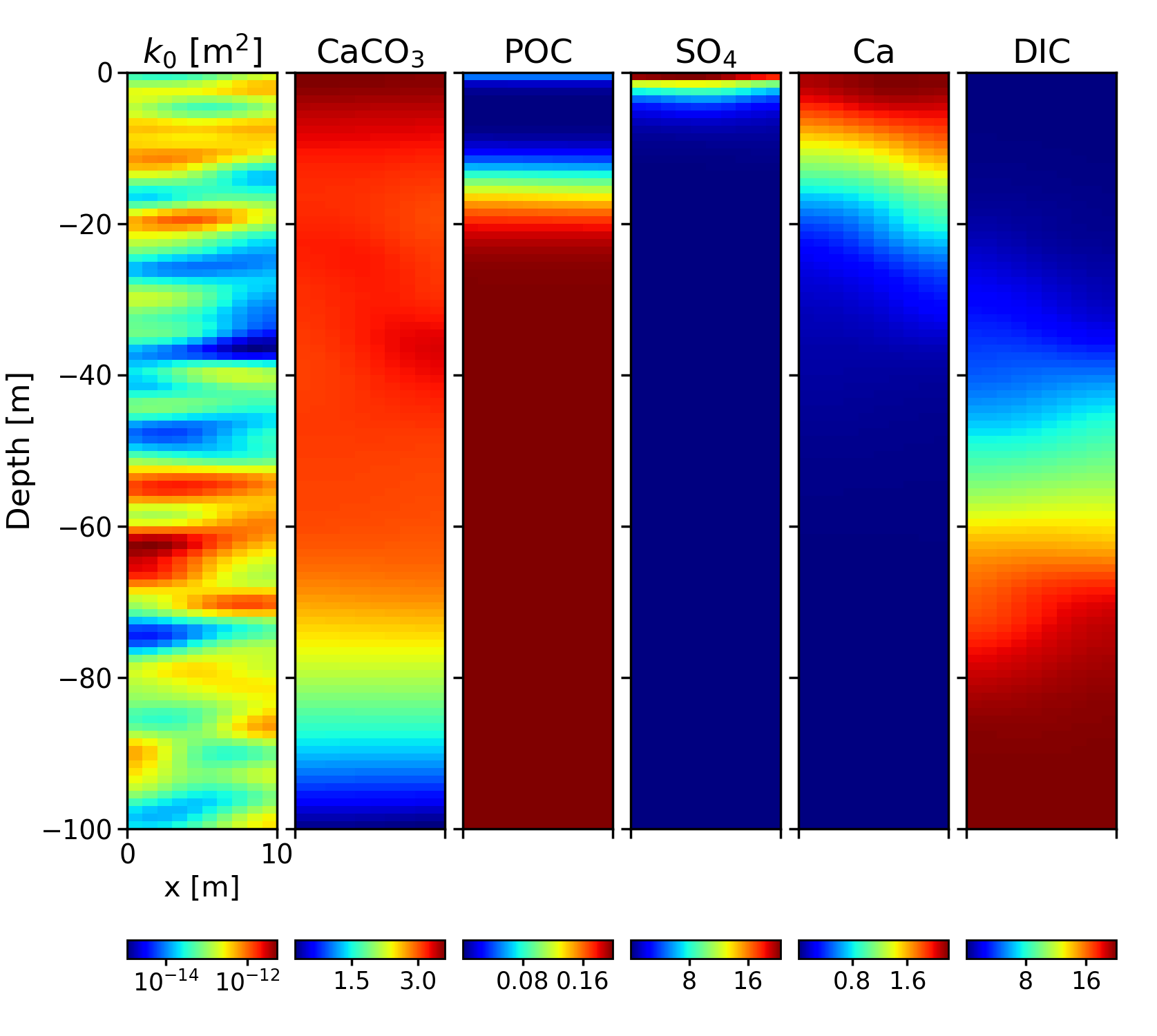}\\
{\small (a) 2D}
\end{minipage}\hfill
\begin{minipage}[b]{0.49\textwidth}\centering
\includegraphics[width=\linewidth]{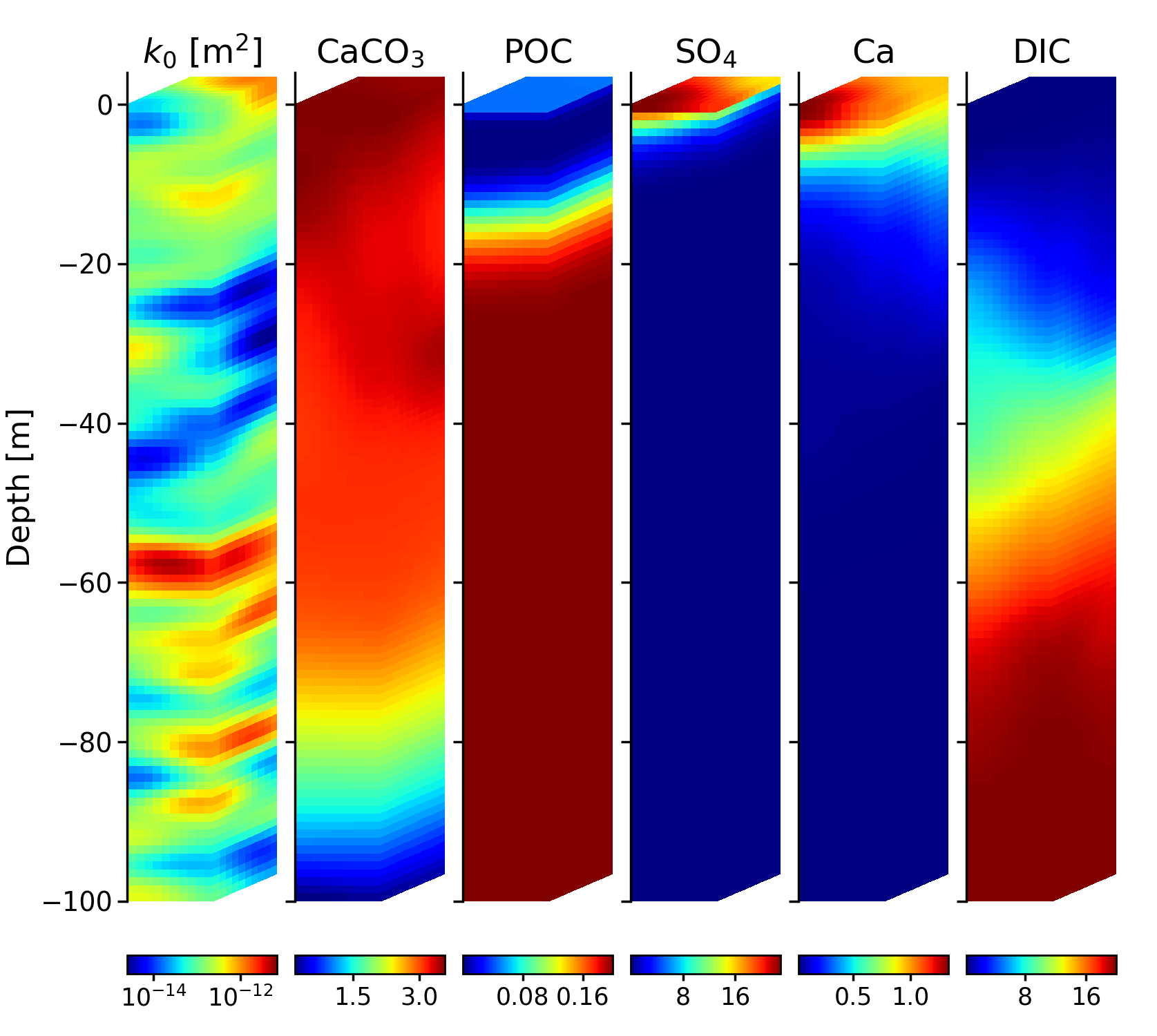}\\
{\small (b) 3D}
\end{minipage}
\caption{Reference solution for the test realization in 2D (a) and 3D (b). In each panel, the initial permeability is followed by the five concentration fields at the final time.}
\label{fig:reference-case3}
\end{figure}

\begin{figure}[h!]
\centering
\begin{minipage}[b]{0.49\textwidth}\centering
\includegraphics[width=\linewidth]{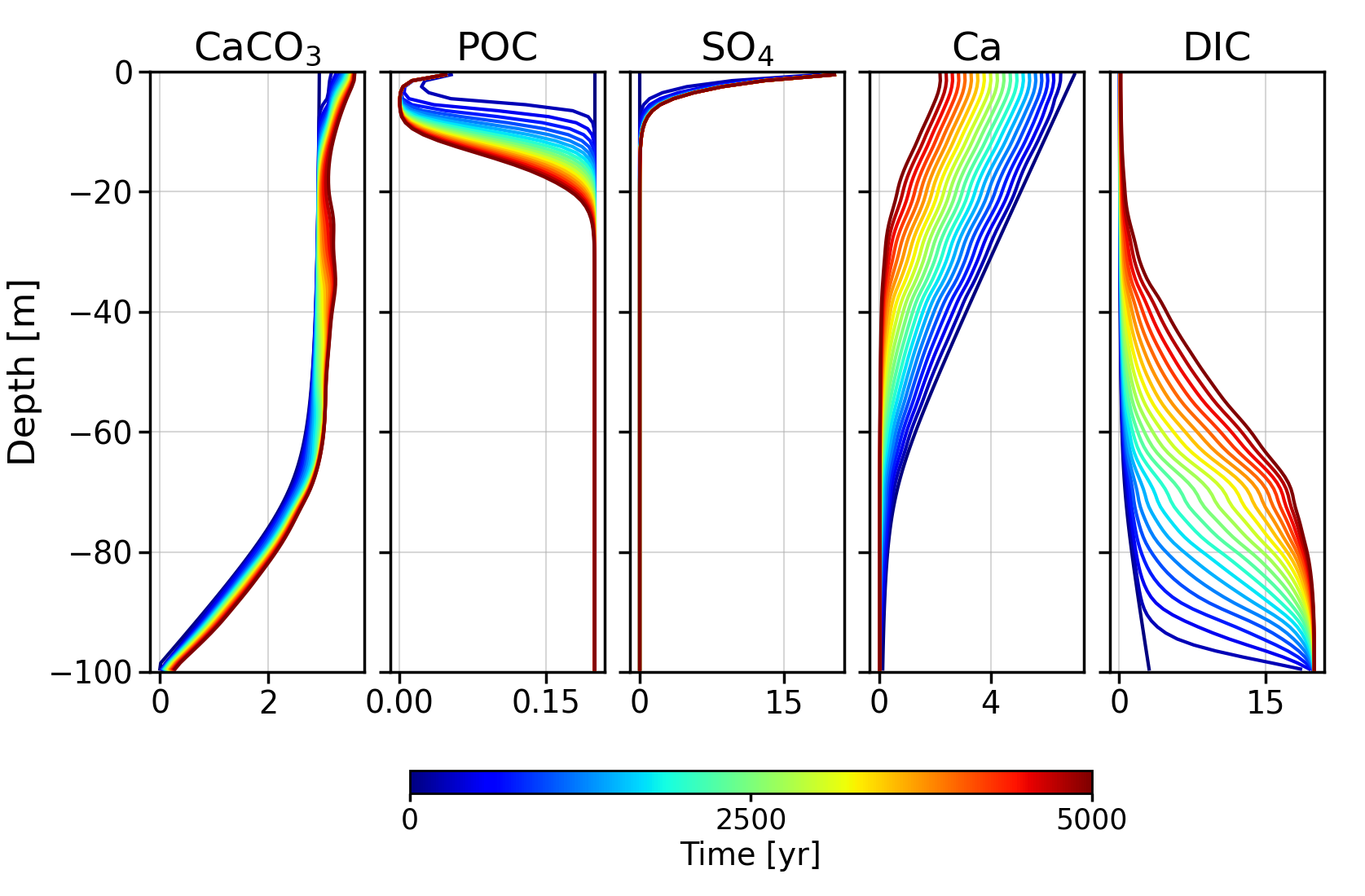}\\
{\small (a) 2D, $x=5$ m}
\end{minipage}\hfill
\begin{minipage}[b]{0.49\textwidth}\centering
\includegraphics[width=\linewidth]{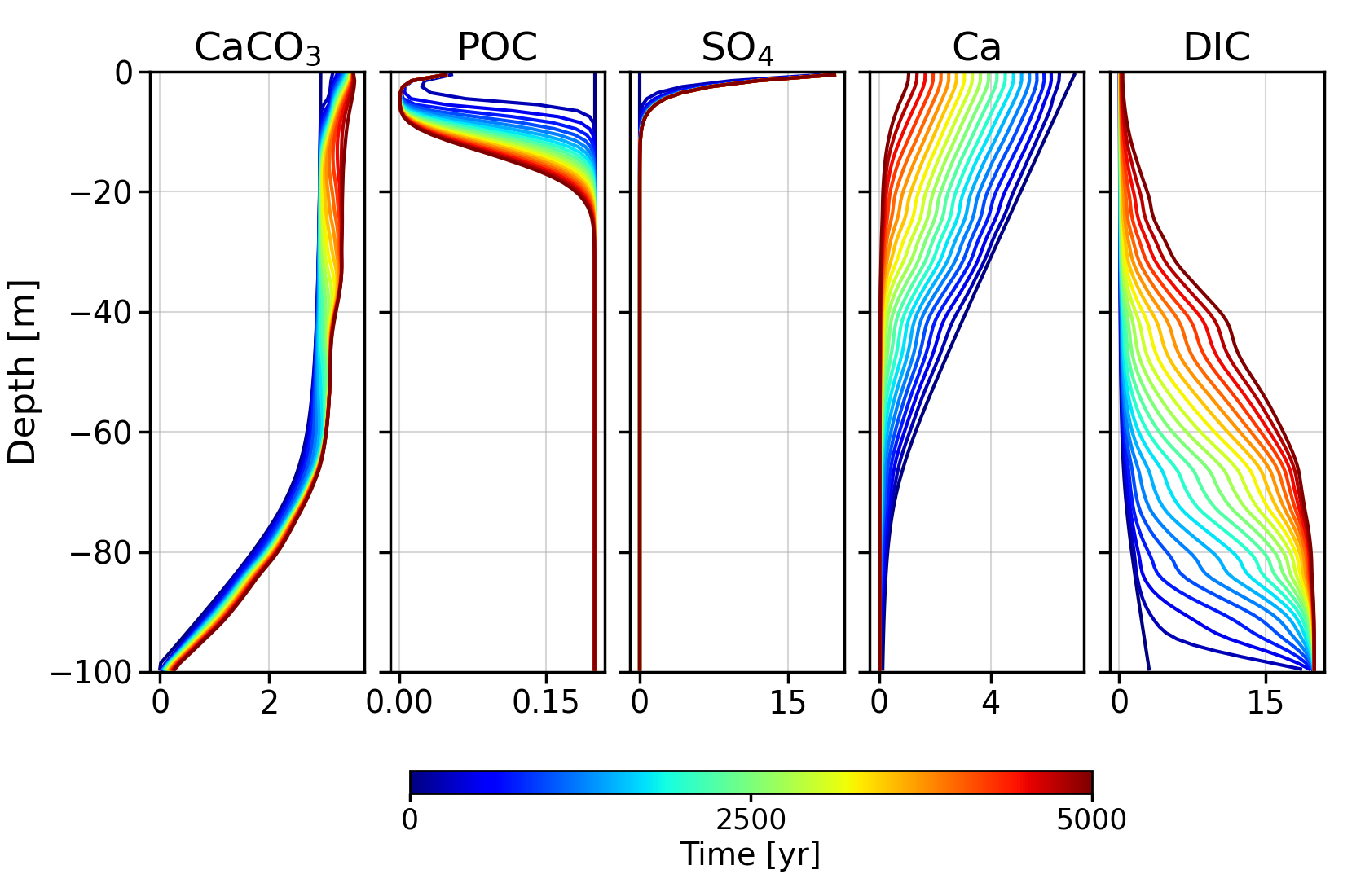}\\
{\small (b) 3D, $(x,y)=(5,5)$ m}
\end{minipage}
\caption{Reference concentration profiles along the central vertical line for the test realization in 2D (a) and 3D (b), every 250 years from $t=0$ to $T=5{,}000$ years.}
\label{fig:reference-case3-slices}
\end{figure}

Figure~\ref{fig:reference-case3} shows the initial permeability and the reference solutions at the final time in 2D and 3D, and Figure~\ref{fig:reference-case3-slices} shows their temporal evolution along the central vertical line.
The reference solutions serve as the spatial and temporal baseline for the surrogate comparisons that follow. 
The implementation uses Python with NumPy and SciPy, and the networks are trained with PyTorch.
The flow and implicit transport systems are solved using either sparse direct factorization or iterative methods.
The iterative configuration uses Jacobi-preconditioned conjugate gradients for flow and BiCGSTAB for transport, initialized with the previous solution and a relative residual tolerance of $10^{-10}$.
Chemistry is advanced by a cellwise nonlinear solve. 
All timing runs were performed on an Apple M2 Max machine with 12 CPU cores and 32 GB of memory. 
Reference timings are reported alongside the learned configurations in Table~\ref{tab:flow-timing}.
Iterative solves reduce the linear-solver cost on the 3D grid, while chemistry remains the dominant cost.

\begin{table}[h!]
\centering
\caption{Final relative $L^2$ concentration errors [\%] for the test realization and evaluation set (mean $\pm$ standard deviation), relative to the full-order reference. Global denotes the global inverse atoms, local the local atoms, and $\mathrm{ms}_{\mathrm{g}}$ and $\mathrm{ms}_{\mathrm{pg}}$ the pure multiscale transport without learning, with the exact multiscale Galerkin or LSPG system.}
\label{tab:flow-errors}
\setlength{\tabcolsep}{4.8pt}
\begin{tabular}{llcccccc}
\hline
Dim. & Run & CaCO$_3$ ($c_1$) & POC ($c_2$) & SO$_4$ ($c_3$) & Ca ($c_4$) & DIC ($c_5$) & All species \\
\hline
\multicolumn{8}{l}{Test realization} \\
\hline
2D & $\mathrm{C}$ & $0.02$ & $0.33$ & $0.01$ & $0.04$ & $0.001$ & $0.01$ \\
 & $\mathrm{CF}$ & $0.02$ & $0.33$ & $0.08$ & $0.09$ & $0.04$ & $0.04$ \\
 & $\mathrm{CT}$, global & $0.56$ & $0.57$ & $1.99$ & $2.06$ & $0.79$ & $0.86$ \\
 & $\mathrm{CT}$, local & $0.08$ & $0.44$ & $0.64$ & $0.29$ & $0.12$ & $0.17$ \\
 & $\mathrm{CT}$, $\mathrm{ms}_{\mathrm{g}}$ & $0.05$ & $0.34$ & $0.02$ & $0.08$ & $0.03$ & $0.03$ \\
 & $\mathrm{CT}$, $\mathrm{ms}_{\mathrm{pg}}$ & $0.08$ & $0.34$ & $0.02$ & $0.09$ & $0.03$ & $0.04$ \\
 & $\mathrm{CFT}$, global & $0.58$ & $0.56$ & $2.15$ & $1.69$ & $0.94$ & $1.00$ \\
 & $\mathrm{CFT}$, local & $0.08$ & $0.44$ & $0.66$ & $0.36$ & $0.13$ & $0.18$ \\
 & $\mathrm{CFT}$, $\mathrm{ms}_{\mathrm{g}}$ & $0.05$ & $0.34$ & $0.05$ & $0.08$ & $0.03$ & $0.03$ \\
 & $\mathrm{CFT}$, $\mathrm{ms}_{\mathrm{pg}}$ & $0.08$ & $0.34$ & $0.05$ & $0.09$ & $0.03$ & $0.04$ \\
\hline
3D & $\mathrm{C}$ & $0.02$ & $0.43$ & $0.005$ & $0.05$ & $0.002$ & $0.01$ \\
 & $\mathrm{CF}$ & $0.09$ & $0.44$ & $0.21$ & $0.62$ & $0.13$ & $0.14$ \\
 & $\mathrm{CT}$, global & $0.50$ & $0.42$ & $1.55$ & $1.80$ & $1.52$ & $1.49$ \\
 & $\mathrm{CT}$, local & $0.40$ & $0.55$ & $0.79$ & $2.69$ & $0.35$ & $0.37$ \\
 & $\mathrm{CT}$, $\mathrm{ms}_{\mathrm{g}}$ & $0.32$ & $0.59$ & $0.44$ & $0.83$ & $0.12$ & $0.16$ \\
 & $\mathrm{CT}$, $\mathrm{ms}_{\mathrm{pg}}$ & $0.40$ & $0.72$ & $0.53$ & $0.99$ & $0.14$ & $0.18$ \\
 & $\mathrm{CFT}$, global & $0.49$ & $0.48$ & $1.69$ & $2.45$ & $1.49$ & $1.47$ \\
 & $\mathrm{CFT}$, local & $0.36$ & $0.56$ & $0.93$ & $3.29$ & $0.44$ & $0.47$ \\
 & $\mathrm{CFT}$, $\mathrm{ms}_{\mathrm{g}}$ & $0.32$ & $0.59$ & $0.44$ & $0.84$ & $0.13$ & $0.16$ \\
 & $\mathrm{CFT}$, $\mathrm{ms}_{\mathrm{pg}}$ & $0.40$ & $0.72$ & $0.53$ & $1.00$ & $0.14$ & $0.18$ \\
\hline
\multicolumn{8}{l}{Evaluation set} \\
\hline
2D & $\mathrm{C}$ & $0.02 \pm 0.01$ & $0.36 \pm 0.11$ & $0.004 \pm 0.002$ & $0.03 \pm 0.01$ & $0.0009 \pm 0.0003$ & $0.01 \pm 0.002$ \\
 & $\mathrm{CF}$ & $0.04 \pm 0.04$ & $0.36 \pm 0.11$ & $0.22 \pm 0.18$ & $0.34 \pm 0.41$ & $0.08 \pm 0.04$ & $0.10 \pm 0.05$ \\
 & $\mathrm{CT}$, global & $0.85 \pm 0.34$ & $0.60 \pm 0.23$ & $2.29 \pm 1.68$ & $2.98 \pm 2.08$ & $1.32 \pm 0.64$ & $1.39 \pm 0.60$ \\
 & $\mathrm{CT}$, local & $0.12 \pm 0.11$ & $0.35 \pm 0.13$ & $0.77 \pm 0.38$ & $0.60 \pm 0.58$ & $0.15 \pm 0.05$ & $0.22 \pm 0.08$ \\
 & $\mathrm{CT}$, $\mathrm{ms}_{\mathrm{g}}$ & $0.07 \pm 0.04$ & $0.36 \pm 0.12$ & $0.04 \pm 0.02$ & $0.11 \pm 0.11$ & $0.03 \pm 0.01$ & $0.03 \pm 0.01$ \\
 & $\mathrm{CT}$, $\mathrm{ms}_{\mathrm{pg}}$ & $0.10 \pm 0.07$ & $0.37 \pm 0.12$ & $0.04 \pm 0.02$ & $0.14 \pm 0.14$ & $0.04 \pm 0.02$ & $0.04 \pm 0.02$ \\
 & $\mathrm{CFT}$, global & $0.86 \pm 0.34$ & $0.59 \pm 0.24$ & $2.36 \pm 1.69$ & $3.00 \pm 2.22$ & $1.41 \pm 0.73$ & $1.48 \pm 0.67$ \\
 & $\mathrm{CFT}$, local & $0.12 \pm 0.10$ & $0.34 \pm 0.13$ & $0.88 \pm 0.58$ & $0.72 \pm 0.70$ & $0.18 \pm 0.07$ & $0.26 \pm 0.12$ \\
 & $\mathrm{CFT}$, $\mathrm{ms}_{\mathrm{g}}$ & $0.07 \pm 0.04$ & $0.36 \pm 0.12$ & $0.05 \pm 0.02$ & $0.11 \pm 0.11$ & $0.03 \pm 0.01$ & $0.03 \pm 0.01$ \\
 & $\mathrm{CFT}$, $\mathrm{ms}_{\mathrm{pg}}$ & $0.10 \pm 0.07$ & $0.37 \pm 0.12$ & $0.05 \pm 0.02$ & $0.14 \pm 0.14$ & $0.04 \pm 0.02$ & $0.04 \pm 0.02$ \\
\hline
3D & $\mathrm{C}$ & $0.02 \pm 0.004$ & $0.39 \pm 0.09$ & $0.005 \pm 0.001$ & $0.05 \pm 0.05$ & $0.001 \pm 0.0002$ & $0.01 \pm 0.001$ \\
 & $\mathrm{CF}$ & $0.08 \pm 0.07$ & $0.41 \pm 0.10$ & $0.31 \pm 0.32$ & $0.83 \pm 0.90$ & $0.19 \pm 0.17$ & $0.20 \pm 0.16$ \\
 & $\mathrm{CT}$, global & $0.59 \pm 0.18$ & $0.65 \pm 0.25$ & $2.19 \pm 1.18$ & $2.41 \pm 1.60$ & $1.79 \pm 0.55$ & $1.78 \pm 0.53$ \\
 & $\mathrm{CT}$, local & $0.42 \pm 0.13$ & $0.69 \pm 0.37$ & $1.64 \pm 0.79$ & $1.30 \pm 0.70$ & $0.39 \pm 0.14$ & $0.51 \pm 0.14$ \\
 & $\mathrm{CT}$, $\mathrm{ms}_{\mathrm{g}}$ & $0.35 \pm 0.09$ & $0.63 \pm 0.31$ & $0.42 \pm 0.19$ & $0.65 \pm 0.35$ & $0.14 \pm 0.03$ & $0.18 \pm 0.03$ \\
 & $\mathrm{CT}$, $\mathrm{ms}_{\mathrm{pg}}$ & $0.45 \pm 0.11$ & $0.78 \pm 0.39$ & $0.53 \pm 0.26$ & $0.80 \pm 0.49$ & $0.16 \pm 0.04$ & $0.21 \pm 0.05$ \\
 & $\mathrm{CFT}$, global & $0.60 \pm 0.18$ & $0.66 \pm 0.26$ & $2.32 \pm 1.45$ & $2.74 \pm 1.99$ & $1.84 \pm 0.57$ & $1.83 \pm 0.55$ \\
 & $\mathrm{CFT}$, local & $0.42 \pm 0.12$ & $0.71 \pm 0.36$ & $1.70 \pm 0.85$ & $1.54 \pm 1.09$ & $0.47 \pm 0.25$ & $0.58 \pm 0.22$ \\
 & $\mathrm{CFT}$, $\mathrm{ms}_{\mathrm{g}}$ & $0.35 \pm 0.09$ & $0.63 \pm 0.31$ & $0.42 \pm 0.19$ & $0.65 \pm 0.35$ & $0.14 \pm 0.03$ & $0.18 \pm 0.03$ \\
 & $\mathrm{CFT}$, $\mathrm{ms}_{\mathrm{pg}}$ & $0.45 \pm 0.11$ & $0.78 \pm 0.39$ & $0.53 \pm 0.26$ & $0.80 \pm 0.49$ & $0.16 \pm 0.04$ & $0.21 \pm 0.05$ \\
\hline
\end{tabular}
\end{table}

\begin{table}[h!]
\centering
\caption{Solver and stage times [s] and speed-ups against the direct and iterative full-order references over the evaluation set (mean $\pm$ standard deviation). 
Global denotes the global inverse atoms, local the local atoms, and $\mathrm{ms}_{\mathrm{g}}$ and $\mathrm{ms}_{\mathrm{pg}}$ the pure multiscale transport without learning, with the exact multiscale Galerkin or LSPG system. }
\label{tab:flow-timing}
\begin{tabular}{llcccc}
\hline
Dim. & Run & Total & Flow & Transport & Chemistry \\
\hline
\multicolumn{6}{l}{Times [s]} \\
\hline
2D & Direct & $14.96 \pm 0.14$ & $0.24 \pm 0.01$ & $0.55 \pm 0.01$ & $14.11 \pm 0.12$ \\
 & Iterative & $14.85 \pm 0.18$ & $0.41 \pm 0.03$ & $0.38 \pm 0.01$ & $14.03 \pm 0.16$ \\
 & $\mathrm{C}$ & $0.86 \pm 0.01$ & $0.22 \pm 0.003$ & $0.52 \pm 0.01$ & $0.07 \pm 0.001$ \\
 & $\mathrm{CF}$ & $0.68 \pm 0.01$ & $0.03 \pm 0.0005$ & $0.52 \pm 0.01$ & $0.07 \pm 0.002$ \\
 & $\mathrm{CFT}$, global & $0.85 \pm 0.01$ & $0.03 \pm 0.0004$ & $0.59 \pm 0.005$ & $0.08 \pm 0.0007$ \\
 & $\mathrm{CFT}$, local & $0.95 \pm 0.01$ & $0.03 \pm 0.0005$ & $0.75 \pm 0.01$ & $0.08 \pm 0.0008$ \\
 & $\mathrm{CFT}$, $\mathrm{ms}_{\mathrm{g}}$ & $1.42 \pm 0.05$ & $0.03 \pm 0.001$ & $1.23 \pm 0.04$ & $0.08 \pm 0.003$ \\
 & $\mathrm{CFT}$, $\mathrm{ms}_{\mathrm{pg}}$ & $1.35 \pm 0.03$ & $0.03 \pm 0.0007$ & $1.17 \pm 0.02$ & $0.08 \pm 0.002$ \\
\hline
3D & Direct & $184.40 \pm 3.55$ & $8.56 \pm 0.14$ & $33.41 \pm 0.51$ & $142.21 \pm 3.02$ \\
 & Iterative & $147.90 \pm 2.03$ & $3.64 \pm 0.13$ & $3.30 \pm 0.05$ & $140.85 \pm 1.94$ \\
 & $\mathrm{C}$ & $42.31 \pm 0.48$ & $8.40 \pm 0.11$ & $33.07 \pm 0.36$ & $0.70 \pm 0.02$ \\
 & $\mathrm{CF}$ & $34.01 \pm 0.12$ & $0.32 \pm 0.002$ & $32.78 \pm 0.11$ & $0.65 \pm 0.02$ \\
 & $\mathrm{CFT}$, global & $3.22 \pm 0.04$ & $0.30 \pm 0.002$ & $1.59 \pm 0.01$ & $0.58 \pm 0.04$ \\
 & $\mathrm{CFT}$, local & $3.72 \pm 0.05$ & $0.30 \pm 0.003$ & $2.39 \pm 0.02$ & $0.59 \pm 0.04$ \\
 & $\mathrm{CFT}$, $\mathrm{ms}_{\mathrm{g}}$ & $8.77 \pm 0.11$ & $0.31 \pm 0.003$ & $7.29 \pm 0.08$ & $0.63 \pm 0.03$ \\
 & $\mathrm{CFT}$, $\mathrm{ms}_{\mathrm{pg}}$ & $8.49 \pm 0.10$ & $0.31 \pm 0.002$ & $7.01 \pm 0.08$ & $0.63 \pm 0.03$ \\
\hline
\multicolumn{6}{l}{Speed-up against the direct reference} \\
\hline
2D & $\mathrm{C}$ & $\times17.48 \pm 0.20$ & $\times1.08 \pm 0.03$ & $\times1.06 \pm 0.02$ & $\times190.50 \pm 4.02$ \\
 & $\mathrm{CF}$ & $\times22.15 \pm 0.35$ & $\times7.62 \pm 0.18$ & $\times1.06 \pm 0.02$ & $\times192.28 \pm 4.63$ \\
 & $\mathrm{CFT}$, global & $\times17.67 \pm 0.22$ & $\times7.36 \pm 0.21$ & $\times0.93 \pm 0.01$ & $\times180.45 \pm 2.29$ \\
 & $\mathrm{CFT}$, local & $\times15.68 \pm 0.21$ & $\times7.00 \pm 0.17$ & $\times0.74 \pm 0.01$ & $\times181.11 \pm 2.03$ \\
 & $\mathrm{CFT}$, $\mathrm{ms}_{\mathrm{g}}$ & $\times10.58 \pm 0.37$ & $\times7.24 \pm 0.29$ & $\times0.45 \pm 0.02$ & $\times183.33 \pm 6.69$ \\
 & $\mathrm{CFT}$, $\mathrm{ms}_{\mathrm{pg}}$ & $\times11.09 \pm 0.25$ & $\times7.42 \pm 0.25$ & $\times0.47 \pm 0.01$ & $\times187.85 \pm 4.49$ \\
\hline
3D & $\mathrm{C}$ & $\times4.36 \pm 0.10$ & $\times1.02 \pm 0.02$ & $\times1.01 \pm 0.02$ & $\times204.07 \pm 6.37$ \\
 & $\mathrm{CF}$ & $\times5.42 \pm 0.10$ & $\times27.17 \pm 0.41$ & $\times1.02 \pm 0.02$ & $\times218.97 \pm 6.36$ \\
 & $\mathrm{CFT}$, global & $\times57.33 \pm 1.44$ & $\times28.29 \pm 0.53$ & $\times21.02 \pm 0.32$ & $\times244.45 \pm 18.16$ \\
 & $\mathrm{CFT}$, local & $\times49.59 \pm 1.22$ & $\times28.07 \pm 0.57$ & $\times13.97 \pm 0.28$ & $\times241.64 \pm 16.44$ \\
 & $\mathrm{CFT}$, $\mathrm{ms}_{\mathrm{g}}$ & $\times21.03 \pm 0.46$ & $\times27.68 \pm 0.49$ & $\times4.58 \pm 0.09$ & $\times227.74 \pm 10.02$ \\
 & $\mathrm{CFT}$, $\mathrm{ms}_{\mathrm{pg}}$ & $\times21.72 \pm 0.55$ & $\times27.41 \pm 0.57$ & $\times4.77 \pm 0.10$ & $\times224.43 \pm 11.70$ \\
\hline
\multicolumn{6}{l}{Speed-up against the iterative reference} \\
\hline
2D & $\mathrm{C}$ & $\times17.36 \pm 0.29$ & $\times1.85 \pm 0.13$ & $\times0.72 \pm 0.01$ & $\times189.45 \pm 4.83$ \\
 & $\mathrm{CF}$ & $\times22.00 \pm 0.37$ & $\times13.09 \pm 0.77$ & $\times0.73 \pm 0.02$ & $\times191.19 \pm 4.41$ \\
 & $\mathrm{CFT}$, global & $\times17.55 \pm 0.25$ & $\times12.64 \pm 0.84$ & $\times0.63 \pm 0.01$ & $\times179.44 \pm 2.48$ \\
 & $\mathrm{CFT}$, local & $\times15.57 \pm 0.32$ & $\times12.02 \pm 0.78$ & $\times0.50 \pm 0.01$ & $\times180.12 \pm 3.38$ \\
 & $\mathrm{CFT}$, $\mathrm{ms}_{\mathrm{g}}$ & $\times10.51 \pm 0.35$ & $\times12.43 \pm 0.82$ & $\times0.31 \pm 0.01$ & $\times182.29 \pm 6.50$ \\
 & $\mathrm{CFT}$, $\mathrm{ms}_{\mathrm{pg}}$ & $\times11.01 \pm 0.25$ & $\times12.75 \pm 0.86$ & $\times0.32 \pm 0.01$ & $\times186.80 \pm 4.62$ \\
\hline
3D & $\mathrm{C}$ & $\times3.50 \pm 0.06$ & $\times0.43 \pm 0.02$ & $\times0.10 \pm 0.002$ & $\times202.12 \pm 5.82$ \\
 & $\mathrm{CF}$ & $\times4.35 \pm 0.06$ & $\times11.56 \pm 0.40$ & $\times0.10 \pm 0.001$ & $\times216.88 \pm 5.88$ \\
 & $\mathrm{CFT}$, global & $\times45.98 \pm 0.89$ & $\times12.04 \pm 0.46$ & $\times2.08 \pm 0.03$ & $\times242.06 \pm 17.01$ \\
 & $\mathrm{CFT}$, local & $\times39.77 \pm 0.71$ & $\times11.94 \pm 0.46$ & $\times1.38 \pm 0.02$ & $\times239.28 \pm 15.22$ \\
 & $\mathrm{CFT}$, $\mathrm{ms}_{\mathrm{g}}$ & $\times16.87 \pm 0.34$ & $\times11.78 \pm 0.38$ & $\times0.45 \pm 0.01$ & $\times225.60 \pm 10.35$ \\
 & $\mathrm{CFT}$, $\mathrm{ms}_{\mathrm{pg}}$ & $\times17.42 \pm 0.33$ & $\times11.66 \pm 0.45$ & $\times0.47 \pm 0.01$ & $\times222.21 \pm 9.68$ \\
\hline
\end{tabular}
\end{table}

To separate the effects of the three learned stages, we label each run by the stages it learns, with $\mathrm{C}$ for chemistry alone, $\mathrm{CF}$ for chemistry and flow, $\mathrm{CT}$ for chemistry and transport, and $\mathrm{CFT}$ for all three, while stages absent from the label retain their full-order solvers.
The coarse grid $\mathcal T_H$ shared by transport and flow divides the depth into $N_H=10$ layers of 10 m, so that the coarse cells are $10\times10$ m in 2D and $10\times10\times10$ m in 3D.
The 11 coarse nodes along the depth give 11 overlapping transport neighborhoods $\omega_j$, and the 11 horizontal coarse faces carry the flow interface bases.
The training data, network architectures, and offline costs of the three learned stages are given in Appendix~\ref{app:training}.
All learned configurations use the same chemistry network, and flow uses MLP trace prediction in 2D and CNN trace prediction in 3D.
Learned transport uses either the global inverse atoms (rows ``global'' in Tables~\ref{tab:flow-errors} and~\ref{tab:flow-timing}), with $K_i=64/256$ in 2D and $K_i=128/512$ in 3D trained on 4,000 fields, or the local atoms with the Galerkin system \eqref{eq:local-galerkin} (rows ``local''), with $K=32$ atoms per species and coarse cell trained on 330 fields and 50 validation fields.
To quantify the effect of learning, rows $\mathrm{ms}_{\mathrm{g}}$ and $\mathrm{ms}_{\mathrm{pg}}$ use the same multiscale space without learning, with the exact multiscale Galerkin or LSPG transport step \eqref{eq:local-galerkin} or \eqref{eq:local-lspg}, assembled from the fine-grid operator at every step. 

Table~\ref{tab:flow-errors} reports the final relative concentration error $e_i^{\mathrm{final}}=100\%\cdot\|\widehat C_i^{N_t}-C_{i,\mathrm{ref}}^{N_t}\|_2/\|C_{i,\mathrm{ref}}^{N_t}\|_2$ of species $i$ for the test realization and the evaluation set, where $\widehat C_i^{N_t}$ and $C_{i,\mathrm{ref}}^{N_t}$ contain the predicted and full-order cell values at $T=N_t\Delta t$.
The all-species error $e_{\mathrm{all}}^{\mathrm{final}}$ is defined in the same way for the vectors that concatenate all five species over all cells.
For the evaluation set, both errors are computed for each realization and reported as the mean and standard deviation.
Chemistry alone ($\mathrm{C}$) gives mean species errors below $1\%$.
Adding learned transport with global inverse atoms ($\mathrm{CT}$) gives mean all-species errors of $1.4\%$ in 2D and $1.8\%$ in 3D.
With all three stages learned ($\mathrm{CFT}$), these errors are $1.5\%$ and $1.8\%$, respectively.
In 3D, the $\mathrm{CFT}$ configuration with global inverse atoms and the augmented CNN flow has mean species errors of $0.60$--$2.7\%$.
The fully learned solver approximates the reference concentrations over 5,000 years with a modest change in aggregate error relative to $\mathrm{CT}$.
With full-order transport ($\mathrm{CF}$), the augmented 3D CNN gives a mean calcium error of $0.83\%$ and a mean all-species error of $0.20\%$.
The local atoms reduce the transport error substantially.
Their mean all-species errors are $0.22\%$ in 2D and $0.51\%$ in 3D for $\mathrm{CT}$, and $0.26\%$ and $0.58\%$ for $\mathrm{CFT}$, compared with $1.4$--$1.8\%$ for the global inverse atoms.
The largest reductions occur in sulfate, calcium, and DIC, which dominate the global-atom error.
For example, the mean DIC error decreases from $1.3\%$ to $0.15\%$ in 2D and from $1.8\%$ to $0.39\%$ in 3D.
For the 3D test realization, the local atoms give a larger Ca error than the global inverse atoms, $3.29\%$ against $2.45\%$ for $\mathrm{CFT}$, but a three times smaller all-species error, $0.47\%$ against $1.47\%$, and over the evaluation set their mean Ca error is also smaller, $1.54\%$ against $2.74\%$.
Without learning, the exact multiscale Galerkin and LSPG transport steps give mean all-species errors of $0.03\%$ and $0.04\%$ in 2D and $0.18\%$ and $0.21\%$ in 3D for $\mathrm{CFT}$.
The exact multiscale errors measure the spatial reduction of the multiscale space, so learning the operator adds about $0.2$--$0.4$ percentage points with local atoms.

Table~\ref{tab:flow-timing} reports online times and speed-ups against the direct and iterative references.
Ratios are computed per realization before averaging.
Setup includes model loading, trace prediction, local extensions, and initial reduced-system construction and solution. It is included in total time and excluded from stage speed-ups. The reported speed-ups exclude offline data generation and training.
For flow, the speed-up comes from the cheap online solve in the mixed multiscale space, whose basis is built once per realization from interface traces and local extensions, without local eigenvalue problems. We note that learning the traces pays off only when these local problems are expensive, for example for larger coarse blocks or finer grids.
In 2D, the complete solver takes 0.85 s with a $1.5\%$ all-species error and an approximately $18\times$ total speed-up.
The chemistry--flow configuration ($\mathrm{CF}$) is fastest in 2D, with a mean total time of 0.68 s, compared with 0.85 s for $\mathrm{CFT}$.
Its mean all-species error is also smaller ($0.10\%$ versus $1.5\%$).
On this grid, learned transport is slower than both full-order transport solvers, so replacing every stage does not give the best accuracy--cost balance.
In 3D, replacing chemistry reduces mean total time from 184.40 to 42.31 s.
Learning flow and transport as well reduces it to 3.22 s, including setup.
The $\mathrm{CFT}$ configuration is fastest in 3D and is $57\times$ faster than the direct reference and $46\times$ faster than the iterative reference.
Relative to the iterative reference, flow, transport, and chemistry speed-ups are $12\times$, $2.1\times$, and $242\times$.
The local atoms assemble and factorize a reduced matrix at every step.
In 2D, where full-order transport is cheap, this makes them slower than both the finite-volume solve and the global inverse atoms, with transport times of $0.75$\,s, $0.52$\,s, and $0.59$\,s.
In 3D, they reduce the transport time from $32.8$\,s to $2.39$\,s, compared with $1.59$\,s for the global inverse atoms, and the fully learned solver takes $3.72$\,s, compared with $3.22$\,s, which is $50\times$ and $40\times$ faster than the direct and iterative references.
The pure multiscale transport assembles and projects the fine-grid operator at every step and takes $7.3$\,s (Galerkin) and $7.0$\,s (LSPG) in 3D, and the fully multiscale solver takes $8.8$\,s and $8.5$\,s.
The multiscale transport is thus faster than the direct solve but slower than the iterative full-order transport ($3.3$\,s), since it still assembles and projects the fine-grid operator.
Learning the transport operator therefore reduces the total time by a further factor of about $2.3$--$2.7$ relative to the multiscale method, at the price of a larger error.
In 2D, the multiscale transport ($1.2$\,s) is slower than full-order transport.
Overall, in 3D the fully learned solver with local atoms is $50\times$ faster than the direct reference and $40\times$ faster than the iterative reference at a mean all-species error of $0.58\%$.
Its online transport step neither assembles nor solves a fine-grid system, so the gain is expected to grow with the cost of the full-order transport, for example on finer grids or with more species.

\section{Conclusions}
\label{sec:conclusions}

We presented the Split-and-Learn framework, which keeps the sequential operator-splitting scheme of the full-order model and replaces chemistry, flow, and transport by data-driven surrogates with their own structure.
Chemistry uses bounded reaction extents with stoichiometric reconstruction, flow uses learned interface traces with a reduced mixed solve, and transport uses global inverse or local operator atoms on a data-driven multiscale space.
{We proved that the split scheme with evolving porosity and permeability is stable for a slowly changing porosity and Lipschitz reactions, with a growth rate that does not depend on the time step, and that the Galerkin multiscale and learned transport steps keep this stability, where the learning adds the error of the reduced operator to the growth rate.
The deviation of the multiscale and learned solvers from the full-order split solver is bounded by the accumulated one-step errors, which are the projection error of the multiscale space, given by the POD truncation and the distance from the snapshots, the error of the learned reduced operator, which enters in proportion to the time step, and the errors of the learned flow and chemistry.}

In 3D, the fully learned solver is $57$ times faster than the direct reference with a mean all-species error of $1.8\%$ for the global inverse atoms, and $50$ times faster with $0.58\%$ for the local atoms.
Learning chemistry gives the largest single saving, since chemistry dominates the full-order cost, and learning flow and transport as well reduces the total time from $42$\,s to $3.2$\,s.
Learning local pieces of the forward transport operator instead of the global inverse map reduces the mean all-species error by a factor of three to six, with $32$ atoms per species and coarse cell and $330$ instead of $4{,}000$ training fields, at the price of a reduced solve at every step.
The comparison with the same multiscale spaces without learning shows where the gains come from.
In 3D, learned transport is $3$ to $5$ times faster than exact multiscale transport because it avoids the online assembly and projection of the fine-grid operator, at the price of a larger error.
In 2D, where full-order transport is cheap, learning chemistry and flow alone gives the best accuracy--cost balance, so the best combination depends on the relative costs of the stages.

The results are limited to the five-species reaction network and to media drawn from the same random-field model as the training data.
Offline construction requires reference data and access to the numerical operators, and fixed spatial bases can limit the accuracy when the medium changes strongly.
Future work will develop a posteriori error estimates along learned trajectories to guide adaptive basis enrichment or selective reference solves.
We will also extend the method to broader reaction networks{, with carbonate speciation, dissolved organic carbon, and further microbial pathways such as methanogenesis and the anaerobic oxidation of methane,} and to larger grids, where the full-order transport and the local flow problems are more expensive and the learned stages are expected to pay off more.

\section*{Acknowledgments}\label{sec:disc}

S. W. Cheung was partially supported by Laboratory Directed Research and Development (LDRD)
Program by the U.S. Department of Energy under Project Number 25-ERD-051 and 
Scientific Discovery through Advanced Computing (SciDAC) program by the U.S. Department of Energy, Office of Science, Office of Advanced Scientific Computing Research,  through the LEADS SciDAC Institute under Project Number SCW1933. 
Lawrence Livermore National Laboratory is operated by Lawrence Livermore National Security, LLC, for the U.S. Department of Energy, National Nuclear Security Administration under Contract DE-AC52-07NA27344. 
IM release number: LLNL-JRNL-2026814.

\bibliographystyle{plain}
\bibliography{lit}

\appendix

\section{Stage estimates for the stability analysis}
\label{app:stage-estimates}

We prove the one-step estimates of the stages used in Theorem~\ref{thm:split-stability}, with the norms and the assumptions (A1) and (A2) of Section~\ref{sec:stability}.

\begin{lemma}[Change of the porosity]
\label{lem:phase-change}
Let (A1) hold. Then $\|C\|_{{\phi^{n+1}}}\leq(1+\mu_{\phi}\Delta t)\|C\|_{{\phi^n}}$ for every concentration field $C$.
\end{lemma}

\begin{proof}
Since $\eta_i(\phi)$ is $\phi$ or $1-\phi$, $|\eta_i(\phi_K^{n+1})-\eta_i(\phi_K^n)|=|\phi_K^{n+1}-\phi_K^n|$, and (A1) gives $\eta_i(\phi_K^{n+1})\leq\eta_i(\phi_K^n)+|\phi_K^{n+1}-\phi_K^n|\leq(1+2\mu_{\phi}\Delta t)\,\eta_i(\phi_K^n)$ in each cell.
Multiplying by $|K|\,C_{i,K}^2$ and summing over the species and the cells gives $\|C\|_{{\phi^{n+1}}}^2\leq(1+2\mu_{\phi}\Delta t)\|C\|_{{\phi^n}}^2$, and $\sqrt{1+2a}\leq1+a$ with $a=\mu_{\phi}\Delta t$ gives the estimate.
\end{proof}

\begin{lemma}[Transport stage]
\label{lem:transport-stage}
Let (A1) hold. Then, for $i=1,\ldots,5$ and all $w$,
\begin{equation}
 (A_i(k^n,\phi^n,q^n)w,w)\geq-\mu_{\phi}\,(M_i(\phi^n)w,w),
 \label{eq:transport-energy}
\end{equation}
and the transport step \eqref{eq:operator-form} with a general right-hand side,
\begin{equation}
 \left[\frac{M_i(\phi^n)}{\Delta t}+A_i(k^n,\phi^n,q^n)\right]w=\frac{M_i(\phi^n)}{\Delta t}v+g,
 \label{eq:stage-system}
\end{equation}
has a unique solution $w$ for all $v$ and $g$, which satisfies
\begin{equation*}
 \|w\|_{{\phi^n}}\leq(1+2\mu_{\phi}\Delta t)\bigl(\|v\|_{{\phi^n}}+\Delta t\,\|M_i(\phi^n)^{-1/2}g\|_2\bigr).
\end{equation*}
\end{lemma}

\begin{proof}
We split $A_i=\mathcal D_i+\mathcal C_i$ into the diffusion matrix, with $(\mathcal D_iw,w)\geq0$, and the upwind {advection} matrix, for which $(\mathcal C_iw,w)\geq\frac12(\mathcal G_iw,w)$ with the discrete divergence $\mathcal G_i=\operatorname{diag}(\sum_{e\subset\partial K}Q_{i,K,e}^n)$ of the phase flux \cite{samarskii2001additive,vabishchevich2012explicit}.
For the solid species, the burial flux is the same on all horizontal faces, so $\mathcal G_i=0$.
For the dissolved species, the discrete balance \eqref{eq:fv-darcy-flux} gives $(\mathcal G_i)_{KK}=-|K|(\phi_K^n-\phi_K^{n-1})/\Delta t$, which can increase the energy only where dissolution opens the pores.
In these cells, (A1) gives $\phi_K^n-\phi_K^{n-1}\leq2\mu_{\phi}\Delta t\,\phi_K^n$, so that, with $x^+=\max(x,0)$,
\begin{equation*}
 (A_iw,w)\geq-\frac12\sum_K|K|\,\frac{(\phi_K^n-\phi_K^{n-1})^+}{\Delta t}\,w_K^2\geq-\mu_{\phi}\sum_K|K|\,\phi_K^n\,w_K^2=-\mu_{\phi}\,(M_i(\phi^n)w,w),
\end{equation*}
which is \eqref{eq:transport-energy}.
Taking the inner product of \eqref{eq:stage-system} with $w$ and using \eqref{eq:transport-energy} and the Cauchy--Schwarz inequality,
\begin{equation*}
 \frac{1-\mu_{\phi}\Delta t}{\Delta t}\|w\|_{{\phi^n}}^2
 \leq\frac{1}{\Delta t}\|w\|_{{\phi^n}}^2+(A_iw,w)
 =\frac{1}{\Delta t}(M_i(\phi^n)v,w)+(g,w)
 \leq\Bigl(\frac{1}{\Delta t}\|v\|_{{\phi^n}}+\|M_i(\phi^n)^{-1/2}g\|_2\Bigr)\|w\|_{{\phi^n}},
\end{equation*}
and $(1-\mu_{\phi}\Delta t)^{-1}\leq1+2\mu_{\phi}\Delta t$ gives the estimate, since $(1-x)(1+2x)=1+x(1-2x)\geq1$ for $0\leq x=\mu_{\phi}\Delta t\leq1/2$.
For $v=0$ and $g=0$ it gives $w=0$, so the matrix of \eqref{eq:stage-system} is nonsingular and the solution is unique.
\end{proof}

By the upwind and two-point fluxes, the matrix of \eqref{eq:stage-system} also has nonpositive off-diagonal entries and positive row sums, at least $|K|\phi_K^{n-1}/\Delta t$ for the dissolved and $|K|(1-\phi_K^n)/\Delta t$ for the solid species.
Hence it is an M-matrix with a nonnegative inverse, and the reference transport step keeps nonnegative concentrations nonnegative.

\begin{lemma}[Reaction stage]
\label{lem:reaction-stage}
Let (A2) hold, and let $m(C)=\Delta t\,R_1(\mathcal R_{\Delta t}(C))$ be the CaCO$_3$ precipitated or dissolved in one step, so that $\Delta\phi^{n+1}=-\gamma_\phi\,m(C^{*,n})$.
Then, for inputs $C$ and $\widetilde C$ whose outputs $\mathcal R_{\Delta t}(C)$ and $\mathcal R_{\Delta t}(\widetilde C)$ have cell values in $\mathcal C$,
\begin{equation*}
 \|\mathcal R_{\Delta t}(C)-\mathcal R_{\Delta t}(\widetilde C)\|_{{\phi^n}}\leq(1+2L_r\Delta t)\|C-\widetilde C\|_{{\phi^n}},\quad
 \|m(C)-m(\widetilde C)\|\leq2\eta_{\min}^{-1/2}L_r\Delta t\,\|C-\widetilde C\|_{{\phi^n}}.
\end{equation*}
\end{lemma}

\begin{proof}
In a cell, the outputs $x=C+\Delta t\,R(x)$ and $\widetilde x=\widetilde C+\Delta t\,R(\widetilde x)$ of the backward Euler reaction step \eqref{eq:reaction-disc} satisfy $|x-\widetilde x|_{{K,\phi^n}}\leq|C-\widetilde C|_{{K,\phi^n}}+L_r\Delta t\,|x-\widetilde x|_{{K,\phi^n}}$, so that
$|x-\widetilde x|_{{K,\phi^n}}\leq(1-L_r\Delta t)^{-1}|C-\widetilde C|_{{K,\phi^n}}\leq(1+2L_r\Delta t)\,|C-\widetilde C|_{{K,\phi^n}}$.
Since $|R_1(x)-R_1(\widetilde x)|\leq\eta_{\min}^{-1/2}|R(x)-R(\widetilde x)|_{{K,\phi^n}}$ and $(1-L_r\Delta t)^{-1}\leq2$ by (A2),
\begin{equation*}
 |m(C)-m(\widetilde C)|=\Delta t\,|R_1(x)-R_1(\widetilde x)|\leq\eta_{\min}^{-1/2}L_r\Delta t\,|x-\widetilde x|_{{K,\phi^n}}\leq2\eta_{\min}^{-1/2}L_r\Delta t\,|C-\widetilde C|_{{K,\phi^n}}.
\end{equation*}
Multiplying the squares by $|K|$ and summing over the cells gives both estimates.
\end{proof}

Lemmas~\ref{lem:phase-change}--\ref{lem:reaction-stage} bound each stage in a fixed medium.
The stages are coupled through the medium, which enters the flow through the permeability $k_U=k(\phi_U)$ and the source $s_U=-|K|\Delta\phi_U/\Delta t$, and the transport through the storage, the diffusion coefficients, and the Darcy fluxes.
{We write the fluxes of a state $U$ as $\mathcal F(U)=T(k_U)\nabla_hp_U$, where $T(k)$ is the transmissibilities $T_{KL}$ and $T_{K,e}$ of \eqref{eq:fv-darcy-flux}--\eqref{eq:fv-flow-matrix}, $\nabla_hp$ are the face differences $p_K-p_L$ and $p_K-p_e$ with the prescribed pressures $p_e$, and $p_U$ solves the pressure system \eqref{eq:fv-flow-matrix}, the cell balance of these fluxes with the source $s_U$.}

\begin{lemma}[Dependence on the medium]
\label{lem:medium}
Let (A1) hold, let $U$ be a state with the boundary data of $U^n$, and let $\widetilde C=\mathcal T_{\Delta t}(C_U,U,\mathcal F(U))$ and $\overline C=\mathcal T_{\Delta t}(C_U,U^n,q^n)$ be the transport of the same concentrations in the two media.
Assume that the transmissibilities of \eqref{eq:fv-darcy-flux} lie in $[{T_{\min}},{T_{\max}}]$ with ${T_{\min}}>0$, that $\phi\mapsto T(k(\phi))$ and the diffusive face coefficients $\phi\mapsto T^D(\phi)$ of \eqref{eq:fv-diffusion} are Lipschitz with constants $L_T$ and $L_D$, and that $|(\nabla_hp^n)_e|\leq G_{{p}}$, $|\widetilde C|\leq C_{\max}$, $|\nabla_h\widetilde C|\leq G_C$, and $|\widetilde C_K-C_{U,K}|\leq c_t\Delta t$ in the Euclidean norm over the species.
Then
\begin{equation}
\begin{gathered}
 \|\mathcal F(U)-q^n\|\leq L_F\bigl(\|\phi_U-\phi^n\|+\Delta t^{-1}\|\Delta\phi_U-\Delta\phi^n\|\bigr),\\
 L_F=\max\{{T_{\max}}{T_{\min}}^{-1}C_P{|K|},\,(1+{T_{\max}}{T_{\min}}^{-1}){G_p}L_T\},
\end{gathered}
 \label{eq:flux-bound}
\end{equation}
and the bound \eqref{eq:medium-bound} holds for 
$\|\widetilde C-\overline C\|_{{\phi^n}}$, 
$\Lambda_\phi=2\bigl(\eta_{\min}^{-1/2}c_t+c_A(1+L_F)\bigr)$, $\Lambda_{\Delta\phi}=2c_AL_F$, and $c_A=\eta_{\min}^{-1/2}c_h\max(L_DG_C,C_{\max})$, where $C_P$ is the constant of the discrete Poincar\'e inequality $(s,e)\leq C_P\|s\|\|\nabla_he\|$, and $c_h$ is a mesh constant. 
{In particular, for the same medium and two fluxes $q$ and $q'$, the transport steps of the same concentrations satisfy $\|\mathcal T_{\Delta t}(C_U,U,q)-\mathcal T_{\Delta t}(C_U,U,q')\|_{{\phi^n}}\leq\Delta t\,L_q\|q-q'\|$ with $L_q=2\eta_{\min}^{-1/2}c_hC_{\max}$.}
\end{lemma}

\begin{proof}
\textit{Flow.}
The pressure difference $e_p=p_U-p^n$ vanishes on the pressure boundaries, and subtracting the two flow problems and taking the inner product with $e_p$ gives
\begin{equation*}
 (T(k_U)\nabla_he_p,\nabla_he_p)=(s_U-s^n,e_p)-\bigl((T(k_U)-T(k^n))\nabla_hp^n,\nabla_he_p\bigr).
\end{equation*}
With $\|(T(k_U)-T(k^n))\nabla_hp^n\|\leq G_{{p}}L_T\|\phi_U-\phi^n\|$ and $\|s_U-s^n\|\leq|K|_{\max}\Delta t^{-1}\|\Delta\phi_U-\Delta\phi^n\|$,
\begin{equation*}
{T}_{\min}\|\nabla_he_p\|\leq C_P|K|\Delta t^{-1}\|\Delta\phi_U-\Delta\phi^n\|+G_{{p}} L_T\|\phi_U-\phi^n\|,
\end{equation*}
and $\mathcal F(U)-q^n=T(k_U)\nabla_he_p+(T(k_U)-T(k^n))\nabla_hp^n$ gives \eqref{eq:flux-bound}.

\noindent\textit{Transport.}
Let $A_i(U)$ and $b_i(U)$ be the transport matrix and boundary vector of the medium of $U$ with the fluxes $\mathcal F(U)$.
Subtracting the transport steps \eqref{eq:operator-form} of $\widetilde C_i$ and $\overline C_i$, the difference $\widetilde w_i=\widetilde C_i-\overline C_i$ solves \eqref{eq:stage-system} with $v=0$ and $g=g_i^S+g_i^A$, where
\begin{equation*}
 g_i^S=\bigl(M_i(\phi^n)-M_i(\phi_U)\bigr)\frac{\widetilde C_i-C_{U,i}}{\Delta t},\quad
 g_i^A=-\bigl(A_i(U)-A_i(U^n)\bigr)\widetilde C_i+b_i(U)-b_i(U^n),
\end{equation*}
and Lemma~\ref{lem:transport-stage}, summed over the species, gives $\|\widetilde w\|_{{\phi^n}}\leq2\Delta t\,\|M(\phi^n)^{-1/2}(g^S+g^A)\|_2$.
Since $|\eta_i(\phi_U)-\eta_i(\phi^n)|=|\phi_U-\phi^n|$ and $\eta_i\geq\eta_{\min}$, the storage part satisfies $\|M(\phi^n)^{-1/2}g^S\|_2\leq\eta_{\min}^{-1/2}c_t\|\phi_U-\phi^n\|$.
The part $g^A$ is the change of the diffusive and advective face fluxes of $\widetilde C$ between the two media, so
\begin{equation*}
 \|M(\phi^n)^{-1/2}g^A\|_2\leq\eta_{\min}^{-1/2}c_h\bigl(L_DG_C\|\phi_U-\phi^n\|+C_{\max}\|\mathcal F(U)-q^n\|\bigr)\leq c_A\bigl(\|\phi_U-\phi^n\|+\|\mathcal F(U)-q^n\|\bigr),
\end{equation*}
and inserting \eqref{eq:flux-bound} gives \eqref{eq:medium-bound}. 
{For the same medium and two fluxes, $g^S$ and the diffusive part of $g^A$ vanish, which gives the bound with $L_q$.}
\end{proof}

\section{One-step errors of the multiscale and learned stages}
\label{app:local-errors}

We bound the one-step errors $\delta_{\mathrm{ms}}^n$, $\delta_\ell^n$, and $\varepsilon_f^n$ of Theorem~\ref{thm:learned-errors}, following the analysis of multiscale methods \cite{efendiev2013gmsfem,vasilyeva2024gmsfemgraph} for the transport.
For one species $i$ and the medium $\kappa=(k,\phi,q)$ of a multiscale or learned state, let $M=M_i(\phi)$, $A=A_i(\kappa)=\mathcal D+\mathcal C$, and $\mathcal L=\mathcal L_i(\kappa)=M/\Delta t+A$ be the matrices of the transport step \eqref{eq:training-system}, with the diffusion matrix $\mathcal D=\mathcal D^T\geq0$ and the upwind advection matrix $\mathcal C$, and let $\|w\|_M^2=(Mw,w)$ and $\|w\|_{\mathcal D}^2=(\mathcal Dw,w)$ be the mass norm, which is the phase norm of the state, and the energy norm, the sum of $T^D_{i,KL}(w_K-w_L)^2$ of \eqref{eq:fv-diffusion} over the faces.
We assume that $\eta_i(\phi_K)\geq\eta_{\min}$ and that the energy bound \eqref{eq:transport-energy} holds in this medium,
\begin{equation}
 (Aw,w)\geq\|w\|_{\mathcal D}^2- \mu_\phi \|w\|_M^2\quad\text{for all }w,\qquad \mu_\phi \Delta t\leq1/2,
 \label{eq:transport-energy-D}
\end{equation}
which, by the proof of Lemma~\ref{lem:transport-stage}, is the case when the porosity increment of the state satisfies (A1) and its fluxes satisfy the cell balance \eqref{eq:fv-darcy-flux}.
Then the solution of $\mathcal Lw=Mv/\Delta t+g$ satisfies the estimate of that lemma in the norm $\|\cdot\|_M$.

\noindent\textit{Multiscale transport step.}
We use the multiscale space of Section~\ref{sec:ms-space}, with the coarse neighborhoods $\omega_j$ of the coarse grid $\mathcal T_H$, the snapshot matrices $X_{i,j}$ with the singular values $\sigma_{i,j,\ell}$ and the first $m$ left singular vectors $U_{i,j}^{(m)}$, the partition-of-unity weights $\chi_j$, and the multiscale functions $\varphi_{i,j,\ell}$, which are the columns of $P_i$ and span $V_i^{\mathrm{ms}}$.
For a fine-grid vector $w\in\mathbb R^{N_h}$ with restrictions $w_j$ to $\omega_j$, let $\Pi_{i,j}$ be the Euclidean orthogonal projection onto the local modes of $\omega_j$, the span of the columns of $U_{i,j}^{(m)}$ and of the constant function, and let
\begin{equation}
 \Pi_iw=\sum_j\chi_j\,\Pi_{i,j}w_j\in V_i^{\mathrm{ms}},\qquad
 \vartheta=w-\Pi_iw=\sum_j\chi_j\,\vartheta_j,\quad{\vartheta}_j=w_j-\Pi_{i,j}w_j,
 \label{eq:ms-projection}
\end{equation}
with zero extension outside $\omega_j$.
The local modes contain the columns of $U_{i,j}^{(m)}$, so $\|(I-\Pi_{i,j})X_{i,j}\|_2\leq\|X_{i,j}-U_{i,j}^{(m)}(U_{i,j}^{(m)})^TX_{i,j}\|_2=\sigma_{i,j,m+1}$, and $\vartheta_j=(I-\Pi_{i,j})X_{i,j}c+(I-\Pi_{i,j})(w_j-X_{i,j}c)$ for any coefficients $c$, with $\|I-\Pi_{i,j}\|_2\leq1$, gives
\begin{equation}
 \|\vartheta_j\|_2\leq\varepsilon_{i,j}(w)=\min_{c}\bigl(\sigma_{i,j,m+1}\|c\|_2+\|w_j-X_{i,j}c\|_2\bigr),
 \label{eq:pod-local-error}
\end{equation}
the sum of the POD truncation error, with $\sigma_{i,j,m+1}^2\leq(1-\rho_{\mathrm{tr}})\|X_{i,j}\|_F^2$ when the retained energy fraction $\rho_{\mathrm{tr}}$ is attained, and of the snapshot error, the distance of $w_j$ from the combinations of the snapshots \cite{efendiev2013gmsfem,calo2016randomized}.
For a snapshot, a column of $X_{i,j}$, $\varepsilon_{i,j}(w)\leq\sigma_{i,j,m+1}$.
The partition-of-unity weights are piecewise linear on the coarse cells, 
so $|\chi_{j,K}-\chi_{j,L}|\leq d_{KL}/H$ on every face, 
 $|e|\,d_{KL}=|K|$ on the uniform grid, 
and $\|\cdot\|_{\mathcal D,\omega_j}$ is the energy norm restricted to the faces of $\omega_j$ \cite{efendiev2013gmsfem, vasilyeva2024gmsfemgraph}.

\begin{lemma}[Projection error]
\label{lem:ms-projection}
Let $n_{\mathrm{ov}}$ be the largest number of coarse neighborhoods $\omega_j$ that contain a cell of one face, $D_{\max}=\max_KD_i(\phi_K)$, and $c_\chi=4d$.
Then, for every $w\in\mathbb R^{N_h}$,
\begin{equation}
 \|\vartheta\|_M^2\leq\sum_j\|\vartheta_j\|_M^2\leq|K|\sum_j\varepsilon_{i,j}(w)^2,\qquad
 \|\vartheta\|_{\mathcal D}^2\leq2n_{\mathrm{ov}}\sum_j\Bigl(\|\vartheta_j\|_{\mathcal D,\omega_j}^2+\frac{c_\chi D_{\max}}{H^2}\,\|\vartheta_j\|_M^2\Bigr).
 \label{eq:ms-projection-error}
\end{equation}
\end{lemma}

\begin{proof}
In a cell, 
$\vartheta_K=\sum_j\chi_{j,K}\vartheta_{j,K}$ 
with 
$0\leq\chi_j\leq1$ and $\sum_j\chi_j=1$, so $\vartheta_K^2\leq\sum_j\chi_{j,K}\vartheta_{j,K}^2\leq\sum_j\vartheta_{j,K}^2$ 
by the Cauchy--Schwarz inequality, and 
$\|\vartheta_j\|_M^2\leq|K|\|\vartheta_j\|_2^2$ 
with \eqref{eq:pod-local-error} gives the mass bounds.
On a face, at most 
$n_{\mathrm{ov}}$ terms of $\vartheta_K-\vartheta_L$ are nonzero, so 
$\|\vartheta\|_{\mathcal D}^2\leq n_{\mathrm{ov}}\sum_j\|\chi_j\vartheta_j\|_{\mathcal D}^2$.
For one $j$, the discrete product rule 
$\chi_{j,K}\vartheta_{j,K}-\chi_{j,L}\vartheta_{j,L}=\chi_{j,K}(\vartheta_{j,K}-\vartheta_{j,L})+(\chi_{j,K}-\chi_{j,L})\,\vartheta_{j,L}$, 
$(a+b)^2\leq2a^2+2b^2$, and $\chi_j\leq1$ give $2\|\vartheta_j\|_{\mathcal D,\omega_j}^2$ from the first parts, which include the faces on $\Gamma_{i,D}$, while on a face between 
$K\subset\omega_j$ and $L\not\subset\omega_j$ only the second part, with $\vartheta_{j,K}$, is present.
For the second parts, the harmonic mean in 
$T^D_{i,KL}=|e|(\eta_iD_i)_{KL}/d_{KL}$ satisfies $(\eta_iD_i)_{KL}\leq2D_{\max}\min(\eta_i(\phi_K),\eta_i(\phi_L))$, so $T^D_{i,KL}(\chi_{j,K}-\chi_{j,L})^2\leq2D_{\max}H^{-2}\,|K_e|\,\eta_i(\phi_{K_e})$ for the cell $K_e$ of the face at which $\vartheta_j$ is taken, and each cell is $K_e$ for at most $2d$ faces.
Hence $\|\chi_j\vartheta_j\|_{\mathcal D}^2\leq2\|\vartheta_j\|_{\mathcal D,\omega_j}^2+2c_\chi D_{\max}H^{-2}\|\vartheta_j\|_M^2$, and summing over $j$ gives the energy bound.
\end{proof}

\begin{lemma}[Multiscale transport step]
\label{lem:ms-step}
Let \eqref{eq:transport-energy-D} hold, let $P_i$ have full column rank, and let $w$ be the solution of the fine-grid system \eqref{eq:training-system} and $w_{\mathrm{ms}}=P_ia_i^*$ the solution of the Galerkin system \eqref{eq:reduced-local} with the same right-hand side $f=Mv/\Delta t+g$, where $v=C_i^n$ and $g=b_i(\kappa)$, so that $\mathcal Lw=f$ and $P_i^T\mathcal LP_i\,a_i^*=P_i^Tf$.
Then the Galerkin system is nonsingular, and
\begin{equation}
 \|w_{\mathrm{ms}}\|_M\leq(1+2 \mu_\phi \Delta t)\bigl(\|v\|_M+\Delta t\,\|M^{-1/2}g\|_2\bigr),\qquad
 \|w-w_{\mathrm{ms}}\|_M\leq c_\nu\bigl(\|\vartheta\|_M^2+\Delta t\,\|\vartheta\|_{\mathcal D}^2\bigr)^{1/2},
 \label{eq:ms-step-error}
\end{equation}
with $\vartheta=w-\Pi_iw$, $c_\nu=3+2\nu$, and 
$\nu=\Delta t\,\|M^{-1/2}\mathcal CM^{-1/2}\|_2$. 
\end{lemma}

\begin{proof}
The Galerkin system is $(\mathcal Lw_{\mathrm{ms}},y)=(f,y)$ for all $y\in V_i^{\mathrm{ms}}$, and $y=w_{\mathrm{ms}}$ gives the energy inequality of the proof of Lemma~\ref{lem:transport-stage} with $w_{\mathrm{ms}}$ in place of $w$, hence the stability estimate and, for $f=0$, $w_{\mathrm{ms}}=0$.
Subtracting the two systems gives $(\mathcal L(w-w_{\mathrm{ms}}),y)=0$ for all $y\in V_i^{\mathrm{ms}}$.
We split $w-w_{\mathrm{ms}}=\vartheta+\zeta$ with $\zeta=\Pi_iw-w_{\mathrm{ms}}\in V_i^{\mathrm{ms}}$ and take $y=\zeta$, so that $(\mathcal L\zeta,\zeta)=-(\mathcal L\vartheta,\zeta)$, and \eqref{eq:transport-energy-D} for the left-hand side and the Cauchy--Schwarz inequality in the mass and energy inner products for the right-hand side give
\begin{equation*}
 \frac{1-\mu_\phi\Delta t}{\Delta t}\|\zeta\|_M^2+\|\zeta\|_{\mathcal D}^2
 \leq\frac{1+\nu}{\Delta t}\|\vartheta\|_M\|\zeta\|_M+\|\vartheta\|_{\mathcal D}\|\zeta\|_{\mathcal D}.
\end{equation*}
By Young's inequality, the right-hand side is at most $\Delta t^{-1}(1+\nu)^2\|\vartheta\|_M^2+(4\Delta t)^{-1}\|\zeta\|_M^2+\frac12\|\vartheta\|_{\mathcal D}^2+\frac12\|\zeta\|_{\mathcal D}^2$, and since $\frac34-\mu_\phi\Delta t\geq\frac14$, multiplying by $4\Delta t$ gives $\|\zeta\|_M^2\leq4(1+\nu)^2\bigl(\|\vartheta\|_M^2+\Delta t\,\|\vartheta\|_{\mathcal D}^2\bigr)$.
The triangle inequality $\|w-w_{\mathrm{ms}}\|_M\leq\|\vartheta\|_M+\|\zeta\|_M$ gives the error estimate.
\end{proof}

The bounds \eqref{eq:ms-projection-error} and \eqref{eq:ms-step-error} of Lemmas~\ref{lem:ms-projection} and~\ref{lem:ms-step}, applied to each species with $w=\widetilde C_{q,i}^{*,n}$, give the multiscale error of Theorem~\ref{thm:learned-errors}.
We measure the error in the $L^2$ norm, which bounds the phase norms by \eqref{eq:norm-equivalence}, and
\begin{equation}
 \|\widehat C_{\mathrm{ms}}^{*,n}-\widetilde C_q^{*,n}\|
 \leq\eta_{\min}^{-1/2}c_\nu\Bigl(\sum_{i,j}\bigl(c_H |K|\,\varepsilon_{i,j}(\widetilde C_{q,i}^{*,n})^2+2n_{\mathrm{ov}}\Delta t\,\|\vartheta_{i,j}\|_{\mathcal D_i,\omega_j}^2\bigr)\Bigr)^{1/2}=\delta_{\mathrm{ms}}^n,
 \label{eq:delta-ms}
\end{equation}
with the local residuals $\vartheta_{i,j}$ of \eqref{eq:ms-projection} and the energy norm $\|\cdot\|_{\mathcal D_i,\omega_j}$ for species $i$ and $w=\widetilde C_{q,i}^{*,n}$, $c_H=1+2n_{\mathrm{ov}}c_\chi D_{\max}\Delta t/H^2$, and the largest $c_\nu$ and $D_{\max}$ over the species.
The snapshot error depends on the test medium and admits no a priori bound, but it can be evaluated a posteriori from a fine-grid solution.
The local error $\varepsilon_{i,j}$ bounds the values of the residual $\vartheta_{i,j}$, and the local energy $\|\vartheta_{i,j}\|_{\mathcal D_i,\omega_j}^2$ measures its discrete gradient, which is the same POD truncation and snapshot error in the energy norm.
A local spectral problem on the snapshots with a reference diffusion matrix, for example for the averaged porosity, would bound its truncation part by the first discarded eigenvalue, as in spectral multiscale spaces \cite{efendiev2013gmsfem,vasilyeva2024gmsfemgraph}, up to the change of $\eta_iD_i(\phi)$.

\noindent\textit{Learned transport step.}
For one species $i$, in the notation of Lemma~\ref{lem:ms-step} with $b=b_i(\kappa)$, the coefficients of the multiscale step $w_{\mathrm{ms}}=P_ia$ and of the learned step $\widehat w=P_i\widehat a$ of the Galerkin local atoms solve
\begin{equation*}
 \mathsf La=\mathsf f,\quad (\mathsf L+\mathsf E_A)\,\widehat a=\mathsf f+\mathsf E_b,\quad 
 \mathsf L=P_i^T\mathcal LP_i,\quad
 \mathsf f=P_i^T\Bigl(\frac{Mv}{\Delta t}+b\Bigr),
\end{equation*}
since the storage terms of \eqref{eq:local-galerkin} are exact and the learning changes only the reduced transport matrix and boundary vector, by the truncation and regression errors $\mathsf E_A$ and $\mathsf E_b$ of the predicted pieces $X_{i,c}$ and $y_{i,c}$ of \eqref{eq:local-atoms}.
With $\mathsf M=P_i^TMP_i$, let $\varepsilon_A=\|\mathsf M^{-1/2}\mathsf E_A \mathsf M^{-1/2}\|_2$ and $\varepsilon_b=\|\mathsf M^{-1/2} \mathsf E_b\|_2$ be these errors relative to the storage, so that $\varepsilon_A$ is a rate that adds to $\mu_\phi$ in the stability factor.

\begin{lemma}[Learned transport step]
\label{lem:learned-step}
Let the assumptions of Lemma~\ref{lem:ms-step} hold, and let $(\mu_\phi+\varepsilon_A)\Delta t\leq1/2$.
Then the learned system is nonsingular, and
\begin{equation}
\begin{aligned}
 \|\widehat w\|_M&\leq\bigl(1+2(\mu_\phi+\varepsilon_A)\Delta t\bigr)\bigl(\|v\|_M+\Delta t\,\|\mathsf M^{-1/2}(P_i^Tb+ \mathsf E_b)\|_2\bigr),\\
 \|\widehat w-w_{\mathrm{ms}}\|_M&\leq(1+2\mu_\phi\Delta t)\,\Delta t\,\bigl(\varepsilon_A\|\widehat w\|_M+\varepsilon_b\bigr).
\end{aligned}
 \label{eq:learned-step-error}
\end{equation}
\end{lemma}

\begin{proof}
For 
$y\in\mathbb R^{r_i}$ and $z=P_iy$, 
we have 
$\|\mathsf M^{1/2}y\|_2=\|z\|_M$, $(\mathsf Ly,y)=(\mathcal Lz,z)\geq(\Delta t^{-1}-\mu_\phi)\|z\|_M^2$ 
by \eqref{eq:transport-energy-D}, and 
$|(\mathsf E_A y,y)|\leq\varepsilon_A\|z\|_M^2$.
Hence the inner product of the learned system with $\widehat a$ gives
\begin{equation*}
 \frac{1-(\mu_\phi+\varepsilon_A)\Delta t}{\Delta t}\|\widehat w\|_M^2
 \leq\frac{1}{\Delta t}(Mv,\widehat w)+(P_i^Tb+ \mathsf E_b,\widehat a)
 \leq\Bigl(\frac{1}{\Delta t}\|v\|_M+\|\mathsf M^{-1/2}(P_i^Tb + \mathsf E_b)\|_2\Bigr)\|\widehat w\|_M,
\end{equation*}
which gives the first estimate, as in the proof of Lemma~\ref{lem:transport-stage}. 
For the second, the learned solution satisfies 
$\mathsf L(\widehat a-a)=\varrho$ with the residual $\varrho = \mathsf E_b - \mathsf E_A \widehat a$, 
and the inner product with 
$y=\widehat a-a$ gives $(\Delta t^{-1}-\mu_\phi)\|z\|_M^2\leq(\varrho,y)\leq\|\mathsf M^{-1/2}\varrho\|_2\,\|z\|_M$ 
for 
$z=\widehat w-w_{\mathrm{ms}}$, with $\|\mathsf M^{-1/2}\varrho\|_2\leq\varepsilon_b+\varepsilon_A\|\widehat w\|_M$.
\end{proof}

The first estimate is the stability of the learned transport step, in which the learning adds its operator error $\varepsilon_A$ to the rate $\mu_\phi$ of the fine-grid step. The second estimate for each species with \eqref{eq:norm-equivalence} gives the learning error of Theorem~\ref{thm:learned-errors},
\begin{equation}
 \|\widehat C^{*,n}-\widehat C_{\mathrm{ms}}^{*,n}\|
 \leq\eta_{\min}^{-1/2}(1+2\mu_\phi\Delta t)\,\Delta t\,\Bigl(\sum_i\bigl(\varepsilon_{A,i}^n\|P_i\widehat a_i^{*,n}\|_{M_i}+\varepsilon_{b,i}^n\bigr)^2\Bigr)^{1/2}=\delta_\ell^n.
 \label{eq:delta-ell}
\end{equation}

\noindent\textit{Operator atoms.}
The errors $\mathsf E_A$ and $\mathsf E_b$ come from a POD of the vectorized pieces followed by a regression of the coefficients.
Let $v=v_{i,c}(\kappa)$ be the piece of \eqref{eq:local-atoms}, which contains $X_{i,c}$ and $y_{i,c}$, let $\bar v$ be the training mean, and let $\mathcal Y$ be the matrix whose columns are the centered training vectors, assembled for the initial media of the training realizations as in Section~\ref{sec:ms-space}, with the singular values $\sigma^{\mathrm a}_k$.
The atoms $\psi_k=\psi_{i,c,k}$, $k=1,\ldots,K_{\mathrm a}$, are the first $K_{\mathrm a}$ left singular vectors of $\mathcal Y$, $\Pi_{\mathrm a}$ is the orthogonal projection onto their span, and the prediction is $\widehat v=\bar v+\sum_{k=1}^{K_{\mathrm a}}\widehat w_k\psi_k$ with the regressed coefficients $\widehat w_k=\widehat w_{i,c,k}$, where $w_k=(\psi_k,v-\bar v)$ are the exact coefficients.

\begin{lemma}[Error of the operator atoms]
\label{lem:atoms}
For every medium $\kappa$,
\begin{equation}
 \|v-\widehat v\|_2^2=\|(I-\Pi_{\mathrm a})(v-\bar v)\|_2^2+\|w-\widehat w\|_2^2,\quad 
 \|(I-\Pi_{\mathrm a})(v-\bar v)\|_2\leq\min_{\tau}\bigl(\sigma^{\mathrm a}_{K_{\mathrm a}+1}\|\tau\|_2+\|v-\bar v-\mathcal Y\tau\|_2\bigr)
 \label{eq:atom-error}
\end{equation}
with the minimum over the coefficients $\tau$ of the columns of $\mathcal Y$.
\end{lemma}

\begin{proof}
Since $v-\widehat v=(I-\Pi_{\mathrm a})(v-\bar v)+\sum_k(w_k-\widehat w_k)\psi_k$ and the atoms are orthonormal, the two terms are orthogonal, which gives the identity.
The bound follows as \eqref{eq:pod-local-error}, from $\|(I-\Pi_{\mathrm a})\mathcal Y\|_2=\sigma^{\mathrm a}_{K_{\mathrm a}+1}$ and $\|I-\Pi_{\mathrm a}\|_2\leq1$.
\end{proof}

The first part is the atom truncation error, which has the same two terms as the local error \eqref{eq:pod-local-error} of the multiscale space, the truncation of the training operators and the distance of the operator of the current medium from their combinations, which includes the change of the medium since the initial time, and the second part is the regression error of the atom coefficients.
For the Galerkin local atoms, $\mathsf E_A=\sum_cS_{i,c}(\widehat X_{i,c}-X_{i,c})S_{i,c}^T$ and $\mathsf E_b=\sum_cS_{i,c}(\widehat y_{i,c}-y_{i,c})$ by \eqref{eq:local-galerkin}.
Let $\varepsilon^{\mathrm a}_{i,c}=\|\widehat v_{i,c}-v_{i,c}\|_2$ be the error of Lemma~\ref{lem:atoms}, which bounds the errors of both pieces, and $\lambda_{\min}(\mathsf M)$ the smallest eigenvalue of $\mathsf M$.
A multiscale function is nonzero in several coarse cells, so its index belongs to several sets $\mathcal J_{i,c}$, at most $n_{\mathcal J}$ of them.
Hence $\sum_c\|S_{i,c}^Ta\|_2^2\leq n_{\mathcal J}\|a\|_2^2$ for $a\in\mathbb R^{r_i}$, and
\begin{equation*}
 \varepsilon_A\leq\frac{n_{\mathcal J}}{\lambda_{\min}(\mathsf M)}\,\max_c\varepsilon^{\mathrm a}_{i,c},\quad 
 \varepsilon_b\leq\Bigl(\frac{n_{\mathcal J}}{\lambda_{\min}(\mathsf M)}\sum_c(\varepsilon^{\mathrm a}_{i,c})^2\Bigr)^{1/2}.
\end{equation*}
The atoms are optimal in the Euclidean norm of the entries of the pieces, while $\varepsilon_A$ and $\varepsilon_b$ are weighted by the storage.

{
\noindent\textit{Remark (LSPG).}
For the LSPG system \eqref{eq:reduced-local}, $w_{\mathrm{ms}}$ minimizes $\|\mathcal L(w-P_ia)\|_2$, so Lemma~\ref{lem:transport-stage} with $v=0$ and \eqref{eq:ms-projection-error} give 
$\|w-w_{\mathrm{ms}}\|_M\leq c_{\mathcal L}\|\vartheta\|_M$ and 
\begin{equation*}
\delta_{\mathrm{ms}}^n = \eta_{\min}^{-1/2} c_{\mathcal L} \Bigl(|K|\sum_{i,j}\varepsilon_{i,j}(\widetilde C_{q,i}^{*,n})^2\Bigr)^{1/2},
\end{equation*} 
with $c_{\mathcal L}=(1+2\mu_\phi\Delta t)\,\eta_{\min}^{-1/2}(1+\nu+\nu_{\mathcal D})$, the diffusion number $\nu_{\mathcal D}=\Delta t\,\|M^{-1/2}\mathcal DM^{-1/2}\|_2$ of the fine grid, and without the energy term.
Lemma~\ref{lem:atoms} holds also for the LSPG pieces $v_{i,c}$, which contain $B_{i,c}$ and $\beta_{i,c}$.
The learned step \eqref{eq:local-lspg} minimizes $\|\widehat\Gamma a-\widehat f\|_2$ with $\widehat\Gamma=\mathcal LP_i+\mathsf E_B$ and $\widehat f=f+\mathsf E_\beta$, where $\mathsf E_B$ and $\mathsf E_\beta$ stack the errors of the predicted pieces over the coarse cells, so that $\|\mathsf E_B\|_2^2+\|\mathsf E_\beta\|_2^2\leq\sum_c(\varepsilon^{\mathrm a}_{i,c})^2$.
The normal equations of the multiscale and learned steps and Lemma~\ref{lem:transport-stage} with $v=0$ give
\begin{equation*}
 \|\widehat w-w_{\mathrm{ms}}\|_M\leq c_M\,\Delta t\,|K|^{-1/2}\bigl(\|\mathsf E_\beta\|_2+\|\mathsf E_B\|_2\|\widehat a\|_2+\sigma_\Gamma^{-1}\|\mathsf E_B\|_2\|\widehat r\|_2\bigr),
\end{equation*}
with $c_M=(1+2\mu_\phi\Delta t)\,\eta_{\min}^{-1/2}$, the smallest singular value $\sigma_\Gamma$ of $\mathcal LP_i$, which satisfies $\sigma_\Gamma^{-1}\leq c_M\Delta t\,(|K|\lambda_{\min}(\mathsf M))^{-1/2}$, and the residual $\widehat r=\widehat f-\widehat\Gamma\widehat a$ of the learned step.
The sum over the species, as in \eqref{eq:delta-ell}, gives $\delta_\ell^n$, which is proportional to $\Delta t$ and has one more term, the operator error times the residual.
}

{
\noindent\textit{Remark (global inverse atoms).}
The global inverse atoms are trained on the encoded fine-grid step \eqref{eq:reduced-global}, so for them the multiscale step in Theorem~\ref{thm:learned-errors} is $\widehat C_{\mathrm{ms},i}^{*,n}=\max\{P_i(G_i(\kappa)a_i^n+g_i(\kappa)),0\}$ with $a_i^n=E_i\widehat C_i^n$, and the learned step replaces $(G_i,g_i)$ by $(I+\widehat J_i,\widehat g_i)$.
The coefficients of the two steps differ by $(\widehat J_i-J_i)a_i^n+\widehat g_i-g_i$, without a reduced solve, so that, since the cutoff does not increase the difference and $\|u\|=|K|^{1/2}\|u\|_2$,
\begin{equation}
 \|\widehat C^{*,n}-\widehat C_{\mathrm{ms}}^{*,n}\|
 \leq|K|^{1/2}\Bigl(\sum_i\|P_i\|_2^2\bigl(\|\widehat J_i-J_i\|_F\|a_i^n\|_2+\|\widehat g_i-g_i\|_2\bigr)^2\Bigr)^{1/2}=\delta_\ell^n,
 \label{eq:delta-ell-global}
\end{equation}
where $\|\widehat J_i-J_i\|_F$ is the error of Lemma~\ref{lem:atoms}, which holds also for the vectorized increment $J_i$ with $\bar v=0$, the $K_i$ vectorized atoms $\Theta_{i,k}$, and the coefficients $\widehat\alpha_{i,k}$, and $\|\widehat g_i-g_i\|_2$ is the regression error of the offset.
The bound has no explicit factor $\Delta t$, but $J_i=G_i-I$ is itself the change of the coefficients over one step.
For the multiscale error, $P_iE_i$ is the Euclidean orthogonal projection onto $V_i^{\mathrm{ms}}$, up to the small shift of the encoder, so $\|(I-P_iE_i)w\|_2\leq\|w-\Pi_iw\|_2\leq(\sum_j\varepsilon_{i,j}(w)^2)^{1/2}$ as in the proof of Lemma~\ref{lem:ms-projection}.
The encoded step differs from the fine-grid step $\widetilde C_{q,i}^{*,n}$ by this projection error and by the encoding error of the current state, carried over one step by the estimate of Lemma~\ref{lem:transport-stage}, so that, when $\widetilde C_q^{*,n}\geq0$,
\begin{equation}
 \|\widehat C_{\mathrm{ms}}^{*,n}-\widetilde C_q^{*,n}\|
 \leq|K|^{1/2}\Bigl(
 \sum_i\Bigl 
 (\bigl(\sum_j\varepsilon_{i,j}(\widetilde C_{q,i}^{*,n})^2\bigr)^{1/2} + 
 c_M\|(I-P_iE_i)\widehat C_i^n\|_2\Bigr)^2
 \Bigr)^{1/2} 
 = \delta_{\mathrm{ms}}^n,
 \label{eq:delta-ms-global}
\end{equation}
with $c_M$ as above, and the energy norm does not enter.
Before the cutoff, the learned step is an affine map of $C_i^n$, as is the step \eqref{eq:local-lspg} of the LSPG local atoms.
The stability of these maps, a bound $1+\widehat\Lambda\Delta t$ for the norm of the linear part induced by $\|\cdot\|_M$, is not guaranteed but can be evaluated with $r_i\times r_i$ matrices along the computed trajectory.
}

\noindent\textit{Multiscale flow step.}
The learned flow of Section~\ref{sec:ml-flow} is the mixed multiscale method on the flux space $\operatorname{range}\Psi$, and its error follows the analysis of mixed multiscale methods \cite{chung2015mixed,chung2016mixed,spiridonov2022mixed}.
Two points differ from these methods.
The modes are the extensions of predicted interface traces and not the eigenvectors of a local spectral problem, and they are built with the initial permeability and kept while the permeability changes.
For a state $U$, let $q=\mathcal F(U)$ and $q_{\mathrm{ms}}=\Psi a_q=\widehat{\mathcal F}(U)$ be the fine-grid fluxes and the multiscale fluxes of \eqref{eq:flow-reduced-system}, $T=T(k_U)$, $\|v\|_{T^{-1}}^2=(T^{-1}v,v)$, $B_h$ the cell divergence, and $R_H$ the sum of the cell values over each coarse block $B$, so that $\mathsf B=R_HB_h\Psi$ and $S=R_Hs_U$.
As in mixed multiscale methods, we assume that each mode has a divergence proportional to the cell volume in each block, $B_h\Psi=Z\mathsf B$ with $Z_{KB}=|K|/|B|$ for $K\subset B$ and the block volume $|B|$, and that every coarse interface has a mode with a nonzero net flux \cite{chung2016mixed,spiridonov2022mixed}, so that \eqref{eq:flow-reduced-system} has a unique solution.
Let $q_r$ be the flux with the smallest norm $\|\cdot\|_{T^{-1}}$ that is supported inside the blocks, has zero normal flux on their boundaries, and satisfies $B_hq_r=s_U-ZS$, and let $\beta=\lambda_{\min}(\mathsf BM_q^{-1}\mathsf B^T)^{1/2}$ be the inf-sup constant of \eqref{eq:flow-reduced-system}.
Then, for every state $U$ and every $b_0\in\mathbb R^{m_q}$,
\begin{equation}
 \|q-q_{\mathrm{ms}}\|_{T^{-1}}
 \leq\|q-q_r-\Psi b_0\|_{T^{-1}}+\beta^{-1}\|S-\mathsf Bb_0\|_2+\|q_r\|_{T^{-1}},
 \label{eq:flow-error}
\end{equation}
and $\|\widehat{\mathcal F}(U)-\mathcal F(U)\|\leq T_{\max}^{1/2}\|q-q_{\mathrm{ms}}\|_{T^{-1}}$, since the transmissibilities are at most $T_{\max}$, which gives $\varepsilon_f^n$ for $U=\widehat U^n$.
This estimate holds for any flux space that satisfies the assumptions above and follows as in the analysis of mixed multiscale methods \cite{chung2016mixed,spiridonov2022mixed}, so we do not repeat the proof. 
The two differences enter only the first two terms, the approximation error of the flux space and of its block divergence.
In \cite{chung2016mixed,spiridonov2022mixed}, these terms are bounded by the eigenvalues of the local spectral problem, and here they consist of the trace POD truncation to $K_f$ modes, the prediction error of the networks $N_{f,\iota}$, and the change of the permeability since the initialization.
The third term, the oscillation of the source $s_U-ZS$ inside the blocks, is bounded as in \cite{spiridonov2022mixed}.
The cell divergence of $q_{\mathrm{ms}}$ is $ZS$, the block average of the porosity rate, so \eqref{eq:transport-energy-D} holds for the transport with the multiscale fluxes when the block-averaged porosity increment satisfies (A1).

\section{Training data and settings}
\label{app:training}

We list here the training data, model settings, and offline costs of the learned stages used in Section~\ref{sec:results}.
Transport and flow use the coarse grid of ten cells along the depth described there.

\noindent\textbf{Chemistry.} We collect 350\,000 local concentration pairs from 15 one-dimensional and 20 two-dimensional reference trajectories, sampling 100 cells per time step.
Two networks with three hidden layers of width 64 are trained jointly for 300 epochs using Adam and a weighted concentration loss.
The offline cost is 516.7 s.
The same model is used in 2D and 3D, without 3D training data.

\noindent\textbf{Transport.} 
\textit{Multiscale representation.}
In each dimension, 40 realizations define a multiscale basis on 11 depth neighborhoods, with at most 20 data-driven modes and a constant mode per neighborhood ($\rho_{\mathrm{tr}}=1-\epsilon$, $\epsilon=10^{-8}$, with the retained-energy criterion defined in Subsection~\ref{sec:ms-space}).
The maximum reduced dimension per species is $11\cdot(20+1)=231$, including the constant modes.
The reduced dimensions are $(229,78,197,231,231)$ in 2D and $(231,79,228,231,231)$ in 3D, giving 966 and 1,000 transport degrees of freedom in total, respectively, compared with 5,000 and 50,000 fine-grid concentration degrees of freedom.
\textit{Operator data and learning.}
Initial-medium operators from 5,000 random fields provide 4,000 fitting and 1,000 validation samples for ridge regression.
We fix $K_i=64/256$ in 2D and $K_i=128/512$ in 3D for solid/dissolved species.
%Detailed tests with varying atom counts, including errors and timings, are presented {in Sections~S1 and~S2 of the Supplementary Material}. 
\textit{Offline cost.}
Offline costs are approximately 304 s in 2D and 2,157 s in 3D, including field generation, snapshot trajectories, basis construction, operator compression, and ridge fitting.
Operator compression accounts for 216 and 681 s, and regression for 17.7 and 23.3 s, respectively.

\noindent\textbf{Flow.}
\textit{Multiscale representation.}
The reconstructed directional traces produce two spatial flux modes per interface in 2D and three in 3D.
Across 11 interfaces, including the pressure boundaries, these give 22 and 33 flux coefficients.
Together with ten coarse-block pressures, the reduced flow system has 32 and 43 unknowns, respectively.
The basis remains fixed while transmissibilities and porosity sources are updated each step.
The realization-specific flux basis is constructed online by local extensions and is included in the total prediction time.
\textit{Trace data and learning.}
Per dimension, 5,000 permeability fields provide 4,000 training, 500 validation, and 500 trace-evaluation samples.
Interface predictors use an MLP with two hidden layers of width 512 in 2D and a CNN with three convolutional layers in 3D.
AdamW trains the 2D MLP for 30 epochs and the 3D CNN for 120 epochs.
The CNN uses the eight symmetries of the square interface to augment the permeability patches and their associated directional trace targets.
Models are selected by the validation error in the predicted trace coefficients.
For each interface, the network takes a local permeability patch as input and outputs coefficients of fixed trace basis vectors learned from the training data.
The count $K_f$ is the number of these vectors and hence the network output dimension.
Their weighted sum reconstructs the fine-face flux traces.
We choose a compression limit of 64 trace basis vectors.
In 2D, each coarse interface contains ten fine faces, with traces from two pressure-gradient directions, giving a flux-trace target vector of $10\times2=20$ entries.
The network predicts 20 coefficients to reconstruct this vector.
In 3D, 100 fine faces and three directions give a flux-trace vector of $100\times3=300$ entries, reconstructed from 64 predicted coefficients.
\textit{Offline cost.}
Training the selected trace predictors takes approximately 11 s for the 2D MLP and 366 s for the augmented 3D CNN, summed over the three interface types.
The full offline costs are approximately 36.5 s in 2D and 701 s in 3D, including permeability-field generation, reference flux-trace generation, trace preprocessing, and training of the selected predictors.

\end{document}